\documentstyle[amscd, leqno]{amsart}
\theoremstyle{plain}
\newtheorem{theorem}{Theorem}[section]
\newtheorem{proposition}[theorem]{Proposition}
\newtheorem{lemma}[theorem]{Lemma}

\newtheorem{corollary}[theorem]{Corollary}
\theoremstyle{definition}

\theoremstyle{remark}
\newtheorem{remark}[theorem]{Remark}
\newtheorem{example}[theorem]{Example}
\newtheorem{question}[theorem]{Question}

\makeatletter

\@addtoreset{equation}{section}
\makeatother

\begin{document}
\title[Singularities and analytic torsion]
{Singularities and analytic torsion}
\author{Ken-Ichi Yoshikawa}
\address{
Department of Mathematics,
Faculty of Science,
Kyoto University,
Kyoto 606-8502, JAPAN}
\email{yosikawa@@math.kyoto-u.ac.jp}
\address{Korea Institute for Advanced Study,
Hoegiro 87, Dongdaemun-gu,
Seoul 130-722, KOREA}

\thanks{The author is partially supported by the Grants-in-Aid 
for Scientific Research (B) 19340016, JSPS}

\begin{abstract}
We prove the logarithmic divergence of equivariant analytic torsion for one-parameter 
degenerations of projective algebraic manifolds, when the coefficient vector bundle is given by 
a Nakano semi-positive vector bundle twisted by the relative canonical bundle.
\end{abstract}

\maketitle
\tableofcontents


\section
{Introduction}
\label{sect:1}
\par
Let ${\mathcal X}$ be a connected projective algebraic manifold of dimension $n+1$ 
and let $C$ be a compact Riemann surface.
Let $\pi\colon{\mathcal X}\to C$ be a surjective holomorphic map with critical locus $\Sigma_{\pi}$ 
and with connected fibers. 
\par
Let $G$ be a compact Lie group acting holomorphically on ${\mathcal X}$ and preserving 
the fibers of $\pi$. Assume that there exists a $G$-equivariant ample line bundle on ${\mathcal X}$.
Set $\Delta=\pi(\Sigma_{\pi})$, $C^{o}=C\setminus\Delta$,
${\mathcal X}^{o}={\mathcal X}|_{\pi^{-1}(C^{o})}$, $\pi^{o}=\pi|_{{\mathcal X}^{o}}$ and
$X_{s}:=\pi^{-1}(s)$ for $s\in C$.
Then $\pi^{o}\colon{\mathcal X}^{o}\to C^{o}$ is a family of projective algebraic manifolds with $G$-action.
\par
Let $T{\mathcal X}/C$ be the $G$-equivariant subbundle of 
$T{\mathcal X}|_{{\mathcal X}\setminus\Sigma_{\pi}}$ defined as
$T{\mathcal X}/C=\ker\pi_{*}|_{{\mathcal X}\setminus\Sigma_{\pi}}$.
Let $h_{\mathcal X}$ be a $G$-invariant K\"ahler metric on ${\mathcal X}$ and 
set $h_{{\mathcal X}/C}=h_{\mathcal X}|_{T{\mathcal X}/C}$.
Let $\omega_{\mathcal X}=\Omega_{\mathcal X}^{n+1}$ be the canonical bundle of ${\mathcal X}$ 
and let $\omega_{{\mathcal X}/C}=\Omega_{\mathcal X}^{n+1}\otimes(\pi^{*}\Omega^{1}_{C})^{-1}$
be the relative canonical bundle of $\pi\colon{\mathcal X}\to C$. 
Let $\xi\to{\mathcal X}$ be a $G$-equivariant holomorphic vector bundle on ${\mathcal X}$ 
equipped with a $G$-invariant Hermitian metric $h_{\xi}$.
We write $\omega_{{\mathcal X}/C}(\xi)=\omega_{{\mathcal X}/C}\otimes\xi$.
We set $\xi_{s}=\xi|_{X_{s}}$ for $s\in C$.
\par
Let $0\in\Delta$ be a critical value of $\pi$. Let $(S,s)$ be a coordinate neighborhood of $C$ 
centered at $0$ such that $S \cap\Delta=\{0\}$. We set  $X=\pi^{-1}(S)$ and $S^{o}=S\setminus\{0\}$.
\par
For $g\in G$ and $s\in S^{o}$, let $\tau_{G}(X_{s},\omega_{X_{s}}(\xi_{s}))(g)$ be the 
equivariant analytic torsion \cite{Bismut95} of $(X_{s},\omega_{X_{s}}(\xi_{s}))$ with respect to 
$h_{X_{s}}=h_{\mathcal X}|_{X_{s}}$ and $h_{\xi_{s}}=h_{\xi}|_{X_{s}}$, where 
$\omega_{X_{s}}=\Omega_{X_{s}}^{n}$ for $s\not=0$. 
The goal of this article is to determine the behavior of 
$\tau_{G}(X_{s},\omega_{X_{s}}(\xi_{s}))(g)$ as $s\to0$, when $(\xi,h_{\xi})$ is 
{\em Nakano semi-positive} on $X$.
(See Sect.\ref{subsect:5.1}  for the notion of Nakano semi-positivity.)
Notice that all $R^{q}\pi_{*}\omega_{X/S}(\xi)$ are locally free in this case
by Takegoshi's torsion freeness theorem \cite{Takegoshi95} and the condition $\dim S=1$.
To express the singularity of $\tau_{G}(X_{s},\omega_{X_{s}}(\xi_{s}))(g)$ as $s\to0$ in more detail,
we briefly recall Gauss maps and semistable reductions.
\par
Let ${\bf P}(T{\mathcal X})^{\lor}$ be the $G$-equivariant projective-space bundle 
such that ${\bf P}(T{\mathcal X})^{\lor}_{x}={\bf P}(T_{x}{\mathcal X})^{\lor}$ 
is the set of $n$-dimensional linear subspaces of $T_{x}{\mathcal X}$ for $x\in{\mathcal X}$. 
The Gauss map 
$\gamma\colon{\mathcal X}\setminus\Sigma_{\pi}\to{\bf P}(T{\mathcal X})^{\lor}$ 
is the section defined as
$\gamma(x)=\ker(\pi_{*})_{x}\in{\bf P}(T_{x}{\mathcal X})^{\lor}$
for $x\in{\mathcal X}\setminus\Sigma_{\pi}$.
Since $\gamma$ extends to a rational map
$\gamma\colon{\mathcal X}\dashrightarrow{\bf P}(T{\mathcal X})^{\lor}$,
there is a resolution $q\colon(\widetilde{\mathcal X},E)\to({\mathcal X},\Sigma_{\pi})$
of the indeterminacy of $\gamma$ with 
$q|_{\widetilde{\mathcal X}\setminus E}\colon\widetilde{\mathcal X}\setminus E
\cong{\mathcal X}\setminus\Sigma_{\pi}$ 
such that $\widetilde{\gamma}=\gamma\circ q$ extends to a holomorphic map 
from $\widetilde{\mathcal X}$ to ${\bf P}(T{\mathcal X})^{\lor}$
and such that $E$ is a normal crossing divisor of $\widetilde{\mathcal X}$. 
Since $\gamma$ is $G$-equivariant, we may assume 
that $G$ acts on $\widetilde{\mathcal X}$ and that $q$ and $\widetilde{\gamma}$ are $G$-equivariant
\cite{BierstoneMilman97}.
We denote by ${\mathcal H}={\mathcal O}_{{\bf P}(T{\mathcal X})^{\lor}}(1)$
the tautological quotient bundle on ${\bf P}(T{\mathcal X})^{\lor}$.
\par
For $g\in G$, let ${\mathcal X}^{g}=\{x\in{\mathcal X};\,g\cdot x=x\}$ be its fixed-point set.
Since $g$ is an isometry of ${\mathcal X}$,
${\mathcal X}^{g}$ is the disjoint union of compact complex submanifolds of ${\mathcal X}$:
$$
{\mathcal X}^{g}={\mathcal X}^{g}_{H}\amalg{\mathcal X}^{g}_{V}.
$$
Here ${\mathcal X}^{g}_{H}$ is a horizontal submanifold, i.e.,
$\pi|_{{\mathcal X}^{g}_{H}}\colon{\mathcal X}^{g}_{H}\to C$ is a flat holomorphic map and 
${\mathcal X}^{g}_{V}$ is a vertical submanifold, i.e., $\pi({\mathcal X}_{V}^{g})$ is a proper subset of $C$. 
Since the $G$-action on $C$ is trivial, one has 
${\mathcal X}^{g}_{V}\subset\Sigma_{\pi}$ and $\pi({\mathcal X}^{g}_{V})\subset\Delta$
by the $G$-equivariance of $\pi$.
Let $\widetilde{\mathcal X}^{g}_{H}\subset\widetilde{\mathcal X}$ be the proper transform 
of ${\mathcal X}^{g}_{H}\subset{\mathcal X}$.
Since $\widetilde{\mathcal X}^{g}_{H}\subset(\widetilde{\mathcal X})^{g}$, we get 
$\widetilde{\gamma}(\widetilde{\mathcal X}^{g}_{H})\subset({\bf P}(T{\mathcal X})^{\lor})^{g}$ 
by the $G$-equivariance of $\widetilde{\gamma}$.
Hence $g\in G$ preserves the fibers of
$(\widetilde{\gamma}^{*}{\mathcal U})|_{\widetilde{\mathcal X}^{g}_{H}}$. 
We set $E_{0}=(\pi\circ q)^{-1}(0)\cap E$ and define 
$$
\begin{aligned}
\alpha_{g}(X_{0},\omega_{X/S}(\xi))
&=
\int_{E_{0}\cap\widetilde{\mathcal X}^{g}_{H}}
\widetilde{\gamma}^{*}\left\{
\frac{{\rm Td}({\mathcal H}^{\lor})^{-1}-1}{c_{1}({\mathcal H}^{\lor})}\right\}\,
q^{*}\{{\rm Td}_{g}(T{\mathcal X}){\rm ch}_{g}(\omega_{\mathcal X}(\xi))\}
\\
&\quad
-\int_{{\mathcal X}^{g}_{V}\cap X_{0}}{\rm Td}_{g}(T{\mathcal X}){\rm ch}_{g}(\omega_{\mathcal X}(\xi)).
\end{aligned}
$$
\par
Let $f\colon(Y,Y_{0})\to(T,0)$ be a semistable reduction of 
$\pi\colon(X,X_{0})\to(S,0)$. We have a commutative diagram, 
where $Y_{0}\subset Y$ is a reduced normal crossing divisor:
$$
\begin{CD}
(Y,Y_{0}=f^{-1}(0))@>F>> (X,X_{0})
\\
@V f VV  @V \pi VV
\\
(T,0) @> \mu >> (S,0).
\end{CD}
$$
By \cite{MourouganeTakayama09},
$R^{q}f_{*}\omega_{Y/T}(F^{*}\xi)$ is a $G$-equivariant locally free sheaf equipped with 
an injective homomorphism
$\varphi\colon R^{q}f_{*}\omega_{Y/T}(F^{*}\xi)\to\mu^{*}R^{q}\pi_{*}\omega_{X/S}(\xi)$ 
of $G$-modules. We regard $R^{q}f_{*}\omega_{Y/T}(F^{*}\xi)$ as a subsheaf of 
$\mu^{*}R^{q}\pi_{*}\omega_{X/S}(\xi)$ by this inclusion.
The Lefschetz trace of the $G$-action on
$(\mu^{*}R\pi_{*}\omega_{X/S}(\xi)/Rf_{*}\omega_{Y/T}(F^{*}\xi))_{0}$ is defined as
$$
{\Bbb L}_{g}\left(
\frac{\mu^{*}R\pi_{*}\omega_{X/S}(\xi)}{Rf_{*}\omega_{Y/T}(F^{*}\xi)}
\right)
=
\sum_{q\geq0}(-1)^{q}
{\rm Tr}\left[g|_{(\mu^{*}R^{q}\pi_{*}\omega_{X/S}(\xi)/R^{q}f_{*}\omega_{Y/T}(F^{*}\xi))_{0}}\right]
$$
for $g\in G$. Now the main result of this article is stated as follows.

\begin{theorem}\label{Theorem1.2}
If $(\xi,h_{\xi})$ is Nakano semi-positive on $X=\pi^{-1}(S)$, then there exist constants
$\nu_{g},c_{g}\in{\bf C}$ such that as $s\to0$
$$
\begin{aligned}
\log\tau_{G}(X_{s},\omega_{X_{s}}(\xi_{s}))(g)
&=
\{
\alpha_{g}(X_{0},\omega_{X/S}(\xi))
+
\frac{1}{\deg\mu}
{\Bbb L}_{g}\left(
\frac{\mu^{*}R\pi_{*}\omega_{X/S}(\xi)}{Rf_{*}\omega_{Y/T}(F^{*}\xi)}
\right)
\}
\log |s|^{2}
\\
&\qquad
+\nu_{g}\,\log(-\log|s|^{2})+c_{g}+O\left(1/\log|s|\right).
\end{aligned}
$$
\end{theorem}

By Theorem~\ref{Theorem1.2}, the logarithmic singularity of
$\log\tau_{G}(X_{s},\omega_{X_{s}}(\xi_{s}))(g)$ is determined by the algebraic term 
${\Bbb L}_{g}(\mu^{*}R\pi_{*}\omega_{X/S}(\xi)/Rf_{*}\omega_{Y/T}(F^{*}\xi))$
measuring the cohomological difference between $(X,X_{0})$ and its semistable reduction $(Y,Y_{0})$ 
and the topological term $\alpha_{g}(X_{0},\omega_{X/S}(\xi))$
arising from the resolution of the Gauss map. See Corollary~\ref{Corollary6.8}
for a formula for $\log\tau_{G}(X_{s},\xi_{s})(g)$ as $s\to0$, when $(\xi,h_{\xi})$ is 
{\em semi-negative in the dual Nakano sense}. 
We remark $\nu_{g}\in{\bf Z}$ if $g=1$.
\par
For some important classes of degenerations, we get the following.

\begin{theorem}\label{Theorem1.4}
Assume that $(\xi,h_{\xi})$ is Nakano semi-positive on $X=\pi^{-1}(S)$.
\begin{itemize}
\item[(1)]
If $X_{0}$ is reduced and the pair $({\mathcal X},X_{0})$ has only log-canonical
singularities, then 
${\Bbb L}_{g}(\mu^{*}R\pi_{*}\omega_{X/S}(\xi)/Rf_{*}\omega_{Y/T}(F^{*}\xi))=0$ in Theorem~\ref{Theorem1.2}.
\item[(2)]
If $X_{0}$ is reduced, normal and has only canonical (equivalently rational) 
singularities, then there exist $c\in{\bf C}$, 
$r\in{\bf Q}_{>0}$, $l\in{\bf Z}_{\geq0}$ such that as $s\to0$,
$$
\log\tau_{G}(X_{s},\omega_{X_{s}}(\xi_{s}))(g)
=
\alpha_{g}(X_{0},\omega_{X/S}(\xi))\log|s|^{2}+c
+
O\left(|s|^{r}(\log|s|)^{l}\right).
$$
\end{itemize}
\end{theorem}

For the definition of (log-)canonical singularities, see Sect.\ref{sect:7.1}.  
Since the pair $({\mathcal X},X_{0})$ has only log-canonical singularities when $\pi\colon(X,X_{0})\to(S,0)$
is semistable, Theorem~\ref{Theorem1.2} and Theorem~\ref{Theorem1.4} (1) are compatible.
If the singularity of $({\mathcal X},X_{0})$ is strictly worse than log-canonical, then
${\Bbb L}_{g}(\mu^{*}R\pi_{*}\omega_{X/S}(\xi)/Rf_{*}\omega_{Y/T}(F^{*}\xi))\not=0$
in general (Sect.\,\ref{sect:8}). 
When $\dim X_{s}=1$, $G=\{1\}$ and $X_{0}$ has at most nodes, asymptotic expansions like 
Theorems~\ref{Theorem1.2} and~\ref{Theorem1.4} (1)
were obtained by Bismut--Bost \cite{BismutBost90} and Wolpert \cite{Wolpert87}, 
where $\nu_{g}\not=0$ in general.
Without the assumption of the curvature of $(\xi,h_{\xi})$, the singularity of analytic torsion was 
determined by Farber \cite{Farber95} when $X_{0}$ is {\em non-singular}. 
In \cite{Yoshikawa09b}, we use Theorem~\ref{Theorem1.2} to prove
the automorphic property of the invariant $\tau_{M}$ introduced
in \cite{Yoshikawa04}, which plays a crucial role to determine an explicit Borcherds product 
expressing $\tau_{M}$ in the case $r(M)\geq18$ (cf. \cite{Yoshikawa09a}).
\par
The strategy for the proof of Theorem~\ref{Theorem1.2} is as follows.
Following \cite[Th.\,5.9]{Bismut97}, \cite[Th.\,1.1]{Yoshikawa07},
we determine the singularity of the equivariant Quillen metric on the equivariant determinant
of the cohomologies of $\xi$ by applying the Bismut immersion formula \cite{Bismut95} to 
the $G$-equivariant embedding $X_{s}\hookrightarrow{\mathcal X}$ (Theorem~\ref{Theorem4.1}).
Then we get Theorem~\ref{Theorem1.2} by studying the behavior of the $L^{2}$-metric 
on $R^{q}\pi_{*}\omega_{X/S}(\xi)$ (Theorem~\ref{Theorem6.8}). 
For this, a theorem of Takegoshi \cite{Takegoshi95} and its extension by
Mourougane-Takayama \cite{MourouganeTakayama09} play a crucial role to 
express the fiberwise harmonic representative of an element of $H^{q}(X,\Omega_{X}^{n+1}(\xi))$,
where the Nakano semi-positivity of $\xi$ is essentially used (Sect.\,\ref{sect:6}).
The asymptotic expansion follows from a theorem of Barlet \cite{Barlet82}.
\par
This article is organized as follows.
In Sect.\ref{sect:2}, we recall Gauss maps.
In Sect.\ref{sect:3}, we recall equivariant Quillen metrics and study their regularity.
In Sect.\ref{sect:4}, we determine the singularity of equivariant Quillen metrics.
In Sect.\ref{sect:5}, we recall the notion of Nakano semi-positivity of vector bundles.
In Sect.\ref{sect:6}, we prove Theorem~\ref{Theorem1.2}.
In Sect.\ref{sect:7}, we prove Theorem~\ref{Theorem1.4}.
In Sect.\ref{sect:8}, we study some examples.
In Sect.\ref{sect:9}, we prove some technical results.
\par{\bf Notation }
For a complex manifold, we set
$d^{c}=\frac{1}{4\pi i}(\partial-\bar{\partial})$.
Hence
$dd^{c}=\frac{1}{2\pi i}\bar{\partial}\partial$.
For a complex manifold $Y$, $A^{p,q}_{Y}$ denotes the vector space of
$C^{\infty}$ $(p,q)$-forms on $Y$.
We set
$\widetilde{A}_{Y}=
\bigoplus_{p\geq0}A^{p,p}_{Y}/{\rm Im}\,\partial+{\rm Im}\,\bar{\partial}$.
For a $G$-equivariant vector bundle $F$ over $Y$
equipped with a $G$-invariant Hermitian metric $h_{F}$, we denote by
$c_{i}(F,h_{F})\in\bigoplus_{p\geq0}A^{p,p}_{Y}$,
${\rm Td}_{g}(F, h_{F}),{\rm ch}_{g}(F,h_{F})\in
\bigoplus_{p\geq0}A^{p,p}_{Y^{g}}$
the $i$-th Chern form, the equivariant Todd form, and 
the equivariant Chern character form of $(F,h_{F})$ with respect to 
the holomorphic Hermitian connection, respectively (cf. \cite{Bismut95}). 
\par
After \cite{Barlet82}, we set
${\mathcal B}(S)=
C^{\infty}(S)\oplus\bigoplus_{r\in{\bf Q}\cap(0,1]}\bigoplus_{k=0}^{n}|s|^{2r}(\log|s|)^{k}C^{\infty}(S)$.
A function of ${\mathcal B}(S)$ is continuous and has an asymptotic expansion at $s=0$.
We write $\phi\equiv_{\mathcal B}\psi$ if $\phi,\psi\in C^{\infty}(S^{o})$ satisfies 
$\phi-\psi\in{\mathcal B}(S)\subset C^{0}(S)$.
Throughout this article, we keep the notation and the assumptions in Sect.\ref{sect:1}.
\par
{\bf Acknowledgements }
We thank Professor Shigeharu Takayama for helpful discussions about 
the singularity of $L^{2}$-metrics and for pointing out some errors and difficulties in the earlier version.
We also thank Professor Vincent Maillot for his comments on the earlier version, 
which improved the formulation of Theorem~\ref{Theorem1.2} and Professor Fumiharu Kato
for helpful discussion about semistable reduction.


\section
{The Gauss map and its equivariant resolution}
\label{sect:2}
\par
Let $\Omega^{1}_{\mathcal X}$ be the holomorphic cotangent bundle of ${\mathcal X}$.
Let $\varPi\colon{\bf P}(\Omega^{1}_{\mathcal X}\otimes\pi^{*}TC)\to{\mathcal X}$ be 
the projective-space bundle associated with $\Omega^{1}_{\mathcal X}\otimes\pi^{*}TC$. 
Let $\varPi^{\lor}\colon{\bf P}(T{\mathcal X})^{\lor}\to{\mathcal X}$ be the dual projective-space 
bundle of ${\bf P}(T{\mathcal X})$, whose fiber ${\bf P}(T_{x}{\mathcal X})^{\lor}$ is the set of 
hyperplanes of $T_{x}{\mathcal X}$ passing through $0_{x}\in T_{x}{\mathcal X}$. 
Since $\dim C=1$, we have
${\bf P}(\Omega^{1}_{\mathcal X}\otimes\pi^{*}TC)=
{\bf P}(\Omega^{1}_{\mathcal X})\cong{\bf P}(T{\mathcal X})^{\lor}$.
\par
Let $x\in{\mathcal X}\setminus\Sigma_{\pi}$.
Let $s$ be a holomorphic local coordinate of $C$ near $\pi(x)\in C$. 
We define the Gauss maps
$\nu\colon{\mathcal X}\setminus\Sigma_{\pi}\to{\bf P}(\Omega^{1}_{\mathcal X}\otimes\pi^{*}TC)$ 
and $\gamma\colon{\mathcal X}\setminus\Sigma_{\pi}\to{\bf P}(T{\mathcal X})^{\lor}$ by
$$
\nu(x)
:=
[d\pi_{x}]
=
\left[
\sum_{i=0}^{n}\frac{\partial(s\circ\pi)}
{\partial z_{i}}(x)\,dz_{i}\otimes
\frac{\partial}{\partial s}
\right],
\qquad
\gamma(x)
:=
[T_{x}{\mathcal X}_{\pi(x)}].
$$
Under the canonical isomorphism
${\bf P}(\Omega^{1}_{\mathcal X}\otimes\pi^{*}TC)\cong{\bf P}(T{\mathcal X})^{\lor}$, 
one has $\nu=\gamma$.
\par
Let ${\mathcal H}:={\mathcal O}_{{\bf P}(T{\mathcal X})^{\lor}}(1)$ be the tautological quotient bundle 
and let ${\mathcal U}$ be the universal hyperplane bundle of ${\bf P}(T{\mathcal X})^{\lor}$. 
We have the exact sequence of $G$-equivariant vector bundles on ${\bf P}(T{\mathcal X})^{\lor}$
$$
{\mathcal S}^{\lor}\colon
0
\longrightarrow 
{\mathcal U}
\longrightarrow
(\varPi^{\lor})^{*}T{\mathcal X}
\longrightarrow
{\mathcal H}
\longrightarrow 
0.
$$
\par
Let $h_{\mathcal U}$ be the Hermitian metric on ${\mathcal U}$ induced from $(\varPi^{\lor})^{*}h_{\mathcal X}$, 
and let $h_{\mathcal H}$ be the Hermitian metric on ${\mathcal H}$ induced from 
$(\varPi^{\lor})^{*}h_{\mathcal X}$ by the $C^{\infty}$-isomorphism ${\mathcal H}\cong{\mathcal U}^{\perp}$.
On ${\mathcal X}\setminus\Sigma_{\pi}$, we have
$(T{\mathcal X}/C,h_{{\mathcal X}/C})=\gamma^{*}({\mathcal U},h_{\mathcal U})$.
\par
Let 
${\mathcal L}:={\mathcal O}_{{\bf P}(\Omega^{1}_{\mathcal X}\otimes\pi^{*}TC)}(-1)
\subset
\varPi^{*}(\Omega^{1}_{\mathcal X}\otimes\pi^{*}TC)$
be the tautological line bundle over ${\bf P}(\Omega^{1}_{\mathcal X}\otimes\pi^{*}TC)$. 
Let $h_{C}$ be a Hermitian metric on $C$.
Let $h_{\Omega^{1}_{\mathcal X}}$ be the Hermitian metric on $\Omega^{1}_{\mathcal X}$ 
induced from $h_{\mathcal X}$.
Let $h_{\mathcal L}$ be the Hermitian metric on ${\mathcal L}$ induced from the metric
$\varPi^{*}(h_{\Omega^{1}_{\mathcal X}}\otimes\pi^{*}h_{C})$ 
by the inclusion 
${\mathcal L}\subset\varPi^{*}(\Omega^{1}_{\mathcal X}\otimes\pi^{*}TC)$. 
\par
Since $\Sigma_{\pi}$ is a proper subvariety of ${\mathcal X}$,
the Gauss maps $\nu$ and $\gamma$ extend to rational maps
$\nu\colon{\mathcal X}\dashrightarrow{\bf P}(\Omega^{1}_{\mathcal X}\otimes\pi^{*}TC)$
and $\gamma\colon{\mathcal X}\dashrightarrow{\bf P}(T{\mathcal X})^{\lor}$.
By \cite[Th.\,13.2]{BierstoneMilman97}, 
there exist a projective algebraic manifold $\widetilde{\mathcal X}$, 
a normal crossing divisor $E\subset\widetilde{\mathcal X}$, 
a birational holomorphic map $q\colon\widetilde{\mathcal X}\to{\mathcal X}$ with $E=q^{-1}(\Sigma_{\pi})$, 
and holomorphic maps
$\widetilde{\nu}\colon\widetilde{\mathcal X}\to{\bf P}(\Omega^{1}_{\mathcal X}\otimes\pi^{*}TC)$ 
and $\widetilde{\gamma}\colon\widetilde{\mathcal X}\to{\bf P}(T{\mathcal X})^{\lor}$ 
with the following properties:
\begin{itemize}
\item[(a)]
The $G$-action on ${\mathcal X}$ lifts to a $G$-action on $\widetilde{\mathcal X}$ and
$q|_{\widetilde{\mathcal X}\setminus E}\colon
\widetilde{\mathcal X}\setminus E\to{\mathcal X}\setminus\Sigma_{\pi}$ 
is a $G$-equivariant isomorphism.
\item[(b)]
$(\pi\circ q)^{-1}(b)$ is a normal crossing divisor of $\widetilde{\mathcal X}$ for all $b\in\Delta$.
\item[(c)]
$\widetilde{\nu}=\nu\circ q$ and $\widetilde{\gamma}=\gamma\circ q$ 
on $\widetilde{\mathcal X}\setminus E$.
\end{itemize}
Then $\widetilde{\nu}=\widetilde{\gamma}$ under the canonical isomorphism
${\bf P}(\Omega^{1}_{\mathcal X}\otimes\pi^{*}TC)\cong{\bf P}(T{\mathcal X})^{\lor}$.

\section
{Equivariant Quillen metrics}
\label{sect:3}
\par
In this section, 
we recall equivariant Quillen metrics and prove its smoothness for smooth projective morphisms. 
For more general treatments including smooth K\"ahler morphisms, we refer to
\cite[III Sects.\,2 and 3]{BGS88}.
In the rest of this article, $\widehat{G}$ denotes the set of equivalence classes of complex irreducible 
representations of $G$. For $W\in\widehat{G}$, the corresponding irreducible character is denoted 
by $\chi_{W}$.

\subsection
{\bf Equivariant analytic torsion and equivariant Quillen metrics}
\label{subsect:3.1}

\subsubsection
{Equivariant analytic torsion}
\label{subsubsect:3.1.1}
\par
Let $V$ be a compact K\"ahler manifold with holomorphic $G$-action. 
Let $h_{V}$ be a $G$-invariant K\"ahler metric on $V$. 
Let $F$ be a $G$-equivariant holomorphic vector bundle on $V$ equipped with a $G$-invariant 
Hermitian metric $h_{F}$.
We set $\overline{V}=(V,h_{V})$ and $\overline{F}=(F,h_{F})$.
Let $A^{p,q}_{V}(F)$ be the vector space of $F$-valued smooth $(p,q)$-forms on $V$. 
We set $S_{F}=\bigoplus_{q\geq0}A^{0,q}_{V}(F)$, which is equipped with the $L^{2}$ metric 
$(\cdot,\cdot)_{L^{2}}$ with respect to $h_{V}$ and $h_{F}$. 
Then $(\cdot,\cdot)_{L^{2}}$ is $G$-invariant with respect to the standard $G$-action on $S_{F}$.
\par
Let $\square_{F}=2(\bar{\partial}_{F}+\bar{\partial}^{*}_{F})^{2}$ be the Laplacian acting on $S_{F}$. 
We denote by $\sigma(\square_{F})$ the spectrum of $\square_{F}$.
Let $K_{F}(\lambda)$ be the eigenspace of $\square_{F}$ 
with eigenvalue $\lambda\in\sigma(\square_{F})$. 
Since $G$ preserves the metrics $h_{V}$ and $h_{F}$, 
$\square_{F}$ commutes with the $G$-action on $S_{F}$. Hence $G$ acts on $K_{F}(\lambda)$.
With respect to the ${\bf Z}$-grading on $K_{F}(\lambda)$ induced from the one on $S_{F}$, 
the number operator $N$ and the supertrace ${\rm Tr}_{\rm s}[\cdot]$ are defined on $K_{F}(\lambda)$ 
(cf. \cite{BGS88}).
For $g\in G$ and $s\in{\bf C}$ with ${\rm Re}\,s\gg0$, we define
$$
\zeta_{G}(g)(s)
:=
\sum_{\lambda\in\sigma(\square_{F})\setminus\{0\}}
\lambda^{-s}\,{\rm Tr}_{\rm s}\,[g\,N|_{K_{F}(\lambda)}].
$$ 
Then $\zeta_{G}(g)(s)$ extends to a meromorphic function on ${\bf C}$ and is holomorphic at $s=0$. 
For $g\in G$, we define
$$
\log\tau_{G}(\overline{V},\overline{F})(g)
:=
-\zeta_{G}'(g)(0).
$$
Then $\log\tau_{G}(\overline{V},\overline{F})(g)$ is called the equivariant analytic torsion of 
$(\overline{V},\overline{F})$ (cf. \cite{Bismut95}).

\subsubsection
{Equivariant Quillen metrics}
\label{subsubsect:3.1.2}
\par
Since $F$ is a $G$-equivariant holomorphic vector bundle on $V$, $G$ acts on 
$H(V,F)=\bigoplus_{q\in{\bf Z}}H^{q}(V,F)$ and preserves its grading. 
One has the isotypical splitting of the ${\bf Z}$-graded vector space
$$
H(V,F)=\bigoplus_{W\in\widehat{G}}{\rm Hom}_{G}(W,H(V,F))\otimes W.
$$
Since $\dim H(V,F)<+\infty$, ${\rm Hom}_{G}(W,H(V,F))=0$ except for finite $W\in\widehat{G}$. 
We set
$$
\lambda_{W}(F)
=
\det{\rm Hom}_{G}(W,H(V,F))\otimes W
:=
\bigotimes_{q\geq0}(\det {\rm Hom}_{G}(W,H^{q}(V,F))\otimes W)^{(-1)^{q}}.
$$
The equivariant determinant of the cohomologies of $F$ is defined as
$$
\lambda_{G}(F):=\prod_{W\in\widehat{G}}\lambda_{W}(F).
$$
Notice that our sign convention is different form the one in \cite[Eq.\,(2.9)]{Bismut95}. 
When ${\rm Hom}_{G}(W,H(V,F))=0$, $\lambda_{W}(F)$ is canonically isomorphic to ${\bf C}$ by definition. 
In this case, the canonical element of $\lambda_{W}(F)$ corresponding to $1\in{\bf C}$ is denoted by 
$1_{\lambda_{W}(F)}$.
A vector $\alpha=(\alpha_{W})_{W\in\widehat{G}}\in\lambda_{G}(F)$ is said to be {\em admissible} 
if $\alpha_{W}\not=0$ for all $W\in\widehat{G}$ and if $\alpha_{W}=1_{\lambda_{W}(F)}$ except for 
finitely many $W\in\widehat{G}$. The set of admissible elements of $\lambda_{G}(F)$ 
is identified with the direct sum $\bigoplus_{W\in\widehat{G}}\lambda_{W}(F)^{\times}$, 
where $\lambda_{W}(F)^{\times}:=\lambda_{W}(F)\setminus\{0\}$
is the set of invertible elements of $\lambda_{W}(F)$.
\par
By Hodge theory, we have an isomorphism of ${\bf Z}$-graded $G$-spaces $H(V,F)\cong K_{F}(0)$. 
The $G$-invariant metric on $H(V,F)$ induced from the $L^{2}$-metric on $K_{F}(0)\subset S_{F}$ 
by this isomorphism is denoted by $h_{H(V,F)}$. Then the isotypical splitting of $H(V,F)$ is orthogonal 
with respect to $h_{H(V,F)}$. Let $\|\cdot\|_{L^{2},\lambda_{W}(F)}$ be the Hermitian metric 
on $\lambda_{W}(F)$ induced from $h_{H(V,F)}$. 
Recall that $\chi_{W}$ is the character of $W\in\widehat{G}$.
For an admissible element
$\alpha=(\alpha_{W})_{W\in\widehat{G}}\in\bigoplus_{W\in\widehat{G}}\lambda_{W}(F)^{\times}$, we set
$$
\log\|\alpha\|^{2}_{Q,\lambda_{G}(F)}(g)
:=
-\zeta_{G}'(g)(0)+
\sum_{W\in\widehat{G}}\frac{\chi_{W}(g)}{\dim W}
\log\|\alpha_{W}\|_{L^{2},\lambda_{W}(F)}^{2}.
$$
The ${\bf C}$-valued function $\log\|\cdot\|^{2}_{Q,\lambda_{G}(F)}(g)$ on
$\bigoplus_{W\in\widehat{G}}\lambda_{W}(F)^{\times}$ is called the equivariant Quillen metric 
on $\lambda_{G}(F)$ with respect to $h_{V}$, $h_{F}$. 
Notice that equivariant Quillen metric makes sense only for admissible elements.
We refer to \cite{RaySinger73}, \cite{Quillen85},  \cite{Bismut95}, \cite{KohlerRoessler01}, \cite{Ma00} 
for more about equivariant analytic torsion, equivariant determinants, and equivariant Quillen metrics.

\subsection
{\bf The smoothness of equivariant Quillen metrics}
\label{subsect:3.2}
\par
Let $M$ be a compact complex manifold with holomorphic $G$-action and 
let $B$ be a compact complex manifold with trivial $G$-action.
Let $\pi\colon M\to B$ be a $G$-equivariant proper surjective flat holomorphic map. 
Assume that there is a $G$-equivariant ample line bundle on $M$. 
Let $F\to M$ be a $G$-equivariant holomorphic vector bundle. 
Let $h_{M}$ be a $G$-invariant K\"ahler metric on $M$ and 
let $h_{F}$ be a $G$-invariant Hermitian metric on $F$. 
We set $M_{b}:=\pi^{-1}(b)$ and $F_{b}:=F|_{M_{b}}$ for $b\in B$. 
Then $G$ preserves the fibers $M_{b}$ and $F_{b}$, 
and $\pi\colon M\to B$ is a family of projective algebraic varieties.
\par
For $W\in\widehat{G}$, we define ${\rm Hom}_{G}(W,R^{q}\pi_{*}F)\otimes W$
to be the sheaf on $B$ associated to the presheaf
$U\mapsto{\rm Hom}_{G}\left(W,H^{q}(\pi^{-1}(U),F|_{\pi^{-1}(U)})\right)\otimes W$.

\subsubsection
{Direct image sheaves}
\label{subsubsect:3.2.1}
\par
Since there is a $G$-equivariant ample line bundle on $M$ by assumption, 
there exist a complex of $G$-equivariant holomorphic vector bundles 
$F_{\bullet}\colon 
0\longrightarrow 
F_{0}\longrightarrow 
F_{1}\longrightarrow
\cdots\longrightarrow 
F_{m}\longrightarrow0$
on $M$ and a homomorphism $i\colon F\to F_{0}$ of $G$-modules with the following conditions:
\begin{itemize}
\item[(i)]
The complex $0\to F\to F_{0}\to\cdots\to F_{m}\to0$ is acyclic.
\item[(ii)]
$H^{q}(M_{b},F_{i}|_{M_{b}})=0$ for all $q>0$ and $i\geq0$. 
\end{itemize}
By (i), (ii),
$\pi_{*}F_{\bullet}
\colon
0\longrightarrow\pi_{*}F_{0}
\longrightarrow\pi_{*}F_{1}
\longrightarrow\cdots\longrightarrow
\pi_{*}F_{m}\longrightarrow0$ 
is a complex of $G$-equivariant locally free sheaves of finite rank over $B$, 
whose cohomology sheaves compute the direct image sheaves $R^{q}\pi_{*}F$, $q\geq0$, i.e.,
\begin{equation}\label{eqn:(3.1)}
R^{q}\pi_{*}F={\mathcal H}^{q}(\pi_{*}F_{\bullet})
:=
\ker\{\pi_{*}F_{q}\to\pi_{*}F_{q+1}\}/{\rm Im}\{\pi_{*}F_{q-1}\to\pi_{*}F_{q}\}.
\end{equation}
Since \eqref{eqn:(3.1)} is an equality of $G$-modules, 
we get for all $W\in\widehat{G}$ and $q\geq0$
\begin{equation}\label{eqn:(3.2)}
{\rm Hom}_{G}(W,R^{q}\pi_{*}F)\otimes{W}
=
{\mathcal H}^{q}({\rm Hom}_{G}(W,\pi_{*}F_{\bullet})\otimes{W}).
\end{equation}
Since ${\rm Hom}_{G}(W,\pi_{*}F_{\bullet})\otimes{W}$ is a complex of locally free ${\mathcal O}_{B}$-modules, 
we deduce from \eqref{eqn:(3.2)} that ${\rm Hom}_{G}(W,R^{q}\pi_{*}F)\otimes{W}$ is a coherent
${\mathcal O}_{B}$-module on $B$. Since
${\mathcal H}^{q}(\pi_{*}F_{\bullet})=
\bigoplus_{W\in\widehat{G}}{\mathcal H}^{q}({\rm Hom}_{G}(W,\pi_{*}F_{\bullet})\otimes{W})$,
we get the isotypical splitting on $B$ 
\begin{equation}\label{eqn:(3.3)}
R^{q}\pi_{*}F=\bigoplus_{W\in\widehat{G}}{\rm Hom}_{G}(W,R^{q}\pi_{*}F)\otimes W.
\end{equation}
Notice that ${\rm Hom}_{G}(W,R^{q}\pi_{*}F)=0$ except for finitely many $W\in\widehat{G}$.

\subsubsection
{Equivariant determinant of cohomologies}
\label{subsubsect:3.2.2}
\par
Recall that for any coherent analytic sheaf ${\mathcal F}$ on $B$, 
one can associate the invertible sheaf $\det{\mathcal F}$ on $B$ 
by \cite{KnudsenMumford76}, \cite[III, Sect.\,3]{BGS88}. We define
$$
\begin{array}{ll}
\lambda_{W}(F)
&:=
\bigotimes_{q\geq0}
\det\left({\rm Hom}_{G}(W,R^{q}\pi_{*}F)\otimes W\right)^{(-1)^{q}},
\\
\lambda_{G}(F)
&:=
\prod_{W\in\widehat{G}}\lambda_{W}(F).
\end{array}
$$
If ${\rm Hom}_{G}(W,R^{q}\pi_{*}F)=0$ for all $q\geq0$,
then $\lambda_{W}(F)$ is canonically isomorphic to ${\mathcal O}_{B}$.
In this case, the canonical section of $\lambda_{W}(F)$ corresponding to 
$1\in H^{0}(B,{\mathcal O}_{B})$ is denoted by $1_{\lambda_{W}(F)}$. 
By \cite[III, Lemma 3.7]{BGS88}, there is a canonical identification
\begin{equation}\label{eqn:(3.4)}
\lambda_{W}(F)_{b}
:=
\lambda_{W}(F)/{\frak m}_{b}\lambda_{W}(F)
=
\lambda_{W}(F|_{M_{b}})
\end{equation}
for all $b\in B$, where ${\frak m}_{b}$ is the maximal ideal of ${\mathcal O}_{B,b}$.
\par
For an open subset $U\subset B$, a holomorphic section $\sigma=(\sigma_{W})_{W\in\widehat{G}}$   
of $\lambda_{G}(F)|_{U}$ is said to be {\em admissible}
if $\sigma_{W}$ is nowhere vanishing on $U$ for all $W\in\widehat{G}$ and 
if $\sigma_{W}=1_{\lambda_{W}(\eta)}$ except for finitely many $W\in\widehat{G}$.

\subsubsection
{The smoothness of equivariant Quillen metrics}
\label{subsubsect:3.2.3}
\par
Let $D\subset B$ be the discriminant locus of $\pi$ and set $B^{o}:=B\setminus D$.
By \eqref{eqn:(3.4)}, $\lambda_{G}(F)|_{B^{o}}$ is equipped with the equivariant Quillen metric 
$\|\cdot\|_{Q,\lambda_{G}(F)}(g)$ with respect to $h_{M}|_{TM/B}$ and $h_{F}$ 
such that for every $b\in B^{o}$,
$$
\|\cdot\|_{Q,\lambda_{G}(F)}(g)(b):=\|\cdot\|_{Q,\lambda_{G}(F_{b})}(g).
$$

\begin{theorem}\label{Theorem3.1}
Let $U\subset B^{o}$ be a small connected open subset with
$\lambda_{W}(F_{i})|_{U}\cong{\mathcal O}_{U}$ for all $i$ and $W\in\widehat{G}$.
Let $\sigma=(\sigma_{W})_{W\in\widehat{G}}$ be an admissible holomorphic section of $\lambda_{G}(F)|_{U}$. 
Then $\log\|\sigma\|_{Q,\lambda_{G}(F)}^{2}(g)\in C^{\infty}(U)$.
\end{theorem}

\begin{pf}
By the definition of $\lambda_{W}(F)$, we have the canonical isomorphism
$$
\varphi_{W}\colon
\lambda_{W}(F)
=
\bigotimes_{q\geq0}
\det{\mathcal H}^{q}({\rm Hom}_{G}(W,\pi_{*}F_{\bullet})\otimes{W})^{(-1)^{q}}
\cong
\bigotimes_{i\geq0}\lambda_{W}(F_{i})^{(-1)^{i}}.
$$
There exists an admissible holomorphic section $\sigma_{i}=(({\sigma}_{i})_{W})_{W\in\widehat{G}}$ of 
$\lambda_{G}(F_{i})|_{U}$ such that $\varphi_{W}(\sigma_{W})=\otimes_{i\geq0}({\sigma}_{i})_{W}^{(-1)^{i}}$
for all $W\in\widehat{G}$. 
\par
Let $h_{F_{i}}$ be a $G$-invariant Hermitian metric on $F_{i}$ and 
let $\|\cdot\|_{Q,\lambda_{G}(F_{i})}(g)$ be the equivariant Quillen metric on $\lambda_{G}(F_{i})$ 
with respect to the $G$-invariant metrics $h_{M/B}=h_{M}|_{TM/B}$ and $h_{F_{i}}$. 
Let $\widetilde{\rm ch}_{g}(F,F_{\bullet};h_{F},h_{F_{\bullet}})\in\widetilde{A}(M^{g})$ 
be the Bott--Chern secondary class \cite[I e), f)]{BGS88} such that
$$
dd^{c}\widetilde{\rm ch}_{g}(F,F_{\bullet};h_{F},h_{F_{\bullet}})
=
\sum_{i\geq0}(-1)^{i}{\rm ch}_{g}(F_{i},h_{F_{i}})-{\rm ch}_{g}(F,h_{F}).
$$
Applying the immersion formula of Bismut \cite[Th.\,0.1]{Bismut95} to 
the immersion $\emptyset\hookrightarrow M_{b}$, $b\in U$, 
we get the following equation of complex-valued functions on $U$
\begin{equation}\label{eqn:(3.5)}
\begin{aligned}
\log\|\sigma\|_{Q,\lambda_{G}(F)}^{2}(g)
&=
\sum_{i\geq0}(-1)^{i}
\log\|\sigma_{i}\|_{Q,\lambda_{G}(F_{i})}^{2}(g)
\\
&\quad
+
\left[
\pi_{*}
\{{\rm Td}_{g}(TM/B,h_{M/B})
\widetilde{\rm ch}_{g}(F,F_{\bullet};h_{F},h_{F_{\bullet}})\}
\right]^{(0)}
\\
&\equiv
\sum_{i\geq0}(-1)^{i}
\log\|\sigma_{i}\|_{Q,\lambda_{G}(F_{i})}^{2}(g)
\mod C^{\infty}(U),
\end{aligned}
\end{equation}
where we used \cite[I, Th.\,1.29 and Cor.\,1.30]{BGS88} to identify 
the Bott--Chern current $T_{g}(F,F_{\bullet};h_{F},h_{F_{\bullet}})$ with the Bott--Chern class
$\widetilde{\rm ch}_{g}(F,F_{\bullet};h_{F},h_{F_{\bullet}})$ in the first equality 
and $[\omega]^{(2d)}$ denotes the component of degree $2d$ of a differential form $\omega$. 
Since the morphism $\pi\colon\pi^{-1}(U)\to U$ is smooth and
since $h^{0}(M_{b},F_{i}|_{M_{b}})$ is a constant function on $U$, we get
$\log\|\sigma_{i}\|_{Q,\lambda_{G}(F_{i})}^{2}(g)\in C^{\infty}(U)$ by \cite[III, Th.\,3.5]{BGS88}. 
This, together with \eqref{eqn:(3.5)}, implies the result.
\end{pf}

\begin{remark}\label{Remark3.2}
The curvature $-dd^{c}\log\|\sigma\|^{2}_{Q,\lambda_{G}(F)}(g)$ was computed by 
Bismut--Gillet--Soul\'e \cite[Th.\,0.1]{BGS88} when $G$ is trivial and by Ma \cite[Th.\,2.12]{Ma00}
when $G$ is general and all direct image sheaves $R^{q}\pi_{*}F$, $q\geq0$ are locally free. 
By Theorem~\ref{Theorem3.1}, the curvature formula of Bismut-Gillet-Soul\'e-Ma
\begin{equation}\label{eqn:(3.6)}
-dd^{c}\log\|\sigma\|_{Q,\lambda_{G}(F)}^{2}(g)
=
\left[\pi_{*}\{{\rm Td}_{g}(TM/B,h_{M/B})\,{\rm ch}_{g}(F,h_{F})\}\right]^{(2)}
\end{equation}
remains valid on $B^{o}$ even if $R^{q}\pi_{*}F$ may not be locally free, 
because there is a dense Zariski open subset $U\subset B^{o}$ over 
which $h^{q}(F_{b})$ is constant for all $q\geq0$.
\end{remark}

\section
{The singularity of equivariant Quillen metrics}
\label{sect:4}
\par
Let $\lambda_{G}(\xi)$ be the equivariant determinant of the cohomologies of $\xi$. 
Then $\lambda_{G}(\xi)|_{C^{o}}$ is equipped with the equivariant Quillen metric
$\|\cdot\|^{2}_{\lambda_{G}(\xi),Q}(\cdot)$ with respect to the $G$-invariant metrics 
$h_{{\mathcal X}/C}$, $h_{\xi}$. 
Let $\sigma$ be an admissible holomorphic section of $\lambda_{G}(\xi)|_{S}$. 
By Theorem~\ref{Theorem3.1}, $\log\|\sigma(s)\|^{2}_{\lambda_{G}(\xi),Q}(g)\in C^{\infty}(S^{o})$ 
for $g\in G$. In this section, we determine the behavior of $\log\|\sigma(s)\|^{2}_{\lambda_{G}(\xi),Q}(g)$ 
as $s\to0$.

\subsection
{The non-twisted case}
\label{sect:4.1}
\par
Let $\Gamma\subset{\mathcal X}\times C$ be the graph of $\pi$. 
Then $\Gamma$ is a smooth divisor on ${\mathcal X}\times C$ preserved by the $G$-action 
on ${\mathcal X}\times C$. 
Let $[\Gamma]$ be the $G$-equivariant holomorphic line bundle on ${\mathcal X}\times C$ 
associated to $\Gamma$. 
Let $\varsigma_{\Gamma}\in H^{0}({\mathcal X}\times C,[\Gamma])$ be the canonical section of 
$[\Gamma]$ such that ${\rm div}(\varsigma_{\Gamma})=\Gamma$. 
We identify ${\mathcal X}$ with $\Gamma$ via the projection $\Gamma\to{\mathcal X}$.
\par
Let $i\colon\Gamma\hookrightarrow{\mathcal X}\times C$ be the inclusion.
Let $p_{1}\colon{\mathcal X}\times C\to{\mathcal X}$ and $p_{2}\colon{\mathcal X}\times C\to C$ 
be the projections.
By the $G$-equivariances of $i$, $p_{1}$, $p_{2}$, we have the exact sequence 
of $G$-equivariant coherent sheaves on ${\mathcal X}\times C$,
\begin{equation}\label{eqn:(4.1)}
0
\longrightarrow
{\mathcal O}_{{\mathcal X}\times C}
([\Gamma]^{-1}\otimes p_{1}^{*}\xi)
@>\otimes s_{\Gamma}>>
{\mathcal O}_{{\mathcal X}\times C}(p_{1}^{*}\xi)
\longrightarrow
i_{*}{\mathcal O}_{\Gamma}(p_{1}^{*}\xi)
\longrightarrow0.
\end{equation}
\par
Let $\lambda_{G}(p_{1}^{*}\xi)$, 
$\lambda_{G}([\Gamma]^{-1}\otimes p_{1}^{*}\xi)$, 
$\lambda_{G}(\xi)$ 
be the equivariant determinants of the direct images 
$R(p_{2})_{*}{\mathcal O}_{{\mathcal X}\times C}(p_{1}^{*}\xi)$, 
$R(p_{2})_{*}{\mathcal O}_{{\mathcal X}\times C}([\Gamma]^{-1}\otimes p_{1}^{*}\xi)$, 
$R\pi_{*}{\mathcal O}_{\mathcal X}(\xi)$, respectively. 
Under the isomorphism $p_{1}^{*}\xi|_{\Gamma}\cong\xi$ 
induced from the identification $p_{1}\colon\Gamma\to{\mathcal X}$, 
the holomorphic vector bundle $\lambda_{G}$ on $C$ defined as
$$
\begin{array}{ll}
\lambda_{G}
&:=
\lambda_{G}\left([\Gamma]^{-1}\otimes p_{1}^{*}\xi\right)
\otimes
\lambda_{G}(p_{1}^{*}\xi)^{-1}
\otimes
\lambda_{G}(\xi)
=
\prod_{W\in\widehat{G}}\lambda_{W},
\\
\lambda_{W}
&:=
\lambda_{W}\left([\Gamma]^{-1}\otimes p_{1}^{*}\xi\right)
\otimes
\lambda_{W}(p_{1}^{*}\xi)^{-1}
\otimes
\lambda_{W}(\xi)
\end{array}
$$
carries the canonical holomorphic section 
$\sigma_{KM}=((\sigma_{KM})_{W})_{W\in\widehat{G}}$ such that 
$(\sigma_{KM})_{W}$ is identified with $1\in H^{0}(C,{\mathcal O}_{C})$ 
under the canonical isomorphism $\lambda_{W}\cong{\mathcal O}_{C}$;
since $\lambda_{W}$ is the determinant of the acyclic complex of coherent sheaves on $C$ 
obtained as the $W$-component of the long exact sequence of direct images associated to \eqref{eqn:(4.1)}, 
$\lambda_{W}$ is canonically isomorphic to ${\mathcal O}_{C}$ 
(\cite{Bismut95}, \cite{BismutLebeau91}, \cite{KnudsenMumford76}). 
Then $\sigma_{\rm KM}$ is admissible.
\par
Let $U\subset S$ be a relatively compact neighborhood of $0\in\Delta$ and set $U^{o}:=U\setminus\{0\}$. 
On $X=\pi^{-1}(S)$, we identify $\pi$ (resp. $d\pi$) with $s\circ\pi$ (resp. $d(s\circ\pi)$). 
Hence $\pi\in{\mathcal O}(X)$ and $d\pi\in H^{0}(X,\Omega^{1}_{X})$ in what follows.
\par
Let $h_{[\Gamma]}$ be a $G$-invariant $C^{\infty}$ Hermitian metric on $[\Gamma]$ such that
\begin{equation}\label{eqn:(4.2)}
h_{[\Gamma]}(\varsigma_{\Gamma},\varsigma_{\Gamma})(w,t)
=
\begin{cases}
\begin{array}{lcr}
|\pi(w)-t|^{2}&\hbox{if}&(w,t)\in\pi^{-1}(U)\times U,
\\
1&\hbox{if}&(w,t)\in({\mathcal X}\setminus X)\times U
\end{array}
\end{cases}
\end{equation}
and let $h_{[\Gamma]^{-1}}$ be the metric on $[\Gamma]^{-1}$ induced from $h_{[\Gamma]}$.
\par
For $g\in G$, let $\|\cdot\|_{Q,\lambda_{G}(\xi)}(g)$ be the equivariant Quillen metric on $\lambda_{G}(\xi)$ 
with respect to $h_{{\mathcal X}/C}$, $h_{\xi}$.
Let $\|\cdot\|_{Q,\lambda_{G}([\Gamma]^{-1}\otimes p_{1}^{*}\xi)}$ 
(resp. $\|\cdot\|_{Q,\lambda_{G}(p_{1}^{*}\xi)}$) be the equivariant Quillen metric 
on $\lambda_{G}([\Gamma]^{-1}\otimes p_{1}^{*}\xi)$ 
(resp. $\lambda_{G}(p_{1}^{*}\xi)$) with respect to 
$h_{\mathcal X}$, $h_{[\Gamma]^{-1}}\otimes h_{\xi}$ (resp. $h_{\mathcal X}$, $h_{\xi}$).
Let $\|\cdot\|_{Q,\lambda_{G}}$ be the equivariant Quillen metric 
on $\lambda_{G}$ defined as the tensor product of those on 
$\lambda_{G}([\Gamma]^{-1}\otimes p_{1}^{*}\xi)$,
$\lambda_{G}(p_{1}^{*}\xi)^{-1}$, 
$\lambda_{G}(\xi)$.
\par
For the germ $(\pi\colon({\mathcal X},X_{0})\to(S,0),\xi)$, we define its topological invariant 
$$
{\frak a}_{g}(X_{0},\xi)
=
\int_{E_{0}\cap\widetilde{\mathcal X}^{g}_{H}}
\widetilde{\gamma}^{*}\{
\frac{1-{\rm Td}({\mathcal H})^{-1}}{c_{1}({\mathcal H})}\}\,
q^{*}\{{\rm Td}_{g}(T{\mathcal X}){\rm ch}_{g}(\xi)\}
-
\int_{{\mathcal X}^{g}_{V}\cap X_{0}}{\rm Td}_{g}(T{\mathcal X}){\rm ch}_{g}(\xi).
$$

\begin{theorem}\label{Theorem4.1}
For $g\in G$, the following identity of functions on $S^{o}$ holds:
$$
\log\|\sigma_{KM}\|^{2}_{Q,\lambda_{G}}(g)\equiv_{\mathcal B}{\frak a}_{g}(X_{0},\xi)\,\log|s|^{2}.
$$
\end{theorem}

Since ${\mathcal X}^{1}_{H}={\mathcal X}$ and ${\mathcal X}^{1}_{V}=\emptyset$ for $g=1$, 
we get \cite[Th.\,1.1]{Yoshikawa07} by Theorem~\ref{Theorem1.1}. 

\begin{pf}
We follow \cite[Sect.\,5]{Bismut97}, \cite[Th.\,6.3]{Yoshikawa04}, \cite[Th.\,5.1]{Yoshikawa07}. 
The proof is quite parallel to that of \cite[Th.\,5.1]{Yoshikawa07}. The major differences come from the fact 
that ${\mathcal X}^{g}$ consists of the horizontal component ${\mathcal X}^{g}_{H}$ and vertical component 
${\mathcal X}^{g}_{V}$ and these two components give different contributions to the singularity of 
$\log\|\sigma_{KM}\|^{2}_{Q,\lambda_{G}}(g)$.
\newline{\bf Step 1 }
Let $[X_{s}]=[\Gamma]|_{X_{s}}$ be the holomorphic line bundle on ${\mathcal X}$ associated to 
the divisor $X_{s}$. The canonical section of $[X_{s}]$ is defined as 
$\varsigma_{s}=\varsigma_{\Gamma}|_{{\mathcal X}\times\{s\}}\in H^{0}({\mathcal X},[X_{s}])$. 
Then ${\rm div}(\varsigma_{s})=X_{s}$. Let $i_{s}\colon X_{s}\hookrightarrow{\mathcal X}$ 
be the natural embedding.
By \eqref{eqn:(4.1)}, we get the exact sequence of $G$-equivariant coherent sheaves on ${\mathcal X}$,
\begin{equation}\label{eqn:(4.3)}
0
\longrightarrow
{\mathcal O}_{\mathcal X}([X_{s}]^{-1}\otimes\xi)
@>\otimes\varsigma_{s}>>
{\mathcal O}_{\mathcal X}(\xi)
\longrightarrow
(i_{s})_{*}{\mathcal O}_{X_{s}}(\xi)
\longrightarrow
0,
\end{equation}
which induces the canonical isomorphism
$(\lambda_{G})_{s}=\lambda_{G}([X_{s}]^{-1}\otimes\xi)\otimes\lambda_{G}(\xi)^{-1}\otimes\lambda_{G}(\xi_{s})$.
\par
Set $h_{[X_{s}]}=h_{[\Gamma]}|_{{\mathcal X}\times\{s\}}$, which is a $G$-invariant Hermitian metric on $[X_{s}]$. 
Let $h_{[X_{s}]^{-1}}$ be the $G$-invariant Hermitian metric on $[X_{s}]^{-1}$ induced from $h_{[X_{s}]}$.
\par
Let $N_{s}=N_{X_{s}/{\mathcal X}}$ (resp. $N_{s}^{*}=N^{*}_{X_{s}/{\mathcal X}}$) 
be the normal (resp. conormal) bundle of $X_{s}$ in ${\mathcal X}$. 
Then $d\pi|_{X_{s}}\in H^{0}(X_{s},N_{s}^{*})$ generates $N^{*}_{s}$ for 
$s\in S^{o}$ and $d\pi|_{X_{s}}$ is $G$-invariant (cf. \cite[Eq.\,(2.2)]{Bismut95}). 
Let $a_{N^{*}_{s}}$ be the $G$-invariant Hermitian metric on $N^{*}_{s}$ 
defined by $a_{N^{*}_{s}}(d\pi|_{X_{s}},d\pi|_{X_{s}})=1$.
Let $a_{N_{s}}$ be the $G$-invariant Hermitian metric on $N_{s}$ induced from $a_{N^{*}_{s}}$. 
We have the equality $c_{1}(N_{s},a_{N_{s}})=0$ for $s\in U^{o}$.
By \cite[Proof of Th.\,5.1 Step 1]{Yoshikawa07}, the $G$-invariant metrics $h_{[X_{s}]^{-1}}\otimes h_{\xi}$ 
and $h_{\xi}$ verify assumption (A) of Bismut \cite[Def.1.5]{Bismut90} with respect to 
$a_{N_{s}}$ and $h_{\xi}|_{X_{s}}$.
\newline{\bf Step 2 }
Let ${\mathcal E}_{s}$ be the exact sequence of $G$-equivariant holomorphic vector bundles on $X_{s}$ 
defined as ${\mathcal E}_{s}\colon0\to TX_{s}\to T{\mathcal X}|_{X_{s}}\to N_{s}\to0$.
By \cite[I]{BGS88}, one has the Bott-Chern class
$\widetilde{\rm Td}_{g}({\mathcal E}_{s};h_{X_{s}},h_{\mathcal X},a_{N_{s}})\in\widetilde{A}_{X^{g}_{s}}$ 
such that
$$
dd^{c}\widetilde{\rm Td}_{g}({\mathcal E}_{s};h_{X_{s}},h_{\mathcal X},a_{N_{s}})
=
{\rm Td}_{g}(TX_{s},h_{X_{s}})\,{\rm Td}(N_{s},a_{N_{s}})|_{X^{g}_{s}}
-{\rm Td}_{g}(T{\mathcal X},h_{\mathcal X}).
$$
Here we used the triviality of the $G$-action on $N_{s}|_{X_{s}^{g}}$ to get the equality
${\rm Td}_{g}(N_{s},a_{N_{s}})={\rm Td}(N_{s},a_{N_{s}})|_{X_{s}^{g}}$.
Set $(X_{H}^{g})_{s}:={\mathcal X}_{H}^{g}\cap X_{s}$.
Applying the embedding formula of Bismut \cite{Bismut95}
(see also \cite[Th.\,5.6]{Bismut97}) to the $G$-equivariant embedding
$i_{s}\colon X_{s}\hookrightarrow{\mathcal X}$ and to the exact sequence 
\eqref{eqn:(4.3)}, we get for all $s\in U^{o}$
\begin{equation}\label{eqn:(4.4)}
\begin{aligned}
\log\|\sigma_{KM}(s)\|^{2}_{Q,\lambda_{G}}(g)
&
=
\int_{({\mathcal X}^{g}_{V}\times\{s\})\amalg({\mathcal X}^{g}_{H}\times\{s\})}
-\frac{{\rm Td}_{g}(T{\mathcal X},h_{\mathcal X})\,{\rm ch}_{g}(\xi,h_{\xi})}{{\rm Td}([\Gamma],h_{[\Gamma]})}
\log h_{[\Gamma]}(\varsigma_{\Gamma},\varsigma_{\Gamma})
\\
&
\quad
-\int_{(X^{g}_{H})_{s}}\frac{\widetilde{\rm Td}_{g}
({\mathcal E}_{s};\,h_{X_{s}},h_{\mathcal X},a_{N_{s}})\,{\rm ch}_{g}(\xi,h_{\xi})}{{\rm Td}(N_{s},a_{N_{s}})}
+C(g),
\end{aligned}
\end{equation}
where $C(g)$ is a topological constant independent of $s\in U^{o}$.
Here we used the triviality of the $G$-action on $[X_{s}]|_{X_{s}^{g}}$ and
the explicit formula for the Bott-Chern current \cite[Rem.\,3.5, especially (3.23), Th.\,3.15, Th.\,3.17]{BGS90}
to get the first term of the right hand side of \eqref{eqn:(4.4)}.
Substituting \eqref{eqn:(4.2)} and $c_{1}(N_{s},a_{N_{s}})=0$ into \eqref{eqn:(4.4)}, we get for $s\in U^{o}$
\begin{equation}\label{eqn:(4.5)}
\begin{aligned}
\,
&
\log\|\sigma_{KM}(s)\|^{2}_{Q,\lambda_{G}}(g)
\equiv_{\mathcal B}
-
\int_{{\mathcal X}^{g}_{H}\times\{s\}}{\rm Td}_{g}(T{\mathcal X},h_{\mathcal X})\,{\rm ch}_{g}(\xi,h_{\xi})
\log|\pi-s|^{2}
\\
&
-\int_{({\mathcal X}^{g}_{V}\cap X_{0})\times\{s\}}{\rm Td}_{g}(T{\mathcal X})\,{\rm ch}_{g}(\xi)\log|s|^{2}
-\int_{(X^{g}_{H})_{s}}
\widetilde{\rm Td}_{g}({\mathcal E}_{s};\,h_{X_{s}},h_{\mathcal X},a_{N_{s}})\,{\rm ch}_{g}(\xi,h_{\xi})
\\
&\equiv_{\mathcal B}
-\{\int_{{\mathcal X}^{g}_{V}\cap X_{0}}{\rm Td}_{g}(T{\mathcal X})\,{\rm ch}_{g}(\xi)\}\,\log|s|^{2}
-
\int_{(X^{g}_{H})_{s}}
\widetilde{\rm Td}_{g}({\mathcal E}_{s};\,h_{X_{s}},h_{\mathcal X},a_{N_{s}})\,{\rm ch}_{g}(\xi,h_{\xi}),
\end{aligned}
\end{equation}
where we used the equality
$h_{[\Gamma]}(\varsigma_{\Gamma},\varsigma_{\Gamma})|_{X_{0}\times\{s\}}=|s|^{2}$
and the fact ${\mathcal X}_{V}^{g}\cap\pi^{-1}(U)={\mathcal X}_{V}^{g}\cap X_{0}$
to get the first equality and \cite[Th.\,9.1]{Yoshikawa07} to get the second equality. 
\newline{\bf Step 3 }
Let $h_{N_{s}}$ be the Hermitian metric on $N_{s}$ induced from $h_{\mathcal X}$ 
by the $C^{\infty}$ isomorphism $N_{s}\cong(TX_{s})^{\perp}$. 
Let $\widetilde{\rm Td}(N_{s};\,a_{N_{s}},h_{N_{s}})\in\widetilde{A}_{X_{s}}$ 
be the Bott--Chern secondary class such that
$$
dd^{c}\widetilde{\rm Td}(N_{s};a_{N_{s}},h_{N_{s}})
=
{\rm Td}(N_{s},a_{N_{s}})-{\rm Td}(N_{s},h_{N_{s}}).
$$
Since $G$ acts trivially on $N_{s}|_{X_{s}^{g}}$, we deduce from
\cite[I, Props.\,1.3.2 and 1.3.4]{GilletSoule90} that
\begin{equation}\label{eqn:(4.6)}
\begin{aligned}
\,&
\widetilde{\rm Td}_{g}({\mathcal E}_{s};\,h_{X_{s}},h_{\mathcal X},a_{N_{s}})=
\\
&
\widetilde{\rm Td}_{g}({\mathcal E}_{s};\,h_{X_{s}},h_{\mathcal X},h_{N_{s}})
+
{\rm Td}_{g}(TX_{s},h_{X_{s}})\,\widetilde{\rm Td}(N_{s};\,a_{N_{s}},h_{N_{s}})=
\\
&
\widetilde{\rm Td}_{g}({\mathcal E}_{s};\,h_{X_{s}},h_{\mathcal X},h_{N_{s}})
+
\left.\gamma^{*}{\rm Td}_{g}({\mathcal U},h_{\mathcal U})\,
\nu^{*}\{
\frac{1-{\rm Td}(-c_{1}({\mathcal L},h_{\mathcal L}))}
{-c_{1}({\mathcal L},h_{\mathcal L})}
\}\,\log\|d\pi\|^{2}\right|_{(X^{g}_{H})_{s}}.
\end{aligned}
\end{equation}
Here we used \cite[Eq.\,(13)]{Yoshikawa07} and the relation
$(TX_{s},h_{X_{s}})=\gamma^{*}({\mathcal U},h_{\mathcal U})|_{X_{s}}$ to get the second equality.
Since
$$
({\mathcal E}_{s},h_{X_{s}},h_{\mathcal X},h_{N_{s}})
=
\gamma^{*}({\mathcal S}^{\lor},h_{\mathcal U},(\varPi^{\lor})^{*}h_{\mathcal X},h_{\mathcal H})|_{X_{s}},
$$
we get by the functorial property of the Bott--Chern class \cite[I]{BGS88}
\begin{equation}\label{eqn:(4.7)}
\widetilde{\rm Td}_{g}({\mathcal E}_{s};\,h_{X_{s}},h_{\mathcal X},h_{N_{s}})
=
\gamma^{*}\widetilde{\rm Td}_{g}
({\mathcal S}^{\lor};\,h_{\mathcal U},(\varPi^{\lor})^{*}h_{\mathcal X},h_{\mathcal H})
|_{(X^{g}_{H})_{s}}.
\end{equation}
Substituting \eqref{eqn:(4.7)} into \eqref{eqn:(4.6)}, we get
\begin{equation}\label{eqn:(4.8)}
\begin{aligned}
\widetilde{\rm Td}_{g}({\mathcal E}_{s};\,h_{X_{s}},h_{\mathcal X},a_{N_{s}})
&
=
\gamma^{*}\widetilde{\rm Td}_{g}
({\mathcal S}^{\lor};\,h_{\mathcal U},(\varPi^{\lor})^{*}h_{\mathcal X},h_{\mathcal H})
|_{(X^{g}_{H})_{s}}+
\\
&\quad
\gamma^{*}{\rm Td}_{g}({\mathcal U},h_{\mathcal U})\,
\nu^{*}\{
\frac{1-{\rm Td}(-c_{1}({\mathcal L},h_{\mathcal L}))}{-c_{1}({\mathcal L},h_{\mathcal L})}\}\,
\log\|d\pi\|^{2}|_{(X^{g}_{H})_{s}}.
\end{aligned}
\end{equation}
Substituting \eqref{eqn:(4.8)} into \eqref{eqn:(4.5)}, we get by the same argument 
as in \cite[p.74 l.1-l.13]{Yoshikawa07}
\begin{equation}\label{eqn:(4.9)}
\begin{aligned}
\,&
\log\|\sigma_{KM}\|^{2}_{Q,\lambda_{G}}(g)
\equiv_{\mathcal B}
-\{\int_{{\mathcal X}^{g}_{V}\cap X_{0}}{\rm Td}_{g}(T{\mathcal X})\,{\rm ch}_{g}(\xi)\}
\,\log|s|^{2}
\\
&\qquad
-(\pi|_{{\mathcal X}^{g}_{H}})_{*}
\left[
\gamma^{*}\widetilde{\rm Td}_{g}
({\mathcal S}^{\lor};\,h_{\mathcal U},(\varPi^{\lor})^{*}h_{\mathcal X},h_{\mathcal H})\,
{\rm ch}_{g}(\xi,h_{\xi})
\right]^{(0)}
\\
&\qquad
-(\pi|_{{\mathcal X}^{g}_{H}})_{*}
\left[
\gamma^{*}{\rm Td}_{g}({\mathcal U},h_{\mathcal U})\,
\nu^{*}\{
\frac{1-{\rm Td}(-c_{1}({\mathcal L},h_{\mathcal L}))}{-c_{1}({\mathcal L},h_{\mathcal L})}\}\,
{\rm ch}_{g}(\xi,h_{\xi})\,
\log\|d\pi\|^{2}
\right]^{(0)}
\\
&
\equiv_{\mathcal B}
-\{\int_{{\mathcal X}^{g}_{V}\cap X_{0}}{\rm Td}_{g}(T{\mathcal X})\,{\rm ch}_{g}(\xi)\}
\,\log|s|^{2}
\\
&
+
(\widetilde{\pi}|_{{\mathcal X}^{g}_{H}})_{*}
\left[
\widetilde{\gamma}^{*}
{\rm Td}_{g}({\mathcal U},h_{\mathcal U})\,
\widetilde{\nu}^{*}\{
\frac{{\rm Td}(-c_{1}({\mathcal L},h_{\mathcal L}))-1}{-c_{1}({\mathcal L},h_{\mathcal L})}\}\,
q^{*}{\rm ch}_{g}(\xi,h_{\xi})\,(q^{*}\log\|d\pi\|^{2})
\right]^{(0)}
\end{aligned}
\end{equation}
where we set $\widetilde{\pi}:=\pi\circ q$.
By \cite[Cor.\,4.6]{Yoshikawa07} applied to the last line of \eqref{eqn:(4.9)}, we get
$$
\begin{aligned}
\,&
\log\|\sigma_{KM}\|^{2}_{Q,\lambda_{G}}(g)
\equiv_{\mathcal B}
\\
&
[\int_{\widetilde{\mathcal X}^{g}_{H}\cap E_{0}}
\widetilde{\gamma}^{*}
\{{\rm Td}_{g}({\mathcal U})\,\frac{{\rm Td}({\mathcal H})-1}{c_{1}({\mathcal H})}\}q^{*}{\rm ch}_{g}(\xi)
-\int_{{\mathcal X}^{g}_{V}\cap X_{0}}{\rm Td}_{g}(T{\mathcal X})\,
{\rm ch}_{g}(\xi)]\,\log|s|^{2}
\\
&=
[\int_{\widetilde{\mathcal X}^{g}_{H}\cap E_{0}}
\widetilde{\gamma}^{*}
\{\frac{1-{\rm Td}({\mathcal H})^{-1}}{c_{1}({\mathcal H})}\}
q^{*}\{{\rm Td}_{g}(T{\mathcal X}){\rm ch}_{g}(\xi)\}]\,\log|s|^{2}
\\
&\qquad
-\{\int_{{\mathcal X}^{g}_{V}\cap X_{0}}{\rm Td}_{g}(T{\mathcal X})\,{\rm ch}_{g}(\xi)\}\,\log|s|^{2}
=
{\frak a}_{g}(X_{0},\xi)\,\log|s|^{2}.
\end{aligned}
$$
Here the first equality follows from the identity
${\rm Td}_{g}({\mathcal U}){\rm Td}({\mathcal H})=(\varPi^{\lor})^{*}{\rm Td}_{g}(T{\mathcal X})$,
which is deduced from the exact sequence
$0\to{\mathcal U}\to(\varPi^{\lor})^{*}T{\mathcal X}\to{\mathcal H}\to0$
on ${\bf P}(T{\mathcal X})^{\lor}$. 
This completes the proof.
\end{pf}

\begin{theorem}\label{Theorem1.1}
For $g\in G$, the following identity holds:
$$
\log\|\sigma\|^{2}_{Q,\lambda_{G}(\xi)}(g)\equiv_{{\mathcal B}}{\frak a}_{g}(X_{0},\xi)\log|s|^{2}.
$$
\end{theorem}

\begin{pf}
There exist admissible holomorphic sections
$$
\alpha=(\alpha_{W})_{W\in\widehat{G}}\in\Gamma(U,\lambda_{G}(p_{1}^{*}\xi)),
\qquad
\beta=(\beta_{W})_{W\in\widehat{G}}\in\Gamma(U,\lambda_{G}([\Gamma]^{-1}\otimes p_{1}^{*}\xi))
$$
such that $\sigma_{KM}=\beta\otimes\alpha^{-1}\otimes\sigma$ on $S$, i.e., 
$(\sigma_{KM})_{W}=\beta_{W}\otimes\alpha_{W}^{-1}\otimes\sigma_{W}$ for all $W\in\widehat{G}$. Then
$$
\begin{aligned}
\log\|\sigma\|^{2}_{Q,\lambda_{G}(\xi)}(g)
&=
\log\|\sigma_{KM}\|^{2}_{Q,\lambda_{G}}(g)
+
\log\|\alpha\|^{2}_{Q,\lambda_{G}(p_{1}^{*}\xi)}(g)
-
\log\|\beta\|^{2}_{Q,\lambda_{G}([\Gamma]\otimes p_{1}^{*}\xi)}(g)
\\
&\equiv_{\mathcal B}
{\frak a}_{g}(X_{0},\xi)\,\log|s|^{2}
\end{aligned}
$$
by Theorems~\ref{Theorem3.1} and \ref{Theorem4.1}. This proves the theorem.
\end{pf}

\begin{corollary}\label{Corollary4.2}
The following equation of $(1,1)$-currents on $S$ holds:
$$
-dd^{c}\log\|\sigma\|_{Q,\lambda_{G}(\xi)}^{2}(g)
=
\pi_{*}
\left[{\rm Td}_{g}(T{\mathcal X}/C,h_{{\mathcal X}/C})\,{\rm ch}_{g}(\xi,h_{\xi})\right]^{(1,1)}
-
{\frak a}_{g}(X_{0},\xi)\,\delta_{0}.
$$
\end{corollary}

\begin{pf}
The result follows from the curvature formula \eqref{eqn:(3.6)} and Theorem~\ref{Theorem1.1}.
\end{pf}

\subsection
{The case of adjoint bundles twisted by the relative canonical bundle}
\label{subsect:4.2}
\par
We set $\Omega_{{\mathcal X}/C}^{1}:=\Omega_{\mathcal X}^{1}/\pi^{*}\Omega_{C}^{1}$ 
and $\Omega_{{\mathcal X}/C}^{q}:=\bigwedge^{q}\Omega_{{\mathcal X}/C}^{1}$.
Let $\omega_{\mathcal X}:=\Omega_{\mathcal X}^{n+1}$ be the canonical line bundle of ${\mathcal X}$ 
and let $\omega_{{\mathcal X}/C}:=\Omega_{\mathcal X}^{n+1}\otimes(\pi^{*}\Omega_{C}^{1})^{-1}$ 
be the relative canonical line bundle of $\pi\colon{\mathcal X}\to C$, which are identified with the dualizing
sheaf of ${\mathcal X}$ and the relative dualizing sheaf of $\pi\colon{\mathcal X}\to C$ respectively.
On ${\mathcal X}\setminus\varSigma_{\pi}$, there is a canonical isomorphism
$\Omega^{n}_{{\mathcal X}/C}\cong\omega_{{\mathcal X}/C}$ induced by the short exact sequence
$0\to\pi^{*}\Omega_{C}^{1}\to\Omega_{\mathcal X}^{1}\to\Omega_{{\mathcal X}/C}^{1}\to0$.
The holomorphic vector bundle $\Omega^{q}_{{\mathcal X}/C}$ on ${\mathcal X}\setminus\varSigma_{\pi}$
is equipped with the Hermitian metric induced from $h_{{\mathcal X}/C}$.
Since $\pi^{*}\Omega_{C}^{1}\subset\Omega_{\mathcal X}^{1}$,
$\omega_{\mathcal X}$ and $\omega_{{\mathcal X}/C}$ are equipped with the Hermitian metrics 
$h_{\omega_{\mathcal X}}$ and $h_{\omega_{{\mathcal X}/C}}$ induced from $h_{\mathcal X}$, respectively. 
Then the canonical isomorphism $\Omega^{n}_{{\mathcal X}/C}\cong\omega_{{\mathcal X}/C}$ is an isometry. 
\par
For $g\in G$, let $\|\cdot\|_{Q,\lambda_{G}(\omega_{{\mathcal X}/C}(\xi))}(g)$ be the equivariant 
Quillen metric on $\lambda_{G}(\omega_{{\mathcal X}/C}(\xi))$ with respect to 
$h_{{\mathcal X}/C}$, $h_{\xi}$, $h_{\omega_{{\mathcal X}/C}}$.
Let $\varsigma$ be a nowhere vanishing holomorphic section of $\lambda(\omega_{{\mathcal X}/C}(\xi))$
on $S$.

\begin{theorem}\label{Theorem4.3}
For $g\in G$, the following identity of functions on $S^{o}$ holds:
$$
\log\|\varsigma\|^{2}_{Q,\lambda_{G}(\omega_{{\mathcal X}/C}(\xi))}(g)
\equiv_{{\mathcal B}}
\alpha_{g}(X_{0},\omega_{X/S}(\xi))\,\log|s|^{2}.
$$
\end{theorem}

\begin{pf}
By Theorem~\ref{Theorem1.1} applied to $\omega_{\mathcal X}(\xi)$, we get on $S$
\begin{equation}\label{eqn:(4.10)}
\log\|\varsigma\|^{2}_{Q,\lambda_{G}(\omega_{\mathcal X}(\xi))}(g)
\equiv_{{\mathcal B}}
{\frak a}_{g}(X_{0},\omega_{\mathcal X}(\xi))\,\log|s|^{2}.
\end{equation}
Since $\pi^{*}\Omega^{1}_{C}$ is generated by $d\pi$ on $X\setminus\varSigma_{\pi}$, 
we get an identification 
$(\omega_{{\mathcal X}/C},h_{\omega_{{\mathcal X}/C}})=
(\omega_{\mathcal X},\|d\pi\|^{-2}h_{\omega_{\mathcal X}})$ on $X$. 
We set $h_{\omega_{\mathcal X}(\xi)}:=h_{\omega_{\mathcal X}}\otimes h_{\xi}$
and $h_{\omega_{{\mathcal X}/C}(\xi)}:=h_{\omega_{{\mathcal X}/C}}\otimes h_{\xi}$.
By the anomaly formula \cite[Th.\,2.5]{Bismut95}, \cite[I, Th.\,0.3]{BGS88}, we get
\begin{equation}\label{eqn:(4.11)}
\begin{aligned}
\,
&
\log\|\varsigma\|^{2}_{Q,\lambda_{G}(\omega_{{\mathcal X}/C}(\xi))}(g)
=
\log\|\varsigma\|^{2}_{Q,\lambda_{G}(\omega_{\mathcal X}(\xi))}(g)
+
\log\frac
{\|\cdot\|^{2}_{Q,\lambda_{G}(\omega_{{\mathcal X}/C}(\xi))}(g)}
{\|\cdot\|^{2}_{Q,\lambda_{G}(\omega_{\mathcal X}(\xi))}(g)}
\\
&=
\log\|\varsigma\|^{2}_{Q,\lambda_{G}(\omega_{\mathcal X}(\xi))}(g)
+
\pi_{*}({\rm Td}_{g}(T{\mathcal X}/C,h_{{\mathcal X}/C})
\widetilde{{\rm ch}}_{g}
(\omega_{\mathcal X}(\xi);h_{\omega_{\mathcal X}(\xi)},h_{\omega_{{\mathcal X}/C}(\xi)})
)^{(0)}.
\end{aligned}
\end{equation}
Here 
$\widetilde{{\rm ch}}_{g}(\omega_{{\mathcal X}}(\xi);
h_{\omega_{\mathcal X}(\xi)},h_{\omega_{{\mathcal X}/C}(\xi)})$
is the Bott-Chern class such that
$$
dd^{c}
\widetilde{{\rm ch}}_{g}(\omega_{\mathcal X}(\xi);
h_{\omega_{\mathcal X}(\xi)},h_{\omega_{{\mathcal X}/C}(\xi)})
=
{\rm ch}_{g}(\omega_{\mathcal X}(\xi),h_{\omega_{\mathcal X}(\xi)})
-
{\rm ch}_{g}(\omega_{\mathcal X}(\xi),h_{\omega_{{\mathcal X}/C}(\xi)})
$$
Since
$(\omega_{{\mathcal X}/C}(\xi),h_{\omega_{{\mathcal X}/C}(\xi)})
=(\omega_{\mathcal X}(\xi),\|d\pi\|^{-2}h_{\omega_{\mathcal X}(\xi)})$ and
$-dd^{c}\log\|d\pi\|^{2}=\gamma^{*}c_{1}({\mathcal L},h_{\mathcal L})$, 
we get  by \cite[I, (1.2.5.1), (1.3.1.2]{GilletSoule90} (see also \cite[Eqs.\,(3.7), (5.5)]{FLY08})
\begin{equation}\label{eqn:(4.12)}
\begin{aligned}
\,&
\widetilde{{\rm ch}}_{g}
(\omega_{\mathcal X}(\xi);h_{\omega_{\mathcal X}(\xi)},h_{\omega_{{\mathcal X}/C}(\xi)})
=
\widetilde{{\rm ch}}_{g}(\omega_{\mathcal X};h_{\omega_{\mathcal X}},\|d\pi\|^{-2}h_{\omega_{\mathcal X}})\,
{\rm ch}_{g}(\xi,h_{\xi})
\\
&=
{\rm ch}_{g}(\omega_{\mathcal X},h_{\omega_{\mathcal X}})
\frac{e^{-\gamma^{*}c_{1}({\mathcal L},h_{\mathcal L})}-1}{-\gamma^{*}c_{1}({\mathcal L},h_{\mathcal L})}\,
(-\log\|d\pi\|^{2})\wedge
{\rm ch}_{g}(\xi,h_{\xi})
\\
&=
-{\rm ch}_{g}(\omega_{\mathcal X}(\xi),h_{\omega_{\mathcal X}(\xi)})
\frac{1-e^{-\gamma^{*}c_{1}({\mathcal L},h_{\mathcal L})}}{\gamma^{*}c_{1}({\mathcal L},h_{\mathcal L})}\,
\log\|d\pi\|^{2}.
\end{aligned}
\end{equation}
Since ${\rm Td}_{g}(T{\mathcal X}/C, h_{{\mathcal X}/C})=\gamma^{*}{\rm Td}_{g}({\mathcal U})$, we get
by \eqref{eqn:(4.11)} and \cite[Cor.\,4.6]{Yoshikawa07}
\begin{equation}\label{eqn:(4.13)}
\begin{aligned}
\,&
\log\frac
{\|\cdot\|^{2}_{Q,\lambda_{G}(\omega_{{\mathcal X}/C}(\xi))}(g)}
{\|\cdot\|^{2}_{Q,\lambda_{G}(\omega_{\mathcal X}(\xi))}(g)}
\\
&=
-\pi_{*}\{
{\rm Td}_{g}(T{\mathcal X}/C, h_{{\mathcal X}/C})
{\rm ch}_{g}(\omega_{\mathcal X}(\xi),h_{\omega_{\mathcal X}}\otimes h_{\xi})
\frac{1-e^{-\gamma^{*}c_{1}({\mathcal L},h_{\mathcal L})}}{\gamma^{*}c_{1}({\mathcal L},h_{\mathcal L})}\,
\log\|d\pi\|^{2}
\}^{(0)}
\\
&\equiv_{\mathcal B}
-\{
\int_{E_{0}\cap\widetilde{\mathcal X}_{H}^{g}}
\widetilde{\gamma}^{*}{\rm Td}_{g}({\mathcal U})q^{*}{\rm ch}_{g}(\omega_{\mathcal X}(\xi))
\frac{1-e^{-\widetilde{\gamma}^{*}c_{1}({\mathcal L})}}{\widetilde{\gamma}^{*}c_{1}({\mathcal L})}\,
\}
\log|s|^{2}
\\
&\equiv_{\mathcal B}
-[
\int_{E_{0}\cap\widetilde{\mathcal X}_{H}^{g}}
\widetilde{\gamma}^{*}\{
{\rm Td}_{g}({\mathcal U})\frac{e^{c_{1}({\mathcal H})}-1}{c_{1}({\mathcal H})}
\}
q^{*}{\rm ch}_{g}(\omega_{\mathcal X}(\xi))
]
\log|s|^{2}.
\end{aligned}
\end{equation}
By \eqref{eqn:(4.10)} and \eqref{eqn:(4.13)}, we get
\begin{equation}\label{eqn:(4.14)}
\begin{aligned}
\,&
\log\|\varsigma\|^{2}_{Q,\lambda_{G}(\omega_{{\mathcal X}/C}(\xi))}(g)
\equiv_{{\mathcal B}}
\\
&
[
{\frak a}_{g}(\pi,X_{0},\omega_{\mathcal X}(\xi))
-
\int_{E_{0}\cap\widetilde{\mathcal X}_{H}^{g}}
\widetilde{\gamma}^{*}
\{{\rm Td}_{g}({\mathcal U})\frac{e^{c_{1}({\mathcal H})}-1}{c_{1}({\mathcal H})}\}
q^{*}{\rm ch}_{g}(\omega_{\mathcal X}(\xi))
]\,
\log|s|^{2}
\\
&=
[
\int_{E_{0}\cap\widetilde{\mathcal X}_{H}^{g}}\widetilde{\gamma}^{*}
\{{\rm Td}_{g}({\mathcal U})
\left(
\frac{{\rm Td}({\mathcal H})-1}{c_{1}({\mathcal H})}
-
\frac{e^{c_{1}({\mathcal H})}-1}{c_{1}({\mathcal H})}
\right)\}
q^{*}{\rm ch}_{g}(\omega_{\mathcal X}(\xi))
]\,
\log|s|^{2}
\\
&\quad
-\{\int_{{\mathcal X}^{g}_{V}\cap X_{0}}{\rm Td}_{g}(T{\mathcal X})\,{\rm ch}_{g}(\omega_{\mathcal X}(\xi))\}\,
\log|s|^{2}
\\
&=
[
\int_{E_{0}\cap\widetilde{\mathcal X}_{H}^{g}}
\widetilde{\gamma}^{*}
\{{\rm Td}_{g}({\mathcal U}){\rm Td}({\mathcal H})
\left(
\frac{1-{\rm Td}({\mathcal H})^{-1}}{c_{1}({\mathcal H})}
-
\frac{e^{c_{1}({\mathcal H})}}{{\rm Td}({\mathcal H})^{2}}
\right)\}
q^{*}{\rm ch}_{g}(\omega_{\mathcal X}(\xi))
]\,
\log|s|^{2}
\\
&\quad
-\{
\int_{{\mathcal X}^{g}_{V}\cap X_{0}}{\rm Td}_{g}(T{\mathcal X})\,{\rm ch}_{g}(\omega_{\mathcal X}(\xi))
\}
\log|s|^{2}
\\
&=
[
\int_{E_{0}\cap\widetilde{\mathcal X}_{H}^{g}}
\widetilde{\gamma}^{*}\left(
\frac{{\rm Td}({\mathcal H}^{\lor})^{-1}-1}{c_{1}({\mathcal H}^{\lor})}
\right)
q^{*}\{{\rm Td}_{g}(T{\mathcal X}){\rm ch}_{g}(\omega_{\mathcal X}(\xi))\}
\}
\,\log|s|^{2}
\\
&\quad
-\{
\int_{{\mathcal X}^{g}_{V}\cap X_{0}}{\rm Td}_{g}(T{\mathcal X})\,{\rm ch}_{g}(\omega_{\mathcal X}(\xi))
]\,
\log|s|^{2}
=
\alpha_{g}(X_{0},\omega_{{\mathcal X}/C}(\xi))\,\log|s|^{2}.
\end{aligned}
\end{equation}
Here the fourth equality follows from the identities
${\rm Td}_{g}({\mathcal U}){\rm Td}({\mathcal H})=(\varPi^{\lor})^{*}{\rm Td}_{g}(T{\mathcal X})$ and
$$
\frac{1-{\rm Td}(x)^{-1}}{x}-e^{x}{\rm Td}(x)^{-2}=\frac{{\rm Td}(-x)^{-1}-1}{-x},
\qquad
{\rm Td}(x)=x/(1-e^{-x}).
$$
This completes the proof.
\end{pf}

\subsection
{Compatibility with the Serre duality}
\label{subsect:4.3}
We check the compatibility of Theorems~\ref{Theorem1.1} and \ref{Theorem4.3} with the Serre duality.
There exists an exact sequence of $G$-equivariant holomorphic vector bundles
$0\to\xi\to \xi_{0}\to\cdots\to \xi_{m}\to 0$ on ${\mathcal X}$ such that 
$H^{q}(X_{s},\xi_{i}|_{X_{s}})=0$ for all $i\geq0$, $q>0$, $s\in C$.
We get the corresponding long exact sequence of $G$-equivariant holomorphic vector bundles
$0\to \omega_{{\mathcal X}/C}(\xi_{m}^{\lor})\to\cdots\to \omega_{{\mathcal X}/C}(\xi_{0}^{\lor})\to 
\omega_{{\mathcal X}/C}(\xi^{\lor})\to0$ on ${\mathcal X}$.
Since $R^{q}\pi_{*}\omega_{{\mathcal X}/C}(\xi_{i}^{\lor})=0$ for $q\not=n$ by the fiberwise Serre duality, 
we get for all $q\geq0$
$$
R^{q}\pi_{*}\omega_{{\mathcal X}/C}(\xi^{\lor})
=
\frac{\ker\{R^{n}\pi_{*}\omega_{{\mathcal X}/C}(\xi_{n-q}^{\lor})\to 
R^{n}\pi_{*}\omega_{{\mathcal X}/C}(\xi_{n-q-1}^{\lor})\}}
{{\rm Im}\{R^{n}\pi_{*}\omega_{{\mathcal X}/C}(\xi_{n-q+1}^{\lor})\to 
R^{n}\pi_{*}\omega_{{\mathcal X}/C}(\xi_{n-q}^{\lor})\}}.
$$
Hence
$\lambda_{W}(\omega_{{\mathcal X}/C}(\xi^{\lor}))^{(-1)^{n}}\cong
\bigotimes_{i\geq0}\lambda_{W}(R^{n}\pi_{*}\omega_{{\mathcal X}/C}(\xi_{i}^{\lor}))^{(-1)^{i}}$ 
for $W\in\widehat{G}$.
Since
$\lambda_{W}(\xi)^{\lor}
\cong
\bigotimes_{i\geq0}(\lambda_{W}(\pi_{*}\xi_{i})^{\lor})^{(-1)^{i}}
\cong
\bigotimes_{i\geq0}\lambda_{W}(R^{n}\pi_{*}\omega_{{\mathcal X}/C}(\xi_{i}^{\lor}))^{(-1)^{i}}$,
there is a canonical isomorphism
$\lambda_{W}(\xi)\otimes\lambda_{W}(\omega_{{\mathcal X}/C}(\xi^{\lor}))^{(-1)^{n}}\cong{\mathcal O}_{C}$.
Let ${\bf 1}_{W}$ be the canonical section of 
$\lambda_{W}(\xi)\otimes\lambda_{W}(\omega_{{\mathcal X}/C}(\xi^{\lor}))^{(-1)^{n}}$ corresponding to 
$1\in H^{0}(C,{\mathcal O}_{C})$. Then ${\bf 1}=({\bf 1}_{W})_{W\in\widehat{G}}$ is the canonical
section of $\lambda_{G}(\xi)\otimes\lambda_{G}(\omega_{{\mathcal X}/C}(\xi^{\lor}))^{(-1)^{n}}$.
By the same argument as in \cite[p.27]{GilletSoule91}, we get
$\log\|{\bf 1}\|_{Q,\lambda_{G}(\xi)\otimes\lambda_{G}(\omega_{{\mathcal X}/C}(\xi^{\lor}))^{(-1)^{n}}}(g)=0$
as a function on $C^{o}$, which, together with Theorems~\ref{Theorem1.1} and \ref{Theorem4.3}, 
implies the equality
\begin{equation}\label{eqn:(4.15)}
{\frak a}_{g}(X_{0},\xi)+(-1)^{n}\alpha_{g}(X_{0},\omega_{X/S}(\xi^{\lor}))=0.
\end{equation}
\par
We check \eqref{eqn:(4.15)}. Let $T{\mathcal X}|_{{\mathcal X}^{g}}=\oplus_{j}E(\theta_{j})$ 
be the decomposition into the eigenbundles with respect to the $g$-action.
Let $\lambda_{0},\ldots,\lambda_{n}$ be the Chern roots of $T{\mathcal X}|_{{\mathcal X}^{g}}$.
We may assume that the element $g$ acts as the multiplication by $e^{\theta_{i}}$ on the line bundle
corresponding to $\lambda_{i}$. 
Since
$e^{-\lambda_{i}}[\lambda_{i}/(1-e^{-\lambda_{i}})]=(-\lambda_{i}/1-e^{-(-\lambda_{i})})$
and
$e^{-\lambda_{i}+\theta_{i}}[1/(1-e^{-\lambda_{i}+\theta_{i}})]=-(1/1-e^{-(-\lambda_{i}+\theta_{i})})$, 
we get
$$
{\rm Td}_{g}(T{\mathcal X}){\rm ch}_{g}(\omega_{\mathcal X}\otimes\xi^{\lor})|_{{\mathcal X}^{g}_{i}}
=
(-1)^{n+1-\dim{\mathcal X}^{g}_{i}}
{\rm Td}_{g}(T{\mathcal X}^{\lor}){\rm ch}_{g}(\xi^{\lor})|_{{\mathcal X}^{g}_{i}}
$$
for every connected component ${\mathcal X}^{g}_{i}$ of ${\mathcal X}^{g}$. 
Set $d_{i}:=\dim{\mathcal X}^{g}_{i}$. Then
\begin{equation}\label{eqn:(4.16)}
\begin{aligned}
\,&
(-1)^{n}\alpha_{g}(X_{0},\omega_{{\mathcal X}/C}(\xi^{\lor}))
\\
&=
(-1)^{n}\sum_{i\in I}
\int_{E_{0}\cap(\widetilde{\mathcal X}^{g}_{H})_{i}}
\widetilde{\gamma}^{*}\{
\frac{{\rm Td}({\mathcal H}^{\lor})^{-1}-1}{c_{1}({\mathcal H}^{\lor})}\}\,
q^{*}\{{\rm Td}_{g}(T{\mathcal X}){\rm ch}_{g}(\omega_{\mathcal X}\otimes\xi^{\lor})\}
\\
&\quad
-(-1)^{n}\sum_{j\in J}
\int_{({\mathcal X}^{g}_{V})_{j}\cap X_{0}}
{\rm Td}_{g}(T{\mathcal X}){\rm ch}_{g}(\omega_{\mathcal X}\otimes\xi^{\lor})
\\
&=
\sum_{i\in I}(-1)^{d_{i}+1}\int_{E_{0}\cap(\widetilde{\mathcal X}^{g}_{H})_{i}}
\widetilde{\gamma}^{*}\{
\frac{{\rm Td}({\mathcal H}^{\lor})^{-1}-1}{c_{1}({\mathcal H}^{\lor})}\}\,
q^{*}\{{\rm Td}_{g}(T{\mathcal X}^{\lor}){\rm ch}_{g}(\xi^{\lor})\}
\\
&\quad
-\sum_{j\in J}(-1)^{d_{j}+1}
\int_{({\mathcal X}^{g}_{V})_{j}\cap X_{0}}{\rm Td}_{g}(T{\mathcal X}^{\lor}){\rm ch}_{g}(\xi^{\lor})
\\
&=
\sum_{i\in I}(-1)^{d_{i}+1}\int_{E_{0}\cap(\widetilde{\mathcal X}^{g}_{H})_{i}}
\sum_{e\geq0}(-1)^{e}
\left[
\widetilde{\gamma}^{*}\{
\frac{{\rm Td}({\mathcal H})^{-1}-1}{c_{1}({\mathcal H})}\}\,
q^{*}\{{\rm Td}_{g}(T{\mathcal X}){\rm ch}_{g}(\xi)\}
\right]^{(2e)}
\\
&\quad
-\sum_{j\in J}(-1)^{d_{j}+1}
\int_{({\mathcal X}^{g}_{V})_{j}\cap X_{0}}\sum_{e\geq0}(-1)^{e}
[{\rm Td}_{g}(T{\mathcal X}){\rm ch}_{g}(\xi)]^{(2e)},
\end{aligned}
\end{equation}
where ${\mathcal X}^{g}_{H}=\amalg_{i\in I}({\mathcal X}^{g}_{H})_{i}$ and 
${\mathcal X}^{g}_{V}=\amalg_{j\in J}({\mathcal X}^{g}_{V})_{j}$
are the decompositions into the connected components. 
Since  ${\mathcal X}^{g}_{H}$ intersects $X_{0}$ properly, we get
$\dim E_{0}\cap(\widetilde{{\mathcal X}}^{g}_{H})_{i}=\dim(\widetilde{\mathcal X}^{g}_{H})_{i}-1$.
On the other hand, since ${\mathcal X}^{g}_{V}$ is contained in the singular fiber of $\pi$, we get
$\dim ({\mathcal X}^{g}_{V})_{j}\cap X_{0}=\dim ({\mathcal X}^{g}_{V})_{j}$ 
when $({\mathcal X}^{g}_{V})_{j}\cap X_{0}\not=\emptyset$. By \eqref{eqn:(4.16)}, we get \eqref{eqn:(4.15)}
$$
\begin{aligned}
(-1)^{n}\alpha_{g}(X_{0},\omega_{{\mathcal X}/C}(\xi^{\lor}))
&=
\sum_{i\in I}\int_{E_{0}\cap(\widetilde{\mathcal X}^{g}_{H})_{i}}
\widetilde{\gamma}^{*}\{
\frac{{\rm Td}({\mathcal H})^{-1}-1}{c_{1}({\mathcal H})}\}\,
q^{*}\{{\rm Td}_{g}(T{\mathcal X}){\rm ch}_{g}(\xi)\}
\\
&\quad
+
\sum_{j\in J}\int_{({\mathcal X}^{g}_{V})_{j}\cap X_{0}}{\rm Td}_{g}(T{\mathcal X}){\rm ch}_{g}(\xi)
=-{\frak a}_{g}(X_{0},\xi).
\end{aligned}
$$

\section
{Nakano semi-positive vector bundles}
\label{sect:5}
\par

\subsection
{Semi-positivity and semi-negativity of vector bundles}
\label{subsect:5.1}
\par
Let $M$ be a connected complex manifold of dimension $n+1$ and 
let $E\to M$ be a holomorphic vector bundle of rank $r$ equipped with 
a Hermitian metric $h_{E}$.
Let $R^{E}=(\nabla^{E})^{2}$ be the curvature of $E$, where $\nabla^{E}$ is
the holomorphic Hermitian connection  of $(E,h_{E})$. Write
$h_{E}(i\,R^{E}(\cdot),\cdot)=\sum_{a,b,\alpha,\beta}R_{\alpha\bar{\beta}a\bar{b}}\,
(e^{\lor}_{\alpha}\otimes\bar{e}_{\beta}^{\lor})\otimes(\theta_{a}\wedge\bar{\theta}_{b})$,
where $\{e_{\alpha}^{\lor}\}$ (resp. $\{\theta_{a}\}$) is a local {\em unitary} frame of 
$E^{\lor}$ (resp. $\Omega^{1}_{M}$).
Then $(E,h_{E})$ is said to be {\em Nakano semi-positive} if 
$\sum_{a,b,\alpha,\beta}R_{a\bar{b}\alpha\bar{\beta}}\zeta_{a}^{\alpha}\bar{\zeta}_{b}^{\beta}
\geq0$ for all $(\zeta_{a}^{\alpha})\in{\bf C}^{r(n+1)}$. 
Similarly, $(E,h_{E})$ is said to be {\em semi-negative in the dual Nakano sense} if 
$\sum_{a,b,\alpha,\beta}R_{a\bar{b}\alpha\bar{\beta}}\zeta_{b}^{\alpha}\bar{\zeta}_{a}^{\beta}
\leq0$ for all $(\zeta_{b}^{\alpha})\in{\bf C}^{r(n+1)}$. Note the difference of indices
in these two definitions.
By \cite[Lemma 4.3]{Siu82}, $(E,h_{E})$ is Nakano semi-positive if and only if
$(E^{\lor},h_{E^{\lor}})$ is semi-negative in the dual Nakano sense,
where $h_{E^{\lor}}$ is the metric on $E^{\lor}$ induced from $h_{E}$. 
We thank Professor Shigeharu Takayama for pointing out the fact that
the dual of a Nakano semi-positive vector bundle is {\em not necessarily} 
Nakano semi-negative but semi-negative in the dual Nakano sense.

\subsection
{Some results of Takegoshi for Nakano semi-positive vector bundles}
\label{subsect:5.2}
We recall two results for Nakano semi-positive vector bundles from \cite{Takegoshi95}.
Let $\varDelta\subset{\bf C}$ be the unit disc. 
Let $\pi\colon M\to\varDelta$ be a proper surjective holomorphic map 
with critical locus $\varSigma$
and set $M_{t}:=\pi^{-1}(t)$ for $t\in\varDelta$.
Assume that $M$ is a K\"ahler manifold with K\"ahler form $\kappa_{M}$.
Set $\kappa_{M_{t}}:=\kappa_{M}|_{M_{t}}$ and $h_{E_{t}}:=h_{E}|_{M_{t}}$.

\begin{theorem}\label{Theorem5.1}
Assume that $(E,h_{E})$ is Nakano semi-positive.
For every $u\in H^{q}(M,\Omega_{M}^{n+1}(E))$, there exists
$\sigma\in H^{0}(M,\Omega^{n+1-q}_{M}(E))$ with
$$
u=[\sigma\wedge\kappa_{M}^{q}],
\qquad
(\pi^{*}dt)\wedge\sigma=0.
$$
In particular, there exists $v\in H^{0}(M\setminus\varSigma,\Omega^{n-q}_{M/\varDelta}(E))$ with
$$
u|_{M\setminus\varSigma}
=
[v\wedge\kappa_{M}^{q}\wedge(\pi^{*}dt)].
$$
For $t\in\varDelta\setminus\pi(\varSigma)$, 
$v|_{M_{t}}\wedge\kappa_{M_{t}}^{q}
\in A^{n,q}_{M_{t}}(\Omega_{M_{t}}^{n}(E|_{M_{t}}))$ 
is a harmonic form with respect to $\kappa_{M_{t}}$, $h_{E_{t}}$.
\end{theorem}

\begin{pf}
See \cite[Th.\,5.2 (i), (ii)]{Takegoshi95}. Notice that
$v|_{M_{t}}\wedge\kappa_{M_{t}}^{q}=*(v|_{M_{t}})$ 
is a harmonic form with respect to $\kappa_{M_{t}}$, $h_{E_{t}}$,
since $v|_{M_{t}}\in H^{0}(M_{t},\Omega_{M_{t}}^{n-q})$
is holomorphic and since the Hodge star operator preserves harmonic forms.
\end{pf}

Theorem~\ref{Theorem5.1} and its extension by Mourougane-Takayama
\cite[Prop.\,4.4]{MourouganeTakayama09} shall play a key role in Sects.\ref{sect:6} and \ref{sect:9} 
to prove Theorem~\ref{Theorem1.2}.

\begin{lemma}\label{Lemma5.2}
If $(E,h_{E})$ is Nakano semi-positive, 
then $R^{q}\pi_{*}\Omega_{M}^{n+1}(E)$ is locally free for all $q\geq0$. 
\end{lemma}

\begin{pf}
Since $R^{q}\pi_{*}\Omega_{M}^{n+1}(E)$ is torsion free by \cite[Th.\,6.5 (i)]{Takegoshi95} and since
$\dim S=1$, we get the result.
\end{pf}

We refer to \cite[Sect.\,4]{Siu82}, \cite[Chap.\,VII]{Demailly09} 
for more about various notions of positivity and negativity of vector bundles 
and \cite{Berndtsson09}, \cite{MourouganeTakayama08}, \cite{MourouganeTakayama09}, \cite{Takegoshi95} 
for more about the direct images of Nakano semi-positive vector bundles twisted by the relative canonical bundle.

\section
{Asymptotic behavior of equivariant analytic torsion}
\label{sect:6}
\par

\subsection
{Set up}
\label{subsect:6.1}
\par
Let $\kappa_{\mathcal X}$ be the K\"ahler form of $h_{\mathcal X}$. 
In the rest of this paper, we assume that 
{\em $(\xi,h_{\xi})$ is Nakano semi-positive on $X$ and that $(S,0)\cong(\varDelta,0)$.}
By Lemma~\ref{Lemma5.2}, $R^{q}\pi_{*}\omega_{X/S}(\xi)$ is locally free on $S$. 
By shrinking $S$ if necessary, we may also assume that $R^{q}\pi_{*}\omega_{X/S}(\xi)$ is a free
${\mathcal O}_{S}$-module on $S$.
By the $G$-equivariance of $R^{q}\pi_{*}\omega_{X/S}(\xi)$, 
${\rm Hom}_{G}(W,R^{q}\pi_{*}\omega_{X/S}(\xi))\otimes W$ is a  vector bundle on $S$.
By definition, 
$$
\lambda_{G}(\omega_{X/S}(\xi))
=
\prod_{W\in\widehat{G}}\bigotimes_{q\geq0}
\det({\rm Hom}_{G}(W,R^{q}\pi_{*}\omega_{X/}(\xi))\otimes W)^{(-1)^{q}}.
$$
\par
Let $r^{q}_{W}\in{\bf Z}_{\geq0}$ be the rank of 
${\rm Hom}_{G}(W,R^{q}\pi_{*}\Omega_{X}^{n+1}(\xi))\otimes{W}$ 
as a free ${\mathcal O}_{S}$-module on $S$. Let
$\{\Psi_{1},\ldots,\Psi_{r^{q}_{W}}\}\subset
H^{0}(S,{\rm Hom}_{G}(W,R^{q}\pi_{*}\Omega_{X}^{n+1}(\xi))\otimes{W})$
be a free basis of the locally free sheaf 
${\rm Hom}_{G}(W,R^{q}\pi_{*}\Omega_{X}^{n+1}(\xi))\otimes{W}$ on $S$.
Define
$$
\sigma_{W}^{q}
:=
(\Psi_{1}\otimes(\pi^{*}ds)^{-1})\wedge\cdots\wedge(\Psi_{r^{q}_{W}}\otimes(\pi^{*}ds)^{-1})
$$
if ${\rm Hom}_{G}(W,R^{q}\pi_{*}\omega_{X/S}(\xi))\not=0$. We set
$\sigma_{W}^{q}:=1_{\lambda_{W}(\xi)}$ if ${\rm Hom}_{G}(W,R^{q}\pi_{*}\omega_{X/S}(\xi))=0$.
Then $\sigma_{W}^{q}$
generates $\det({\rm Hom}_{G}(W,R^{q}\pi_{*}\omega_{X/S}(\xi))\otimes W)$ on $S$.

\subsection
{Semistable reduction}
\label{subsect:6.2}
\par
Let $T$ be the unit disc in ${\bf C}$. For $0<\epsilon<1$, we set $T(\epsilon):=\{t\in T;\,|t|<\epsilon\}$
and $T^{o}:=T\setminus\{0\}$.
By the semistable reduction theorem \cite[Chap.\,II]{Mumford73}, there is a diagram
$$
\begin{CD}
(Y,Y_{0})@>r>> (X\times_{S}T,X_{0}) @>{\rm pr}_{1}>>  (X,X_{0})
\\
@V f VV @V{\rm pr}_{2}VV  @V \pi VV
\\
(T,0) @>{\rm id}>> (T,0) @>\mu>> (S,0)
\end{CD}
$$
such that $Y_{0}=f^{-1}(0)$ is a reduced normal crossing divisor.
Here $\mu\colon(T,0)\to(S,0)$ is given by $\mu(t)=t^{\nu}$ for some $\nu\in{\bf Z}_{>0}$ and
$r\colon Y\to X\times_{S}T$ is a projective resolution. Set $F={\rm pr}_{1}\circ r\colon Y\to X$.
Since $(F^{*}\xi,h_{F^{*}\xi}:=F^{*}h_{\xi})$ is Nakano semi-positive, we may assume 
by Lemma~\ref{Lemma5.2} that
$R^{q}f_{*}\Omega^{n+1}_{Y}(F^{*}\xi)$ is a free ${\mathcal O}_{T}$-module.
Since $f|_{Y\setminus Y_{0}}\colon Y\setminus Y_{0}\to T^{o}$ is $G$-equivariant, 
$R^{q}f_{*}\Omega^{n+1}_{Y}(F^{*}\xi)|_{T^{o}}$ is a $G$-equivariant holomorphic vector bundle on $T^{o}$.

\begin{lemma}\label{Lemma6.1}
The $G$-action on $R^{q}f_{*}\Omega^{n+1}_{Y}(F^{*}\xi)|_{T^{o}}$
extends to a holomorphic $G$-action on $R^{q}f_{*}\Omega^{n+1}_{Y}(F^{*}\xi)$.
\end{lemma}

\begin{pf}
Let $\rho'\colon Z'\to X\times_{S}T$ be a $G$-equivariant resolution and 
set $\varpi':={\rm pr}_{2}\circ\rho'$ and $\varPi':={\rm pr}_{1}\circ\rho'$. 
By the $G$-equivariance of $\varpi'\colon Z'\to T$, 
$R^{q}\varpi'_{*}\Omega^{n+1}_{Z'}(\varPi^{'*}\xi)$ is a $G$-equivariant holomorphic vector bundle on $T$. 
Since $Y$ is birational to $Z'$, we get an isomorphism
$R^{q}f_{*}\Omega^{n+1}_{Y}(F^{*}\xi)\cong R^{q}\varpi'_{*}\Omega^{n+1}_{Z'}(\varPi^{'*}\xi)$
of holomorphic vector bundles on $T$ by \cite[Th.\,6.9 (i)]{Takegoshi95}, 
which induces the desired $G$-action on $R^{q}f_{*}\Omega^{n+1}_{Y}(F^{*}\xi)$.
\end{pf}

\subsection
{Estimate for the $L^{2}$-metric for the semistable family}
\label{subsect:6.3}
\par
We write $t$ for the coordinate of $T\cong\varDelta$ centered at $0$.
Set $\kappa_{T}=i\,dt\wedge d\bar{t}$, which is a K\"ahler form on $T$. 
By \cite[Sect.\,4.1]{MourouganeTakayama09},
$$
\kappa_{Y}:=F^{*}\kappa_{\mathcal X}+i\,f^{*}\kappa_{T}
$$
is a $C^{\infty}$ $(1,1)$-form on $Y$, which is a K\"ahler form only on $Y\setminus Y_{0}$.
To get an estimate of the $L^{2}$-metric on $R^{q}f_{*}\omega_{Y/T}(F^{*}\xi)$ with respect to
the degenerate K\"ahler form $\kappa_{Y}$, we need an analogue of Theorem~\ref{Theorem5.1}
for the Nakano semi-positive vector bundle $(F^{*}\xi,F^{*}h_{\xi})$ on $(Y,\kappa_{Y})$.
Such an extension of Theorem~\ref{Theorem5.1} was given by Mourougane-Takayama
\cite[Prop.\,4.4]{MourouganeTakayama09}. However, we can not apply it to our situation at one 
for the following reason:
Set
$$
\varSigma:={\rm Sing}(X\times_{S}T),
\qquad
U:=(X\times_{S}T)\setminus\varSigma.
$$
We may assume by \cite[Chap.\,II]{Mumford73} that
$r\colon Y\setminus r^{-1}(\varSigma)\to U$ is an isomorphism and that $r^{-1}(\varSigma)\subset Y_{0}$ 
is a normal crossing divisor. However, it is not clear from the construction in \cite[Chap.\,II]{Mumford73}
if $r\colon Y\to X\times_{S}T$ is a {\em composite of blowing-ups with non-singular centers} 
disjoint from $U$. Since this condition is essential in the construction of a sequence of K\"ahler forms 
on $Y$ approximating $\kappa_{Y}$ (cf. \cite[Proof of Prop.\,4.4 Step 1]{MourouganeTakayama09}), 
we can {\em not} apply the arguments in \cite[Prop.\,4.4]{MourouganeTakayama09} to the bundle
$(F^{*}\xi,h_{F^{*}\xi})$ on $(Y,\kappa_{Y})$ at once. In stead of applying it to $(Y,\kappa_{Y})$, 
we apply it to a manifold $Z$ dominating $Y$, whose construction is as follows.
\par
By \cite[Th.\,13.2]{BierstoneMilman97}, there exists a resolution of the singularity
$r'\colon W\to X\times_{S}T$, which is a composite of blowing-ups with non-singular centers disjoint 
from $U$. We apply \cite[Lemma 1.3.1]{AKMW02} by setting $X_{1}=W$ and $X_{2}=Y$.
As a result, there exist a projective resolution $\rho\colon Z\to X\times_{S}T$ 
and a birational morphism $\varphi\colon Z\to Y$ such that $\rho$ is a composite of blowing-ups 
with non-singular centers disjoint from $U$ and $\rho=r\circ\varphi$. In particular, 
$Z\setminus\rho^{-1}(\varSigma)\cong Y\setminus r^{-1}(\varSigma)\cong U$ and $Z_{0}={\varphi}^{-1}(Y_{0})$
is a normal crossing divisor .
We set $\varpi:={\rm pr}_{2}\circ\rho\colon Z\to T$ and $\varPi:={\rm pr}_{1}\circ\rho\colon Z\to X$.
\par
Regarding $X\times_{S}T$ as a hypersurface of $X\times T$, we get on $Y$
(cf. \cite[Sect.\,4.1]{MourouganeTakayama09})
$$
\kappa_{Y}=r^{*}(\kappa_{X}+\kappa_{T}).
$$
Similarly, $\kappa_{Z}:=\rho^{*}(\kappa_{X}+\kappa_{T})$ is a degenerate K\"ahler form on $Z$
with $\kappa_{Z}=\varphi^{*}\kappa_{Y}$.
Since $Z$ is obtained from $X\times_{S}T$ by a composite of blowing-ups with non-singular centers,
we deduce from e.g. \cite[Prop.\,12.4]{Demailly09},
\cite[Proof of Prop.\,4.4 Step 1]{MourouganeTakayama09} the existence of
a sequence of K\"ahler forms $\{\kappa_{Z,k}\}_{k\geq1}$ on $Z$ such that 
$\kappa_{Z,k}=\kappa_{Z}$ on $Z\setminus\varpi^{-1}(T(\frac{1}{k}))$.
(In fact, we can assume $[\kappa_{Z,k}]|_{Z\setminus Z_{0}}=[\kappa_{Z}]|_{Z\setminus Z_{0}}$
for $k>1$ by an appropriate construction of $\kappa_{Z,k}$. See Sect.\,\ref{subsect:9.2} below.)

\begin{proposition}\label{Proposition:Takegoshi}
For every $\Theta\in H^{q}(Y,\Omega_{Y}^{n+1}(F^{*}\xi))$, there exists a holomorphic differential form 
$\theta\in H^{0}(Y,\Omega_{Y}^{n+1-q}(F^{*}\xi))$ such that
$$
\Theta|_{Y\setminus Y_{0}}=[\theta\wedge\kappa_{Y}^{q}]|_{Y\setminus Y_{0}}
\in H^{q}(Y\setminus Y_{0},\Omega_{Y}^{n+1}(F^{*}\xi)),
\qquad
\theta\wedge f^{*}dt=0.
$$
\end{proposition}

\begin{pf}
Since $\varphi^{-1}(Y)=Z$, we get
$\varphi^{*}\Theta\in H^{q}(Z,\Omega_{Z}^{n+1}(\varPi^{*}\xi))$.
By Theorem~\ref{Theorem5.1} applied to the Nakano semi-positive vector bundle
$(\varPi^{*}\xi,\varPi^{*}h_{\xi})$ on the K\"ahler manifold $(Z,\kappa_{Z,k})$, 
there exists a holomorphic differential form 
$\theta_{k}\in H^{0}(Z,\Omega_{Z}^{n+1-q}(\varPi^{*}\xi))$ such that
\begin{equation}\label{eqn:(9.1)}
\varphi^{*}\Theta=[\theta_{k}\wedge\kappa_{Z,k}^{q}],
\qquad\qquad
(\varpi^{*}ds)\wedge\theta_{k}=0.
\end{equation}
Let ${\mathcal W}\subset T$ be an open subset  such that 
${\mathcal W}\subset T\setminus T(\frac{1}{k})$ for all $k>1$.
Since $\kappa_{Z,k}=\kappa_{Z}$ on $Z\setminus\varpi^{-1}(T(\frac{1}{k}))$,
we get $\kappa_{Z,k}=\kappa_{Z}$ on $\varpi^{-1}({\mathcal W})$ for all $k>1$.
By \cite[Proof of Prop.\,4.4 Step 3]{MourouganeTakayama09}, the equality 
$\kappa_{Z,k}|_{\varpi^{-1}({\mathcal W})}=\kappa_{Z,l}|_{\varpi^{-1}({\mathcal W})}$ implies that
$\theta_{k}|_{\varpi^{-1}({\mathcal W})}=\theta_{l}|_{\varpi^{-1}({\mathcal W})}$ and
hence $\theta_{k}=\theta_{l}$ for all $k,l>1$. We set $\theta_{\infty}:=\theta_{k}$. Then
\begin{equation}\label{eqn:(9.2)}
\varphi^{*}\Theta=[\theta_{\infty}\wedge\kappa_{Z,k}^{q}],
\qquad\qquad
(\varpi^{*}ds)\wedge\theta_{\infty}=0
\end{equation}
for all $k>1$ by \eqref{eqn:(9.1)}. 
Since $\kappa_{Z,k}=\kappa_{Z}$ on $Z\setminus\varpi^{-1}(T(\frac{1}{k}))$,
we get the equality of cohomology classes 
\begin{equation}\label{eqn:(9.3)}
\varphi^{*}\Theta|_{Z\setminus\varpi^{-1}(T(\frac{1}{k}))}
=
[\theta_{\infty}\wedge\kappa_{Z}^{q}]|_{Z\setminus\varpi^{-1}(T(\frac{1}{k}))}
\end{equation}
for all $k>1$. Since $k>1$ is arbitrary, we get by \eqref{eqn:(9.3)}
\begin{equation}\label{eqn:(9.4)}
\varphi^{*}\Theta|_{Z\setminus Z_{0}}=[\theta_{\infty}\wedge\kappa_{Z}^{q}]|_{Z\setminus Z_{0}},
\qquad
Z_{0}:=\varpi^{-1}(0).
\end{equation}
(In fact, $\varphi^{*}\Theta=[\theta_{\infty}\wedge\kappa_{Z}^{q}]$. See Sect.~\ref{subsect:9.3} below.)
\par
Since $\varphi\colon Z\to Y$ induces an isomorphism between $Z\setminus\rho^{-1}(\varSigma)$ and
$Y\setminus r^{-1}(\varSigma)$, we get
$(\varphi^{-1})^{*}\theta_{\infty}\in H^{0}(Y\setminus Y_{0},\Omega_{Y}^{n+1-q}(F^{*}\xi))$.
Let $\omega_{Y}$ be a K\"ahler form on $Y$. 
Since $\theta_{\infty}\wedge\overline{\theta_{\infty}}\wedge\varphi^{*}\omega_{Y}^{q}$
is a $C^{\infty}$ top form on $Z$, we get $(\varphi^{-1})^{*}\theta_{\infty}\in L^{2}_{\rm loc}(Y)$.
In particular, $(\varphi^{-1})^{*}\theta_{\infty}$ extends to a holomorphic differential form on $Y$.
Hence there exists $\theta\in H^{0}(Y,\Omega_{Y}^{n+1-q}(F^{*}\xi))$ such that
$(\varphi^{-1})^{*}\theta_{\infty}=\theta|_{Y\setminus Y_{0}}$. 
Since $\kappa_{Z}=\varphi^{*}\kappa_{Y}$, we get
$$
\Theta|_{Y\setminus Y_{0}}
=
(\varphi^{-1})^{*}\varphi^{*}\Theta|_{Z\setminus Z_{0}}
=
(\varphi^{-1})^{*}[\theta_{\infty}\wedge\kappa_{Z}^{q}]|_{Z\setminus Z_{0}}
=
[\theta\wedge\kappa_{Y}^{q}]|_{Y\setminus Y_{0}}.
$$
Since $\varpi=f\circ\varphi$ and hence $(\varphi^{-1})^{*}\varpi^{*}=f^{*}$, we get
$(f^{*}ds)\wedge\theta=0$ by the relation $(\varpi^{*}ds)\wedge\theta_{\infty}=0$.
This completes the proof.
\end{pf}

In fact, we get $\Theta=[\theta\wedge\kappa_{Y}^{q}]$ in Proposition~\ref{Proposition:Takegoshi}.
See Sect.\,\ref{subsect:9.3}.
\par
Let $\{\Theta_{1},\ldots,\Theta_{r_{W}^{q}}\}$ be a basis of
${\rm Hom}_{G}(W,R^{q}f_{*}\Omega^{n+1}_{Y}(F^{*}\xi))\otimes W$ as a free ${\mathcal O}_{T}$-module. 
Shrinking $T$ if necessary, we may assume 
$\Theta_{\alpha}\in H^{q}(Y,\Omega^{n+1}_{Y}(F^{*}\xi))$.
By Proposition~\ref{Proposition:Takegoshi}, there exist holomorphic differential forms
$\theta_{\alpha}\in H^{0}(Y,\Omega_{Y}^{n-q+1}(F^{*}\xi))$ and
$\Xi_{\alpha}\in H^{0}(Y\setminus Y_{0},\Omega_{Y/T}^{n-q}(F^{*}\xi))$ 
such that
\begin{equation}\label{eqn:(6.1)}
\Theta_{\alpha}|_{Y\setminus Y_{0}}=[\theta_{\alpha}\wedge\kappa_{Y}^{q}]|_{Y\setminus Y_{0}}
\in H^{q}(Y\setminus Y_{0},\Omega_{Y}^{n+1}(F^{*}\xi)),
\end{equation}
\begin{equation}\label{eqn:(6.2)}
\theta_{\alpha}|_{Y\setminus Y_{0}}=\Xi_{\alpha}\wedge f^{*}dt.
\end{equation}
For $t\in T^{o}$, we set
$$
H_{\alpha\bar{\beta}}(t):=
(\Theta_{\alpha}\otimes(f^{*}dt)^{-1}|_{Y_{t}},\Theta_{\beta}\otimes(f^{*}dt)^{-1}|_{Y_{t}})_{L^{2}}.
$$
Then $H(t)=(H_{\alpha\bar{\beta}}(t))$ is a positive-definite $r_{W}^{q}\times r_{W}^{q}$-Hermitian
matrix.
Set $\kappa_{Y_{t}}:=\kappa_{Y}|_{Y_{t}}$ for $t\in T^{o}$. 
Since $\Xi_{\alpha}\wedge\kappa_{Y}^{q}|_{Y_{t}}$ is the harmonic representative of
its class $[\Xi_{\alpha}\wedge\kappa_{Y}^{q}]|_{Y_{t}}$
with respect to the metrics $\kappa_{Y_{t}}$, $h_{F^{*}\xi}|_{Y_{t}}$, we get for all $t\in T^{o}$
$$
H_{\alpha\bar{\beta}}(t)
=
([\Xi_{\alpha}\wedge\kappa_{Y}^{q}]|_{Y_{t}},[\Xi_{\beta}\wedge\kappa_{Y}^{q}]|_{Y_{t}})_{L^{2}}
=
\int_{Y_{t}}
i^{(n-q)^{2}}h_{F^{*}\xi}(\Xi_{\alpha}\wedge\overline{\Xi}_{\beta})|_{Y_{t}}
\wedge\kappa_{Y_{t}}^{q}.
$$
Here $h_{F^{*}\xi}(\Xi_{\alpha}\wedge\overline{\Xi}_{\beta})$ is defined as follows:
Let $\{e_{1},\ldots,e_{r}\}$ be a local frame of $F^{*}\xi$ and
let $\{e_{1}^{\lor},\ldots,e_{r}^{\lor}\}$ be its dual frame of $F^{*}\xi^{\lor}$. 
We can express locally
$\Xi_{\alpha}=\sum_{i}\Xi_{\alpha,i}\otimes e_{i}$ and
$h_{F^{*}\xi}=\sum_{i,j}h_{i\bar{j}}e_{i}^{\lor}\otimes\bar{e}_{j}^{\lor}$,
where $\Xi_{\alpha,i}$ is a local holomorphic section of $\Omega_{Y/T}^{q}$ 
and $h_{i\bar{j}}$ is a local $C^{\infty}$ function. We define
$$
h_{F^{*}\xi}(\Xi_{\alpha}\wedge\overline{\Xi}_{\beta})
:=
\sum_{i,j}h_{i\bar{j}}\Xi_{\alpha,i}\wedge\overline{\Xi}_{\beta,j}
\in
C^{\infty}(Y\setminus Y_{0},\Omega_{Y/T}^{q}\wedge\overline{\Omega}_{Y/T}^{q}).
$$
Then $h_{F^{*}\xi}(\Xi_{\alpha}\wedge\overline{\Xi}_{\beta})|_{Y_{t}}
=h_{F^{*}\xi}(\Xi_{\alpha}|_{Y_{t}}\wedge\overline{{\Xi}_{\beta}|_{Y_{t}}})
\in A^{q,q}_{Y_{t}}$ for all $t\in T^{o}$.

\begin{lemma}\label{Lemma6.2}
One has $H_{\alpha\bar{\beta}}(t)\in\bigoplus_{m=0}^{n}(\log|t|)^{m}\,C^{\infty}(T)$, so that
there exist constants $a_{\alpha\bar{\beta};m}\in{\bf C}$ with
$$
H_{\alpha\bar{\beta}}(t)
=
\sum_{m=0}^{n}a_{\alpha\bar{\beta};m}\,(\log|t|^{2})^{m}+O\left(|t|(\log|t|)^{n}\right)
\qquad
(t\to0).
$$
In particular, $\det H(t)\in\bigoplus_{m=0}^{n}(\log|t|)^{m}\,C^{\infty}(T)$ and
there exist constants $c_{m}\in{\bf R}$, $0\leq m\leq nr_{W}^{q}$ with
$$
\det H(t)=\sum_{m=0}^{nr_{W}^{q}}c_{m}\,(\log|t|^{2})^{m}+O\left(|t|(\log|t|)^{nr_{W}^{q}}\right)
\qquad
(t\to0).
$$
\end{lemma}

\begin{pf}
By \eqref{eqn:(6.2)}, we get on $Y\setminus Y_{0}$
\begin{equation}\label{eqn:(6.3)}
f^{*}(i\,dt\wedge d\bar{t})\wedge
\{i^{(n-q)^{2}}h_{F^{*}{\xi}}(\Xi_{\alpha}\wedge\overline{\Xi}_{\beta})
\wedge\kappa_{Y}^{q}\}
=
i^{(n-q+1)^{2}}h_{F^{*}{\xi}}(\theta_{\alpha}\wedge\overline{\theta}_{\beta})
\wedge\kappa_{Y}^{q}.
\end{equation}
\par
Let $\{{\mathcal V}_{\lambda}\}_{\lambda\in\Lambda}$ be an open covering of $Y$ with $\#\Lambda<+\infty$
such that there is a system of coordinates $(z_{0},\ldots,z_{n})$ on $V_{\lambda}$ with 
$f|_{{\mathcal V}_{\lambda}}(z)=z_{0}\cdots z_{k}$. Here $k$ depends on $\lambda\in\Lambda$.
Let $\{\varrho_{\lambda}\}_{\lambda\in\Lambda}$ be a partition of unity subject to the covering
$\{{\mathcal V}_{\lambda}\}_{\lambda\in\Lambda}$.
On ${\mathcal V}_{\lambda}$, we define
$$
\tau:=\frac{1}{k}\sum_{j=0}^{k}(-1)^{j-1}z_{j}\,
dz_{0}\wedge\cdots\wedge dz_{j-1}\wedge dz_{j+1}
\wedge\cdots\wedge dz_{k}\wedge dz_{k+1}\wedge\cdots\wedge dz_{n}.
$$
Then $\pi^{*}(dt/t)\wedge\tau=dz_{0}\wedge\cdots\wedge dz_{n}$ on ${\mathcal V}_{\lambda}$.
Since $\theta_{\alpha}$ and $\theta_{\beta}$ are holomorphic $n-q+1$-forms on $Y$ 
and since $\kappa_{Y}\in A^{1,1}_{Y}$,
there exists $B_{\alpha\bar{\beta}}(z)\in C^{\infty}_{0}({\mathcal V}_{\lambda})$ such that
\begin{equation}\label{eqn:(6.4)}
\begin{aligned}
\varrho_{\lambda}(z)h_{F^{*}{\xi}}(\theta_{\alpha}\wedge\overline{\theta}_{\beta})\wedge\kappa_{Y}^{n-q}
|_{{\mathcal V}_{\lambda}}
&=
B_{\alpha\bar{\beta}}(z)\,dz_{0}\wedge\cdots\wedge dz_{n}\wedge
\overline{dz_{0}\wedge\cdots\wedge dz_{n}}
\\
&=
(-1)^{n}B_{\alpha\bar{\beta}}(z)\,\left(\frac{\tau}{f^{*}t}\right)
\wedge\overline{\left(\frac{\tau}{f^{*}t}\right)}\wedge
f^{*}(dt\wedge d\bar{t}).
\end{aligned}
\end{equation}
Comparing \eqref{eqn:(6.3)} and \eqref{eqn:(6.4)}, we get
\begin{equation}\label{eqn:(6.5)}
\int_{Y_{t}\cap{\mathcal V}_{\lambda}}
i^{(n-q)^{2}}\,\varrho_{\lambda}(z)h_{F^{*}{\xi}}(\Xi_{\alpha}\wedge\overline{\Xi}_{\beta})|_{Y_{t}}
\wedge\kappa_{Y_{t}}^{q}
=
|t|^{-2}\int_{Y_{t}\cap{\mathcal V}_{\lambda}}
i^{(n-q)^{2}}\,B_{\alpha\bar{\beta}}(z)\,\tau\wedge\overline{\tau}.
\end{equation}
Since
$\tau/f^{*}t=
\frac{1}{k}\sum_{j=0}^{k}(-1)^{j-1}
\frac{dz_{0}}{z_{0}}\wedge\cdots\wedge\frac{dz_{j-1}}{z_{j-1}}
\wedge\frac{dz_{j+1}}{z_{j+1}}\wedge\cdots\wedge
\frac{dz_{k}}{z_{k}}\wedge dz_{k+1}\wedge\cdots\wedge dz_{n}$ 
has only logarithmic singularities, there exists by \eqref{eqn:(6.5)} a constant $C_{0}>0$ such that 
\begin{equation}\label{eqn:(6.6)}
\left|
\int_{Y_{t}\cap{\mathcal V}_{\lambda}}
\varrho_{\lambda}(z)h_{F^{*}{\xi}}(\Xi_{\alpha}\wedge\overline{\Xi}_{\beta})|_{Y_{t}}\wedge\kappa_{Y_{t}}^{q}
\right|
\leq
C_{0}\,(-\log|t|^{2})^{k}
\end{equation}
for all $t\in T^{o}$. On the other hand, we deduce from \cite[p.166 Th.\,4bis.]{Barlet82} that
$f_{*}(B_{\alpha\bar{\beta}}\,\tau\wedge\overline{\tau})\in
C^{\infty}(T)\oplus\bigoplus_{m=0}^{n}|t|^{2}(\log|t|)^{m}\,C^{\infty}(T)$
is of the following form
\begin{equation}\label{eqn:(6.7)}
\begin{aligned}
\int_{Y_{t}\cap{\mathcal V}_{\lambda}}
i^{(n-q)^{2}}\,B_{\alpha\bar{\beta}}(z)\,\tau\wedge\overline{\tau}
&=
a_{\alpha\bar{\beta}}^{(\lambda)}
+
b_{\alpha\bar{\beta}}^{(\lambda)}\,t
+
c_{\alpha\bar{\beta}}^{(\lambda)}\,\bar{t}
+
\sum_{m=0}^{k}a_{\alpha\bar{\beta};m}^{(\lambda)}(t)\,|t|^{2}(\log|t|^{2})^{m}
\end{aligned}
\end{equation}
as $t\to0$, where $a_{\alpha\bar{\beta}}^{(\lambda)}$, $b_{\alpha\bar{\beta}}^{(\lambda)}$, 
$c_{\alpha\bar{\beta}}^{(\lambda)}$ are constants and $a_{\alpha\bar{\beta};m}^{(\lambda)}(t)\in C^{\infty}(T)$. 
Comparing \eqref{eqn:(6.5)}, \eqref{eqn:(6.6)}, \eqref{eqn:(6.7)}, we get
$a_{\alpha\bar{\beta}}^{(\lambda)}=b_{\alpha\bar{\beta}}^{(\lambda)}=c_{\alpha\bar{\beta}}^{(\lambda)}=0$ and
\begin{equation}\label{eqn:(6.8)}
\int_{Y_{t}\cap{\mathcal V}_{\lambda}}
i^{(n-q)^{2}}\,\varrho_{\lambda}(z)h_{F^{*}{\xi}}(\Xi_{\alpha}\wedge\overline{\Xi}_{\beta})|_{Y_{t}}
\wedge\kappa_{Y_{t}}^{q}
=
\sum_{m=0}^{k}a_{\alpha\bar{\beta};m}^{(\lambda)}(t)\,(\log|t|^{2})^{m}.
\end{equation}
Since 
$f_{*}(\varrho_{\lambda}h_{F^{*}{\xi}}(\Xi_{\alpha}\wedge\overline{\Xi}_{\beta})\wedge\kappa_{Y}^{q})\in
\bigoplus_{m=0}^{n}(\log|t|)^{m}\,C^{\infty}(T)$ by \eqref{eqn:(6.8)}, we get
$$
H_{\alpha\bar{\beta}}(t)
=
\sum_{\lambda\in\Lambda}
i^{(n-q)^{2}}f_{*}(\varrho_{\lambda}\,h_{F^{*}{\xi}}(\Xi_{\alpha}\wedge\overline{\Xi}_{\beta})
\wedge\kappa_{Y}^{q})
\in\bigoplus_{m=0}^{n}(\log|t|)^{m}\,C^{\infty}(T).
$$
This completes the proof.
\end{pf}

\begin{lemma}\label{Lemma6.3}
There exists $C>0$ such that the following inequality holds
for all $u=(u_{1},\ldots,u_{r_{W}^{q}})\in{\bf C}^{r_{W}^{q}}$ and $t\in T^{o}$:
$$
\|\sum_{\alpha}u_{\alpha}\Theta_{\alpha}\otimes(f^{*}dt)^{-1}|_{Y_{t}}\|_{L^{2}}^{2}
\geq
C\,|u|^{2}.
$$
In particular, the following inequality holds for all $t\in T^{o}$
$$
\det H(t)\geq C^{r_{W}^{q}}.
$$
\end{lemma}

\begin{pf}
By \cite[Th.\,6.9 (i)]{Takegoshi95},
$\{(\varphi^{*}\Theta_{1})\otimes(\varpi^{*}dt)^{-1},\ldots,
(\varphi^{*}\Theta_{r_{W}^{q}})\otimes(\varpi^{*}dt)^{-1}\}$ 
is a basis of the free ${\mathcal O}_{T}$-module
${\rm Hom}_{G}(W,R^{q}\varpi_{*}\omega_{Z/T}(\varphi^{*}F^{*}\xi))\otimes W$.
Since $\Theta_{\alpha}|_{Y\setminus Y_{0}}=[\Xi_{\alpha}\wedge\kappa_{Y}^{q}\wedge(f^{*}dt)]$
by \eqref{eqn:(6.1)}, \eqref{eqn:(6.2)}, we get
$\varphi^{*}\Theta_{\alpha}|_{Z\setminus Z_{0}}=
[(\varphi^{*}\Xi_{\alpha})\wedge\kappa_{Z}^{q}\wedge(\varpi^{*}dt)]$.
\par
Since $\rho\colon Z\to X\times_{T}S$ is a composite of blowing-ups with non-singular centers
disjoint from $U$ and hence satisfies \cite[Sect.\,2.3 (1), (2)]{MourouganeTakayama09}, 
there is a constant $C_{0}>0$ by  \cite[Lemmas 4.7 and 4.8]{MourouganeTakayama09} 
such that for all $u=(u_{\alpha})\in{\bf C}^{r_{W}^{q}}$ and $t\in T^{o}$,
$$
\|\sum_{\alpha}u_{\alpha}(\varphi^{*}\Xi_{\alpha})\wedge\kappa_{Z}^{q}|_{Z_{t}}\|_{L^{2}}^{2}
=
\|\sum_{\alpha}u_{\alpha}\varphi^{*}\Xi_{\alpha}|_{Z_{t}}\|_{L^{2}}^{2}
\geq
C_{0}\,|u|^{2}.
$$
Since $\kappa_{Z}=\varphi^{*}\kappa_{Y}$ and hence
$$
\|\sum_{\alpha}u_{\alpha}\Theta_{\alpha}\otimes(f^{*}dt)^{-1}|_{Y_{t}}\|_{L^{2}}^{2}
=
\|\sum_{\alpha}u_{\alpha}\Xi_{\alpha}\wedge\kappa_{Y}^{q}|_{Y_{t}}\|_{L^{2}}^{2}
=
\|\sum_{\alpha}u_{\alpha}(\varphi^{*}\Xi_{\alpha})\wedge\kappa_{Z}^{q}|_{Z_{t}}\|_{L^{2}}^{2},
$$
we get the first inequality. 
Since ${}^{t}vH(t)\bar{v}\geq C_{0}\|v\|^{2}$ for all $v\in{\bf C}^{r_{W}^{q}}$ and $t\in T^{o}$
by the first inequality, the smallest eigenvalue of $H(t)$ is greater than or equal to $C_{0}>0$,
which implies $\det H(t)\geq C_{0}^{r_{W}^{q}}$.
\end{pf}

\begin{proposition}\label{Proposition6.5}
There exist an integer $\nu_{W}^{q}\in{\bf Z}_{\geq0}$ and a constant $c_{W}^{q}\in{\bf R}$ such that
as $t\to0$,
$$
\log\det H(t)
=
\nu_{W}^{q}\log\left(-\log|t|^{2}\right)+c_{W}^{q}+O\left(1/\log|t|\right).
$$
\end{proposition}

\begin{pf}
By the second statement of Lemma~\ref{Lemma6.3}, there exists non-zero $c_{m}$ in Lemma~\ref{Lemma6.2}.
Let $c_{N}(\log|t|)^{N}$ be the leading term of the expansion of $\det H(t)$. Namely,
$c_{N}\not=0$ and $c_{m}=0$ for $m>N$ in Lemma~\ref{Lemma6.2}. Then we get  by Lemma~\ref{Lemma6.2}
$$
\det H(t)=c_{N}(\log|t|)^{N}\left(1+O\left(1/\log|t|\right)\right)
\qquad
(t\to0).
$$
The result follows from this estimate.
\end{pf}

\begin{proposition}\label{Proposition6.6}
Let $\varsigma\in\Gamma(T,R^{q}f_{*}\omega_{Y/T}(F^{*}\xi))$. Then 
$R^{q}f_{*}\omega_{Y/T}(F^{*}\xi)/{\mathcal O}_{T}\varsigma$ is a free ${\mathcal O}_{T}$-module 
near $t=0$ if and only if
$$
\log\|\varsigma(t)\|_{L^{2}}=O\left(\log(-\log|t|)\right)
\qquad
(t\to0).
$$
In particular, $\tau\in\Gamma(T^{o},R^{q}f_{*}\omega_{Y/T}(F^{*}\xi))$ extends to a holomorphic
section defined on $T$ if and only if there exists $e\in{\bf Z}_{\geq0}$ with
$$
\log\|\tau(t)\|_{L^{2}}=e\,\log|t|+O\left(\log(-\log|t|)\right)
\qquad
(t\to0).
$$
\end{proposition}

\begin{pf}
Set $l_{q}:=h^{q}(Y_{t},\omega_{Y_{t}}(F^{*}\xi))$.
If $R^{q}f_{*}\omega_{Y/T}(F^{*}\xi)/{\mathcal O}_{T}\varsigma$ is free around $t=0$,
then there is a basis $\{\Theta_{1},\ldots,\Theta_{l_{q}}\}$ of $R^{q}f_{*}\Omega_{Y}^{n+1}(F^{*}\xi)$
around $t=0$ with $\varsigma=\Theta_{1}\otimes(f^{*}dt)^{-1}$. The desired estimate
$C\leq\|\varsigma(t)\|_{L^{2}}\leq C(-\log|t|)^{nl_{q}}$ for all $t\in T^{o}$ follows from
\eqref{eqn:(6.7)} and Lemma~\ref{Lemma6.3} after shrinking $T$ if necessary.
\par
Assume $\log\|\varsigma(t)\|_{L^{2}}=O\left(\log(-\log|t|)\right)$ as $t\to0$.
Let $\{\Theta_{1},\ldots,\Theta_{l_{q}}\}$ be a basis of $R^{q}f_{*}\Omega_{Y}^{n+1}(F^{*}\xi)$
on $T(\epsilon)$ for some $\epsilon>0$. On $T(\epsilon)$, we can express
$\varsigma(t)=\sum_{\alpha}c_{\alpha}(t)\,\Theta_{\alpha}\otimes(f^{*}dt)^{-1}|_{Y_{t}}$,
where $c_{\alpha}(t)\in{\mathcal O}(T(\epsilon))$.
Let $\nu\in{\bf Z}_{\geq0}$ be such that $\varsigma(t)/t^{\nu}\in{\mathcal O}(T(\epsilon))$ and
$\varsigma(t)/t^{\nu+1}\not\in{\mathcal O}(T(\epsilon))$. There is a basis
$\{\Theta'_{1},\ldots,\Theta'_{l_{q}}\}$ of $R^{q}f_{*}\Omega_{Y}^{n+1}(F^{*}\xi)$ with
$\Theta'_{1}\otimes(f^{*}dt)^{-1}=t^{-\nu}\varsigma(t)$. Since
$C\leq\|\Theta'_{1}\otimes(f^{*}dt)^{-1}|_{Y_{t}}\|_{L^{2}}\leq C(-\log|t|)^{nl_{q}}$
by \eqref{eqn:(6.7)} and Lemma~\ref{Lemma6.3}, 
$\log\|\varsigma(t)\|_{L^{2}}=-\nu\log|t|+O(\log(-\log|t|))$ as $t\to0$. 
Since $\log\|\varsigma(t)\|_{L^{2}}=O\left(\log(-\log|t|)\right)$ as $t\to0$ by assumption, we get $\nu=0$.
Hence $\varsigma=\Theta'_{1}\otimes(f^{*}dt)^{-1}$ is a part of a basis of 
$R^{q}f_{*}\Omega_{Y}^{n+1}(F^{*}\xi)$.
\end{pf}

\subsection
{Comparison of the $L^{2}$-metrics}
\label{subsect:6.4}
\par
We recall the following result of Mourougane-Takayama.

\begin{proposition}\label{Proposition:MourouganeTakayama}
There is a natural injection
\begin{equation}\label{eqn:(6.9)}
\varphi\colon
R^{q}f_{*}\omega_{Y/T}(F^{*}\xi)\hookrightarrow\mu^{*}R^{q}\pi_{*}\omega_{X/S}(\xi)
\end{equation}
with the following properties:
\begin{itemize}
\item[(1)]
$\mu^{*}R^{q}\pi_{*}\omega_{X/S}(\xi)/\varphi(R^{q}f_{*}\omega_{Y/T}(F^{*}\xi))$ 
is a torsion sheaf on $T$ supported at $0$.
\item[(2)]
$\varphi$ preserves the $L^{2}$-metrics, i.e.,
$$
\varphi^{*}\mu^{*}h_{R^{q}\pi_{*}\omega_{X/S}(\xi)}=h_{R^{q}f_{*}\omega_{Y/T}(F^{*}\xi)}.
$$
\end{itemize}
Here $h_{R^{q}\pi_{*}\omega_{X/S}(\xi)}$ (resp. $h_{R^{q}f_{*}\omega_{Y/T}(F^{*}\xi)}$)
is the $L^{2}$-metric on $R^{q}\pi_{*}\omega_{X/S}(\xi)$
(resp. $R^{q}f_{*}\omega_{Y/T}(F^{*}\xi)$)
with respect to $\kappa_{\mathcal X}$, $h_{\xi}$ (resp. $\kappa_{Y}$, $F^{*}h_{\xi}$).
\end{proposition}

\begin{pf}
See \cite[Lemmas 3.3 and 4.2]{MourouganeTakayama09}.
\end{pf}

We remark that since $\varphi|_{T^{o}}$ is $G$-equivariant and $\varphi$ is defined on $T$,
$\varphi$ is $G$-equivariant on $T$.
\par
Let ${\frak m}_{0}=t\,{\mathcal O}_{T}\subset{\mathcal O}_{T}$ be the ideal sheaf of $0\in T$. 
For $q\geq0$ and $W\in\widehat{G}$, we set
\begin{equation}\label{eqn:(6.10)}
\delta_{W}^{q}
:=
\dim_{{\mathcal O}_{T}/{\frak m}_{0}}
\frac{{\rm Hom}_{G}(W,\mu^{*}R^{q}\pi_{*}\omega_{X/S}(\xi))\otimes W}
{{\rm Hom}_{G}(W,R^{q}f_{*}\omega_{Y/T}(F^{*}\xi))\otimes W}
\in{\bf Z}_{\geq0}.
\end{equation}

\begin{theorem}\label{Theorem6.8}
By choosing the basis $\{\Psi_{1},\ldots,\Psi_{r_{W}^{q}}\}$ and  $\{\Theta_{1},\ldots,\Theta_{r_{W}^{q}}\}$
appropriately,
there exist integers $e_{1},\ldots,e_{r_{W}^{q}}\geq0$ such that the $r_{W}^{q}\times r_{W}^{q}$-Hermitian 
matrix 
$G(s):=((\Psi_{\alpha}\otimes(\pi^{*}ds)^{-1}|_{X_{s}},\Psi_{\beta}\otimes(\pi^{*}ds)^{-1}|_{X_{s}})_{L^{2}})$ 
is expressed as follows:
$$
G(\mu(t))
=
D(t)\cdot H(t)\cdot\overline{D(t)},
\qquad
D(t)={\rm diag}(t^{-e_{1}},\ldots,t^{-e_{\rho_{q}}}).
$$
In particular, as $s\to 0$,
$$
\log\|\sigma_{W}^{q}(s)\|_{L^{2}}^{2}
=
-\frac{\delta_{W}^{q}}{\deg\mu}\,\log |s|^{2}+\nu_{W}^{q}\log\left(-\log|s|^{2}\right)+c_{W}^{q}+O\left(1/\log|s|\right).
$$
\end{theorem}

\begin{pf}
By choosing the basis $\{\Psi_{\alpha}\}$ and $\{\Theta_{\alpha}\}$ suitably,
there exists by Proposition~\ref{Proposition:MourouganeTakayama} (1) 
an integer  $e_{\alpha}\in{\bf Z}_{\geq0}$ for all $1\leq\alpha\leq r_{W}^{q}$ 
such that the following equality holds on $T$:
$$
\mu^{*}(\Psi_{\alpha}\otimes(\pi^{*}dt)^{-1})|_{Y_{t}}
=t^{-e_{\alpha}}\varphi(\Theta_{\alpha}\otimes(f^{*}dt)^{-1}|_{Y_{t}}).
$$
Then we get
\begin{equation}\label{eqn:(6.11)}
\begin{aligned}
G_{\alpha\bar{\beta}}(\mu(t))
&=
({\Psi}_{\alpha}\otimes(\pi^{*}ds)^{-1}|_{X_{t^{\nu}}},{\Psi}_{\beta}\otimes(\pi^{*}ds)^{-1}|_{X_{t^{\nu}}})_{L^{2}}
\\
&=
h_{R^{q}\pi_{*}\omega_{X/S}(\xi)}
({\Psi}_{\alpha}\otimes(\pi^{*}ds)^{-1},{\Psi}_{\beta}\otimes(\pi^{*}ds)^{-1})(\mu(t))
\\
&=
t^{-e_{\alpha}}\bar{t}^{-e_{\beta}}
h_{R^{q}\pi_{*}\omega_{X/S}(\xi)}
(\varphi({\Theta}_{\alpha}\otimes(f^{*}dt)^{-1}),
\varphi({\Theta}_{\beta}\otimes(f^{*}dt)^{-1}))(\mu(t))
\\
&=
t^{-e_{\alpha}}\bar{t}^{-e_{\beta}}
\varphi^{*}\mu^{*}h_{R^{q}\pi_{*}\omega_{X/S}(\xi)}
({\Theta}_{\alpha}\otimes(f^{*}dt)^{-1},
{\Theta}_{\beta}\otimes(f^{*}dt)^{-1})(t)
\\
&=
t^{-e_{\alpha}}\bar{t}^{-e_{\beta}}
h_{R^{q}f_{*}\omega_{Y/T}(F^{*}\xi)}
({\Theta}_{\alpha}\otimes(f^{*}dt)^{-1},
{\Theta}_{\beta}\otimes(f^{*}dt)^{-1})(t)
\\
&=
t^{-e_{\alpha}}\bar{t}^{-e_{\beta}}H_{\alpha\bar{\beta}}(t),
\end{aligned}
\end{equation}
where the fifth equality follows from Proposition~\ref{Proposition:MourouganeTakayama}\,(2). 
This proves the first equality of Proposition~\ref{Theorem6.8}. Since
\begin{equation}\label{eqn:(6.12)}
\begin{aligned}
\sum_{\alpha}e_{\alpha}
&=
\dim_{{\mathcal O}_{T}/{\frak m}_{0}}
\bigoplus_{\alpha}
{\mathcal O}_{T}\Psi_{\alpha}\otimes(\pi^{*}ds)^{-1}/
{\mathcal O}_{T}\varphi(\Theta_{\alpha}\otimes(f^{*}dt)^{-1})
\\
&=
\dim_{{\mathcal O}_{T}/{\frak m}_{0}}
\frac{{\rm Hom}_{G}(W,\mu^{*}R^{q}\pi_{*}\omega_{X/S}(\xi))\otimes W}
{{\rm Hom}_{G}(W,R^{q}f_{*}\omega_{Y/S}(F^{*}\xi))\otimes W}
=
\delta_{W}^{q},
\end{aligned}
\end{equation}
we get by \eqref{eqn:(6.11)}, \eqref{eqn:(6.12)}
\begin{equation}\label{eqn:(6.13)}
\det G(\mu(t))=|t|^{-2\delta_{W}^{q}}\det H(t).
\end{equation}
Since $|t|=|s|^{\frac{1}{\nu}}$ and $\|\sigma_{W}^{q}(s)\|_{L^{2}}^{2}=\det G(s)$, 
the second equality of Proposition~\ref{Theorem6.8} follows from
\eqref{eqn:(6.13)} and Proposition~\ref{Proposition6.5}.
\end{pf}

We remark that one can get Theorem~\ref{Theorem6.8} by using the method of variation of Hodge structures
when $(\xi,h_{\xi})$ is a trivial Hermtian line bundle on ${\mathcal X}$ \cite[Sect.\,2.2]{Yoshikawa10}.
For an application of Theorem~\ref{Theorem6.8} to the curvature of $L^{2}$-metric, see \cite{Yoshikawa10}.
\par
We define
$$
{\Bbb L}_{g}\left(
\frac{\mu^{*}R\pi_{*}\omega_{X/S}(\xi)}{Rf_{*}\omega_{Y/T}(F^{*}\xi)}
\right)
:=
\sum_{q\geq0}(-1)^{q}{\rm Tr}
[g|_{(\mu^{*}R^{q}\pi_{*}\omega_{X/S}(\xi)/R^{q}f_{*}\omega_{Y/T}(F^{*}\xi))_{0}}].
$$

\begin{lemma}\label{Lemma6.9}
The following identity holds
$$
{\Bbb L}_{g}(\mu^{*}R\pi_{*}\omega_{X/S}(\xi)/Rf_{*}\omega_{Y/T}(F^{*}\xi))
=
\sum_{W\in\widehat{G}}\sum_{q\geq0}
(-1)^{q}\frac{\chi_{W}(g)}{\dim W}\delta_{W}^{q}.
$$
\end{lemma}

\begin{pf}
Since
$$
\frac{\mu^{*}R^{q}\pi_{*}\omega_{X/S}(\xi)}{R^{q}f_{*}\omega_{Y/T}(F^{*}\xi)}
=
\bigoplus_{W\in\widehat{G}}
\frac{{\rm Hom}_{G}(W,\mu^{*}R^{q}\pi_{*}\omega_{X/S}(\xi))\otimes W}
{{\rm Hom}_{G}(W,R^{q}f_{*}\omega_{Y/T}(F^{*}\xi))\otimes W},
$$
we get
\begin{equation}\label{eqn:(6.15)}
\begin{aligned}
\,&
{\rm Tr}\left[g|_{(\mu^{*}R^{q}\pi_{*}\omega_{X/S}(\xi)/R^{q}f_{*}\omega_{Y/T}(F^{*}\xi))_{0}}\right]
\\
&=
\sum_{W\in\widehat{G}}
\left(
\dim_{{\mathcal O}_{T}/{\frak m}_{0}}
\frac{{\rm Hom}_{G}(W,\mu^{*}R^{q}\pi_{*}\omega_{X/S}(\xi))}
{{\rm Hom}_{G}(W,R^{q}f_{*}\omega_{Y/T}(F^{*}\xi))}
\right)
{\rm Tr}[g|_{W}]
\\
&=
\sum_{W\in\widehat{G}}
\left(
\dim_{{\mathcal O}_{T}/{\frak m}_{0}}
\frac{{\rm Hom}_{G}(W,\mu^{*}R^{q}\pi_{*}\omega_{X/S}(\xi))}
{{\rm Hom}_{G}(W,R^{q}f_{*}\omega_{Y/T}(F^{*}\xi))}
\right)
\chi_{W}(g)
\\
&=
\sum_{W\in\widehat{G}}
\left(
\dim_{{\mathcal O}_{T}/{\frak m}_{0}}
\frac{{\rm Hom}_{G}(W,\mu^{*}R^{q}\pi_{*}\omega_{X/S}(\xi))\otimes W}
{{\rm Hom}_{G}(W,R^{q}f_{*}\omega_{Y/T}(F^{*}\xi))\otimes W}
\right)
\frac{\chi_{W}(g)}{\dim W}
=
\sum_{W\in\widehat{G}}\delta_{W}^{q}\frac{\chi_{W}(g)}{\dim W},
\end{aligned}
\end{equation}
from which the result follows.
\end{pf}

\subsection
{Proof of Theorem~\ref{Theorem1.2} }
\label{subsect:6.5}
\par
By the definition of $\sigma_{W}$, $\sigma:=(\sigma_{W})_{W\in\widehat{G}}$ is an admissible 
section of $\lambda_{G}(\omega_{X/S}(\xi))=\lambda_{G}(\omega_{{\mathcal X}/C}(\xi))|_{S}$.
By Theorem~\ref{Theorem4.3}, we get
\begin{equation}\label{eqn:(6.16)}
\log\left\|
\sigma(s)
\right\|_{\lambda_{G}(\omega_{{\mathcal X}/C}(\xi)),Q}^{2}(g)
=
\alpha_{g}(X_{0},\omega_{X/S}(\xi))\,\log|s|^{2}
+O\left(1\right)
\end{equation}
as $s\to0$. On the other hand, we get by Proposition~\ref{Theorem6.8} and Lemma~\ref{Lemma6.9}
\begin{equation}\label{eqn:(6.17)}
\begin{aligned}
\,&
\log\left\|
\sigma(s)
\right\|_{\lambda_{G}(\omega_{{\mathcal X}/C}(\xi)),Q}^{2}(g)
\\
&=
\log\tau_{G}(X_{s},\omega_{X_{s}}(\xi_{s}))(g)
+
\sum_{q\geq0,\,W\in\widehat{G}}(-1)^{q}\frac{\chi_{W}(g)}{\dim W}
\log\|\sigma^{q}_{W}(s)\|_{L^{2}}^{2}
\\
&=
\log\tau_{G}(X_{s},\omega_{X_{s}}(\xi_{s}))(g)+
\\
&\quad
\sum_{q,\,W}(-1)^{q}\frac{\chi_{W}(g)}{\dim W}
\left\{
-\frac{\delta_{W}^{q}}{\deg\mu}\,\log |s|^{2}
+
\nu_{W}^{q}\log(-\log |s|^{2})
+
c_{W}^{q}
+
O\left(\frac{1}{\log|s|}\right)
\right\}
\\
&=
\log\tau_{G}(X_{s},\omega_{X_{s}}(\xi_{s}))(g)
-
\frac{1}{\deg\mu}
{\Bbb L}_{g}(\mu^{*}R\pi_{*}\omega_{X/S}(\xi)/Rf_{*}\omega_{Y/T}(F^{*}\xi))\,\log |s|^{2}
\\
&\quad
+(\sum_{q,\,W}(-1)^{q}\frac{\chi_{W}(g)}{\dim W}\nu_{W}^{q})\,\log(-\log |s|^{2})
+
\sum_{q,\,W}(-1)^{q}\frac{\chi_{W}(g)}{\dim W}c_{W}^{q}
+
O\left(\frac{1}{\log|s|}\right).
\end{aligned}
\end{equation}
Comparing \eqref{eqn:(6.16)} and \eqref{eqn:(6.17)}, we get
\begin{equation}\label{eqn:(6.18)}
\log\tau_{G}(X_{s},\omega_{X_{s}}(\xi_{s}))(g)
=
\beta_{g}\,\log |s|^{2}+\nu_{g}\log(-\log |s|^{2})+c_{g}+O(1/\log|s|)
\end{equation}
as $s\to0$, where we used the following notation in \eqref{eqn:(6.18)}
\begin{equation}\label{eqn:(6.19)}
\begin{array}{ll}
\beta_{g}
&:=
\alpha_{g}(X_{0},\omega_{X/S}(\xi))
+
\frac{1}{\deg\mu}{\Bbb L}_{g}\left(\mu^{*}R\pi_{*}\omega_{X/S}(\xi)/Rf_{*}\omega_{Y/T}(F^{*}\xi)\right),
\\
\nu_{g}
&:=
\sum_{q\geq0,\,W\in\widehat{G}}(-1)^{q+1}\chi_{W}(g)\nu_{W}^{q}/\dim W,
\\
c_{g}
&:=
\sum_{q\geq0,\,W\in\widehat{G}}(-1)^{q+1}\chi_{W}(g)c_{W}^{q}/\dim W.
\end{array}
\end{equation}
This completes the proof.
\qed

\begin{corollary}\label{Corollary6.8}
Let $(E,h_{E})$ be a holomorphic Hermitian vector bundle on ${\mathcal X}$.
If $(E,h_{E})$ is semi-negative in the dual Nakano sense on $X$, then as $s\to0$
$$
(-1)^{n+1}\log\tau_{G}(X_{s},E_{s})(g)=\beta_{g}\,\log|s|^{2}+\nu_{g}\,\log(-\log|s|^{2})+c_{g}+O(1/\log|s|).
$$
Here $\beta_{g},\nu_{g},c_{g}$ are constants defined by the formula \eqref{eqn:(6.19)} by setting $\xi=E^{\lor}$.
\end{corollary}

\begin{pf}
Let $\square_{E_{s}}^{p,q}$ denote the Laplacian acting on $A^{p,q}_{X_{s}}(\xi)$ 
and let $*$ be the Hodge star operator. Since
$*\square_{E_{s}}^{0,q}*^{-1}=\square_{E_{s}^{\lor}}^{n,n-q}
=\square_{\Omega_{X_{s}}^{n}(E_{s}^{\lor})}^{0,n-q}$,
we get the relation
\begin{equation}\label{eqn:(6.20)}
\log\tau_{G}(X_{s},E_{s})(g)
=
(-1)^{n+1}\log\tau_{G}(X_{s},\omega_{X_{s}}(E_{s}^{\lor}))(g).
\end{equation}
Since $(E^{\lor},h_{E^{\lor}})$ is Nakano semi-positive on $X$, 
the result follows from Theorem~\ref{Theorem1.2} and \eqref{eqn:(6.20)}.
\end{pf}

\begin{corollary}\label{Corollary6.9}
If $(\xi,h_{\xi})$ is Nakano semi-positive on $X$ 
and if $X_{0}$ is a reduced normal crossing divisor of $X$, then as $s\to0$
$$
\log\tau_{G}(X_{s},\omega_{X_{s}}(\xi_{s}))(g)
=
\alpha_{g}(X_{0},\omega_{X/S}(\xi))\log|s|^{2}
+
\nu_{g}\log(-\log|s|^{2})+c_{g}+O(\frac{1}{\log|s|}).
$$
\end{corollary}

\begin{pf}
We get ${\Bbb L}_{g}(\mu^{*}R\pi_{*}\omega_{X/S}(\xi)/Rf_{*}\omega_{Y/T}(F^{*}\xi))=0$
in Theorem~\ref{Theorem1.2}, since $\pi\colon X\to S$ is a semistable degeneration.
The result follows from Theorem~\ref{Theorem1.2}.
\end{pf}

\section
{(Log-)Canonical singularities and analytic torsion}
\label{sect:7}

\subsection
{(Log-)Canonical singularities}
\label{sect:7.1}
\par
Let $V$ be an $n$-dimensional normal projective variety with {\em locally free} dualizing sheaf $\omega_{V}$. 
Set $V_{\rm reg}:=V\setminus{\rm Sing}\,V$ and let $i\colon V_{\rm reg}\hookrightarrow V$ be the inclusion.
Then $\omega_{V}=i_{*}\Omega^{n}_{V_{\rm reg}}$. 
The zero divisor of a holomorphic section of $\omega_{V}$ is called a canonical divisor of $V$ 
and is denoted by $K_{V}$. Then $V$ has only {\em canonical} (resp. {\em log-canonical}) singularities
if there exist a resolution $\varphi\colon\widetilde{V}\to V$ 
and an $\varphi$-exceptional normal crossing divisor 
$E=\sum_{i\in I}a_{i}\,E_{i}\subset\widetilde{V}$ such that
$K_{\widetilde{V}}=\varphi^{*}K_{V}+E$ and $a_{i}\geq0$ (resp. $a_{i}\geq-1$)
for all $i\in I$. Here $E_{i}$ are irreducible and reduced divisors of $\widetilde{V}$. 
If $V$ has only canonical singularities, then 
$\varphi^{*}\colon H^{0}(V,\omega_{V})\to H^{0}(\widetilde{V},\omega_{\widetilde{V}})$ 
is an isomorphism. 
In particular, every element of $H^{0}(V,\omega_{V})$ is square integrable with respect to 
any Hermitian metric on $V_{\rm reg}$.
\par
Let $W$ be a smooth projective manifold and let $D\subset W$ be a divisor.
A birational morphism $\varphi\colon\widetilde{W}\to W$ between smooth projective manifolds 
is called an embedded resolution of $D$ if $\varphi$ is an isomorphism 
between $\widetilde{W}\setminus \varphi^{-1}({\rm Sing}\,D)$ 
and $W\setminus{\rm Sing}\,D$ such that $\widetilde{D}$, 
the proper transform of $D$, is smooth. 
The pair $(W,D)$ has only canonical (resp. log-canonical) singularities if there exist an embedded resolution 
$\varphi\colon\widetilde{W}\to W$ of $D$ and an $\varphi$-exceptional normal crossing divisor
$E=\sum_{i\in I}a_{i}\,E_{i}\subset\widetilde{W}$ such that $\varphi^{-1}(D)$ is a normal crossing divisor 
and such that $K_{\widetilde{W}}+\widetilde{D}=\varphi^{*}(K_{W}+D)+E$ with
$a_{i}\geq0$ (resp. $a_{i}\geq-1$) for all $i\in I$.
If $D$ is a reduced normal crossing divisor of $W$, then the pair $(W,D)$ has only log-canonical singularities.
\par
If $X_{0}$ is reduced and normal and has only canonical (resp. log-canonical) singularities, 
then the pair $(X,X_{0})$ has only canonical (resp. log-canonical) singularities by \cite{Stevens88}, 
\cite[Th.\,7.9]{Kollar97} (resp. \cite[Th.\,7.5]{Kollar97}).
By \cite[Cor.\,11.13]{Kollar97}, the condition that $X_{0}$ has only canonical singularities
is equivalent to the one that $X_{0}$ has only {\em rational} singularities since $X_{0}$ is reduced and normal.
We refer to e.g. \cite{Kollar97} and the references therein for more about related notions of singularities.

\subsection
{Integration along fibers and (log-)canonical singularities}
\label{sect:7.2}
\par
Let $\varrho\in A^{n+1,n+1}_{X}$ and $\chi\in A^{n+1,n}_{X}$.
Write $\varrho=\pi^{*}(ds\wedge d\bar{s})\wedge R$ and $\chi=(\pi^{*}ds)\wedge K$, where
$R\in C^{\infty}(X\setminus{\rm Sing}\,X_{0},\Omega_{X/S}^{n}\wedge\overline{\Omega_{X/S}^{n}})$
and
$K\in C^{\infty}(X\setminus{\rm Sing}\,X_{0},\Omega_{X/S}^{n}\wedge\overline{\Omega_{X/S}^{n}})$.
Define ${\mathcal R}(s)\in C^{\infty}(S^{o})$ and ${\mathcal K}(s)\in C^{\infty}(S^{o})$ as
$$
{\mathcal R}(s):=\int_{X_{s}}R|_{X_{s}},
\qquad
{\mathcal K}(s):=\int_{X_{s}}K|_{X_{s}},
$$ 
so that $\pi_{*}(\varrho)={\mathcal R}(s)\,ds\wedge d\bar{s}$ and $\pi_{*}(\chi)={\mathcal K}(s)\,ds$.

\begin{lemma}\label{Lemma7.1}
Assume that $X_{0}$ is a reduced divisor of ${\mathcal X}$.
\begin{itemize}
\item[(1)]
If the pair $({\mathcal X},X_{0})$ has only canonical singularities,
then ${\mathcal R}(s),{\mathcal K}(s)\in{\mathcal B}(S)$ and
${\mathcal R}(0)=\int_{(X_{0})_{\rm reg}}R|_{(X_{0})_{\rm reg}}$,
${\mathcal K}(0)=\int_{(X_{0})_{\rm reg}}K|_{(X_{0})_{\rm reg}}$.
\item[(2)]
If the pair $({\mathcal X},X_{0})$ has only log-canonical singularities, then
there exists $C>0$ such that $|{\mathcal R}(s)|\leq C\,(-\log|s|)^{n}$.
\end{itemize}
\end{lemma}

\begin{pf}
{\bf Step 1 }
Since the pair $({\mathcal X},X_{0})$ has only canonical (resp. log-canonical) singularities,
there exist an embedded resolution $\varphi\colon{\mathcal Z}\to{\mathcal X}$ of $X_{0}$
and an $\varphi$-exceptional normal crossing divisor $E\subset{\mathcal Z}$ such that 
$K_{\mathcal Z}+\widetilde{X}_{0}=\varphi^{*}(K_{\mathcal X}+X_{0})+E$,
$E=\sum_{i\in I}a_{i}D_{i}$. Here $\widetilde{X}_{0}\subset {\mathcal Z}$ is the proper 
transform of $X_{0}$, $D_{i}$ are irreducible and reduced divisors of ${\mathcal Z}$
and $a_{i}\geq0$ (resp. $a_{i}\geq-1$) for all $i\in I$.
We may assume that $Z_{0}:=\varphi^{-1}(X_{0})$ is a normal crossing divisor of ${\mathcal Z}$
with $\widetilde{X}_{0}\cup E\subset Z_{0}$.
Set $\pi':=\pi\circ\varphi$. Then $Z_{0}=(\pi')^{-1}(0)$ and $Z_{s}:=(\pi')^{-1}(s)\cong X_{s}$ for $s\in S^{o}$.
\par
Let ${\mathcal V}\subset X$ be an open subset, 
on which there is a nowhere vanishing holomorphic $n+1$-form $\Theta$. 
By an argument using partition of unity, it suffices to prove
the assertion when $\varrho$ and $\chi$ are supported in ${\mathcal V}$.
In what follows, we assume ${\rm supp}\,\varrho\subset{\mathcal V}$ and 
${\rm supp}\,\chi\subset{\mathcal V}$.
\par
Let $p\in Z_{0}\cap\varphi^{-1}({\mathcal V})$. Since $X_{0}$ is a reduced divisor of $X$,
there is a coordinate neighborhood $(U,(z_{0},\cdots,z_{n}))$ centered at $p$ satisfying
$\varphi(U)\subset{\mathcal V}$ and one of the following (a), (b): 
Write $\pi'|_{U}(z)=z_{0}^{e_{0}}\cdots z_{n}^{e_{n}}$, $e_{i}\in{\bf Z}_{\geq0}$.
\begin{itemize}
\item[(a)]
If $p\in\widetilde{X}_{0}$, then
$Z_{s}\cap U=\{z\in U;\,z_{0}z_{1}^{e_{1}}\cdots z_{n}^{e_{n}}=s\}$ and $\widetilde{X}_{0}\cap U=\{z_{0}=0\}$.
\item[(b)]
If $p\not\in\widetilde{X}_{0}$, then
$Z_{s}\cap U=\{z\in U;\,
z_{0}^{e_{0}}z_{1}^{e_{1}}\cdots z_{n}^{e_{n}}=s\}$, $e_{0}>0$ and $\widetilde{X}_{0}\cap U=\emptyset$.
\end{itemize}
Since $e_{0}=1$ in case (a), we get the following expression on $U$: 
\begin{equation}\label{eqn:(7.1)}
(\pi')^{*}\left(\frac{ds}{s}\right)
=
\varphi^{*}\pi^{*}\left(\frac{ds}{s}\right)
=
\begin{cases}
\begin{array}{lll}
\frac{dz_{0}}{z_{0}}+\sum_{i\geq1}e_{i}\frac{dz_{i}}{z_{i}}
&{\rm if }&p\in\widetilde{X}_{0},
\\
\sum_{i\geq0}e_{i}\frac{dz_{i}}{z_{i}}
&{\rm if }&p\not\in\widetilde{X}_{0}.
\end{array}
\end{cases}
\end{equation}
%
\newline{\bf Step 2 }
Since $\Theta/\pi^{*}s\in H^{0}({\mathcal V},\omega_{\mathcal X}(X_{0}))$, 
$\varphi^{*}(\Theta/\pi^{*}s)$ is a meromorphic canonical form on 
${\mathcal Z}\cap\varphi^{-1}({\mathcal V})$ 
with at most logarithmic pole along $\widetilde{X}_{0}\cap\varphi^{-1}({\mathcal V})$ 
(resp. $\widetilde{Z}_{0}\cap\varphi^{-1}({\mathcal V})$)
by the assumption 
$K_{\mathcal Z}+\widetilde{X}_{0}=\varphi^{*}(K_{\mathcal X}+X_{0})+\sum_{i\in I}a_{i}D_{i}$,
$a_{i}\geq0$ (resp. $a_{i}\geq-1$).
Hence we get the following expression on $U$
\begin{equation}\label{eqn:(7.2)}
\varphi^{*}\left(\frac{\Theta}{\pi^{*}s}\right)
=
\begin{cases}
\begin{array}{llll}
a(z)\,
\frac{dz_{0}}{z_{0}}\wedge dz_{1}\wedge\cdots\wedge dz_{n}
&{\rm if}&p\in\widetilde{X}_{0}\cap U,
&({\mathcal X},X_{0})\hbox{ is C},
\\
b(z)\,dz_{0}\wedge\cdots\wedge dz_{n}
&{\rm if}&p\not\in\widetilde{X}_{0}\cap U,
&({\mathcal X},X_{0})\hbox{ is C},
\\
a'(z)\,\frac{dz_{0}}{z_{0}}\wedge\cdots\wedge\frac{dz_{n}}{z_{n}}
&{\rm if}&p\in Z_{0}\cap U,
&({\mathcal X},X_{0})\hbox{ is LC},
\end{array}
\end{cases}
\end{equation}
where $a(z),b(z),a'(z)\in{\mathcal O}(U)$.
Here we wrote C (resp. LC) for canonical (resp. log-canonical).
Let $\Xi\in\Gamma({\mathcal V}\setminus X_{0},\Omega^{n}_{X/S})$ be such that 
$\Theta=\Xi\wedge\pi^{*}ds$. Since
\begin{equation}\label{eqn:(7.3)}
\varphi^{*}(\Theta/\pi^{*}s)
=
\varphi^{*}\Xi\wedge(\pi')^{*}(ds/s),
\end{equation}
$\varphi^{*}\Xi$ is expressed as follows on $U\setminus Z_{0}$ 
by \eqref{eqn:(7.1)}, \eqref{eqn:(7.2)}, \eqref{eqn:(7.3)}:
\begin{equation}\label{eqn:(7.4)}
\varphi^{*}\Xi
=
\begin{cases}
\begin{array}{llll}
a(z)\,dz_{1}\wedge\cdots\wedge dz_{n}\mod(\pi')^{*}ds
&{\rm if}&p\in\widetilde{X}_{0}\cap U,
&({\mathcal X},X_{0})\hbox{ is C},
\\
\frac{z_{0}}{e_{0}}b(z)\,dz_{1}\wedge\cdots\wedge dz_{n}
\mod(\pi')^{*}ds
&{\rm if}&p\not\in\widetilde{X}_{0}\cap U,
&({\mathcal X},X_{0})\hbox{ is C},
\\
a'(z)\,\frac{dz_{1}}{z_{1}}\wedge\cdots\wedge\frac{dz_{n}}{z_{n}}
\mod(\pi')^{*}ds
&{\rm if}&p\in Z_{0}\cap U,
&({\mathcal X},X_{0})\hbox{ is LC}.
\end{array}
\end{cases}
\end{equation}
If $({\mathcal X},X_{0})$ has only canonical singularities, then
we set $c(z):=a(z)$ in case (a) and $c(z):=z_{0}b(z)/e_{0}$ in case (b). 
If $({\mathcal X},X_{0})$ has only log-canonical singularities, then we set $c(z):=a'(z)$.
We always have $c(z)\in{\mathcal O}(U)$. 
\par
By \eqref{eqn:(7.4)}, we get for all $F\in C_{0}^{\infty}(U)$ and 
$s\in S^{o}$
$$
\int_{Z_{s}\cap U}
F\,\varphi^{*}(\Xi\wedge\overline{\Xi})|_{Z_{s}\cap U}
=
\int_{Z_{s}\cap U}F(z)\,|c(z)|^{2}\,
dz_{1}\wedge\cdots\wedge dz_{n}\wedge\overline{dz_{1}\wedge\cdots\wedge dz_{n}}.
$$
\newline{\bf Step 3 }
{\em In Steps 3,4, we assume that $({\mathcal X},X_{0})$ has only canonical singularities.} 
Since the integrand is a $C^{\infty}$ $(n,n)$-form on $U$, we get by \cite[Th.\,1]{Barlet82}
\begin{equation}\label{eqn:(7.5)}
\begin{aligned}
\lim_{s\to0}\int_{Z_{s}\cap U}F\,\varphi^{*}(\Xi\wedge\overline{\Xi})
&=
\int_{Z_{0}\cap U}F(z)\,|c(z)|^{2}\,
dz_{1}\wedge\cdots\wedge dz_{n}
\wedge\overline{dz_{1}\wedge\cdots\wedge dz_{n}}
\\
&=
\int_{\widetilde{X}_{0}\cap U}F(z)\,|c(z)|^{2}\,
dz_{1}\wedge\cdots\wedge dz_{n}
\wedge\overline{dz_{1}\wedge\cdots\wedge dz_{n}}
\\
&=
\int_{\widetilde{X}_{0}\cap U}F\,\varphi^{*}(\Xi\wedge\overline{\Xi}).
\end{aligned}
\end{equation}
Here we get the second equality as follows. 
In case (a), $(Z_{0}\setminus\widetilde{X}_{0})\cap U$ is defined locally by
the equation $z_{1}^{e_{1}}\cdots z_{n}^{e_{n}}=0$.
Hence one of $dz_{1},\ldots,dz_{n}$ vanishes on $(Z_{0}\setminus\widetilde{X}_{0})\cap U$,
which implies the second equality of \eqref{eqn:(7.5)}. 
In case (b), let $x\in(Z_{0}\setminus\widetilde{X}_{0})\cap U$. 
Then one of $z_{0},\ldots,z_{n}$ vanishes on a neighborhood $W$ of $x$
in $(Z_{0}\setminus\widetilde{X}_{0})\cap U$. 
If $z_{0}|_{W}=0$, then $c_{i}|_{W}=(z_{0}b_{i}/e_{0})|_{W}=0$.
If $z_{j}|_{W}=0$ for some $j>0$, then $dz_{j}|_{W}=0$. 
Since $(Z_{0}\setminus\widetilde{X}_{0})\cap U$ is covered by these $W$, 
we get
$$
\int_{(Z_{0}\setminus\widetilde{X}_{0})\cap U}
F(z)\,|c(z)|^{2}\,dz_{1}\wedge\cdots\wedge dz_{n}\wedge\overline{dz_{1}\wedge\cdots\wedge dz_{n}}=0.
$$
\newline{\bf Step 4 }
Let $\{U_{\alpha}\}_{\alpha\in A}$ be a covering of $\pi^{-1}({\mathcal V})$ with $\#A<\infty$ 
such that on $U_{\alpha}$, there is a system of coordinates satisfying (a) or (b) as in Step 2. 
Let $\{\chi_{\alpha}\}_{\alpha\in A}$ be a partition of unity subject to the covering $\{U_{\alpha}\}_{\alpha\in A}$. 
\par
Since $\varrho\in A^{n+1,n+1}_{\mathcal X}$ and $\chi\in A^{n+1,n}_{\mathcal X}$ have compact support 
in ${\mathcal V}$, there exist $m\in C^{\infty}({\mathcal V})$ and $B\in A^{0,q}_{\mathcal X}$ with 
${\rm supp}\,B\subset{\mathcal V}$ such that 
$$
\varrho
=
(-1)^{n}m\,\Theta\wedge\overline{\Theta}
=
\pi^{*}(ds\wedge d\bar{s})\wedge m\,\Xi\wedge\overline{\Xi},
\qquad
\chi=(-1)^{n}\Theta\wedge B=(\pi^{*}ds)\wedge\Xi\wedge B.
$$
Since $\pi_{*}\varrho=\pi_{*}(m\,\Xi\wedge\overline{\Xi})\,ds\wedge d\bar{s}$ 
and $\pi_{*}\chi=\pi_{*}(\Xi\wedge B)\,ds$,
we get ${\mathcal R}(s)=\pi_{*}(m\,\Xi\wedge\overline{\Xi})$ and ${\mathcal K}(s)=\pi_{*}(\Xi\wedge B)$. 
Since
\begin{equation}\label{eqn:(7.7)}
\begin{aligned}
\lim_{s\to0}{\mathcal R}(s)
&=
\lim_{s\to0}
\int_{X_{s}}(m\,\Xi\wedge\overline{\Xi})|_{X_{s}}
=
\lim_{s\to0}\sum_{\alpha\in A}
\int_{Z_{s}\cap U_{\alpha}}
\chi_{\alpha}\,\varphi^{*}(m\,\Xi\wedge\overline{\Xi})
|_{Z_{s}\cap U_{\alpha}}
\\
&=
\sum_{\alpha\in A}
\int_{\widetilde{X}_{0}\cap U_{\alpha}}
\chi_{\alpha}\,\varphi^{*}(m\,\Xi\wedge\overline{\Xi})
|_{\widetilde{X}_{0}\cap U_{\alpha}}
=
\int_{\widetilde{X}_{0}}
\varphi^{*}(m\,\Xi\wedge\overline{\Xi})
|_{\widetilde{X}_{0}}
\\
&=
\int_{(X_{0})_{\rm reg}}
m\,\Xi\wedge\overline{\Xi}|_{(X_{0})_{\rm reg}}
={\mathcal R}(0)
\end{aligned}
\end{equation}
by \eqref{eqn:(7.5)}, we get ${\mathcal R}(s)\in C^{0}(S)$ by \eqref{eqn:(7.7)}. 
Since $|s|^{2}{\mathcal R}(s)\in{\mathcal B}(S)$ by \cite[Lemma 9.2]{Yoshikawa07}, we get
${\mathcal R}(s)\in C^{0}(S)\cap |s|^{-2}{\mathcal B}(S)={\mathcal B}(S)$. 
\par
Let $B\in A^{0,n}_{\mathcal X}$ with ${\rm supp}\,B\subset{\mathcal B}$. 
Since $B\wedge\pi^{*}d\bar{s}=h\,\overline{\Theta}$ with some $h\in C_{0}^{\infty}({\mathcal V})$, we get
$B=h\,\overline{\Xi}\mod\pi^{*}d\bar{s}$.
Since $\Xi\wedge B=h\Xi\wedge\overline{\Xi}\mod\pi^{*}ds$ with $h\in C_{0}^{\infty}({\mathcal V})$
and hence ${\mathcal K}(s)=\int_{X_{s}}h\Xi\wedge\overline{\Xi}|_{X_{s}}$, we get by \eqref{eqn:(7.7)}
\begin{equation}\label{eqn:(7.8)}
\lim_{s\to0}{\mathcal K}(s)
=
\int_{(X_{0})_{\rm reg}}\Xi\wedge B|_{(X_{0})_{\rm reg}}
={\mathcal K}(0)
\end{equation}
and ${\mathcal K}(s)\in{\mathcal B}(S)$.
This proves (1).
\newline{\bf Step 5 }
Assume that $({\mathcal X},X_{0})$ has only log-canonical singularities. 
By \eqref{eqn:(7.4)}, there exists a constant $C>0$ such that
$$
\begin{aligned}
\left|\int_{Z_{s}\cap U}F\,\varphi^{*}(\Xi\wedge\overline{\Xi})\right|
&\leq
C\,\left|
\int_{z\in U,\,z_{0}^{e_{1}}\cdots z_{n}^{e_{n}}=s}
\frac{dz_{1}}{z_{1}}\wedge\cdots\wedge\frac{dz_{n}}{z_{n}}
\wedge
\overline{\frac{dz_{1}}{z_{1}}\wedge\cdots\wedge\frac{dz_{n}}{z_{n}}}
\right|
\\
&\leq
C\,(-\log|s|)^{n}.
\end{aligned}
$$
This, together with the equality
${\mathcal R}(s)=\sum_{\alpha\in A}\int_{Z_{s}\cap U_{\alpha}}
\chi_{\alpha}\,\varphi^{*}(m\,\Xi\wedge\overline{\Xi})|_{Z_{s}\cap U_{\alpha}}$, 
implies the desired estimate $|{\mathcal R}(s)|\leq C\,(-\log|s|)^{n}$. 
This proves (2).
\end{pf}

\subsection
{The case of canonical singularities}
\label{sect:7.3}
\par
We assume that {\em $X_{0}$ is reduced and the pair $({\mathcal X},X_{0})$ has only canonical singularities}.
Set $l_{q}:=h^{q}(X_{s},\Omega_{X_{s}}^{n}(\xi_{s}^{\lor}))$.

\begin{theorem}\label{Theorem7.2}
The $L^{2}$ metric on $R^{q}\pi_{*}\omega_{X/S}(\xi)|_{S^{o}}$ extends to a continuous Hermitian metric 
of class ${\mathcal B}(S)$ on $R^{q}\pi_{*}\omega_{X/S}(\xi)|_{S}$ for all $q\geq0$.
\end{theorem}

\begin{pf}
After shrinking $S$ if necessary, we get by Lemma~\ref{Lemma5.2} sections
$\Psi_{1},\ldots,\Psi_{l_{q}}\in H^{q}(X,\Omega^{n+1}_{X}(\xi))$ such that
$R^{q}\pi_{*}\Omega^{n+1}_{X}(\xi)|_{S}={\mathcal O}_{S}\Psi_{1}\oplus\cdots\oplus{\mathcal O}_{S}\Psi_{l_{q}}$. 
We set
\begin{equation}\label{eqn:(7.9)}
g_{\alpha\bar{\beta}}(s)
:=
(\Psi_{\alpha}\otimes(\pi^{*}ds)^{-1}|_{X_{s}},
\Psi_{\beta}\otimes(\pi^{*}ds)^{-1}|_{X_{s}})_{L^{2}}
\end{equation}
for $s\in S^{o}$. Then $g_{\alpha\bar{\beta}}\in C^{\infty}(S^{o})$ and $\det(g_{\alpha\bar{\beta}})>0$ on $S^{o}$.
It suffices to prove that $g_{\alpha\bar{\beta}}\in{\mathcal B}(S)$ and $\det(g_{\alpha\bar{\beta}})>0$ on $S$. 
\par
There exist $\psi_{\alpha}\in H^{0}(X,\Omega_{X}^{n-q+1}(\xi))$ and 
$\Xi_{\alpha}\in H^{0}(X\setminus{\rm Sing}\,X_{0},\Omega_{X/S}^{n-q}(\xi))$ by Theorem~\ref{Theorem5.1} 
such that 
$$
\Psi_{\alpha}=[\psi_{\alpha}\wedge\kappa_{\mathcal X}^{q}],
\qquad
\psi_{\alpha}|_{X\setminus{\rm Sing}\,X_{0}}=\Xi_{\alpha}\wedge\pi^{*}ds.
$$ 
In $H^{q}(X\setminus X_{0},\omega_{X/S}(\xi))$, we get the equality
\begin{equation}\label{eqn:(7.10)}
\Psi_{\alpha}\otimes(\pi^{*}ds)^{-1}=[\Xi_{\alpha}\wedge\kappa_{\mathcal X}^{q}],
\qquad
\alpha=1,\ldots,l_{q}
\end{equation}
under the canonical identification
$\omega_{X/S}(\xi)|_{X\setminus X_{0}}=\Omega_{X/S}^{n}(\xi)|_{X\setminus X_{0}}$.
Since $\Xi_{\alpha}$ is holomorphic and hence $\Xi_{\alpha}\wedge\kappa_{\mathcal X}^{q}|_{X_{s}}$ 
is the harmonic representative of the cohomology class
$\Psi_{\alpha}\otimes(\pi^{*}ds)^{-1}|_{X_{s}}$ for $s\in S^{o}$, 
we deduce from \eqref{eqn:(7.9)}, \eqref{eqn:(7.10)} that
\begin{equation}\label{eqn:(7.11)}
\pi_{*}(i^{(n-q+1)^{2}}h_{\xi}(\psi_{\alpha}\wedge\overline{\psi}_{\beta})\wedge\kappa_{\mathcal X}^{q})
=
i\,g_{\alpha\bar{\beta}}(s)\,ds\wedge d\bar{s}.
\end{equation}
Since $\psi_{\alpha}\in H^{0}(X,\Omega_{X}^{n-q+1}(\xi))$, we see that
$h_{\xi}(\psi_{\alpha}\wedge\overline{\psi_{\beta}})\wedge\kappa_{\mathcal X}^{q}$ 
is a $C^{\infty}$ $(n+1,n+1)$-form on $X$. 
By Lemma~\ref{Lemma7.1} (1) and \eqref{eqn:(7.11)}, we get
$g_{\alpha\bar{\beta}}(s)\in{\mathcal B}(S)$.
\par
By \eqref{eqn:(6.13)} and the second inequality of Lemma~\ref{Lemma6.3}, 
there exists $C>0$ such that
\begin{equation}\label{eqn:(7.12)}
\begin{aligned}
\,&
\det\left(
{\Psi}_{\alpha}\otimes(\pi^{*}ds)^{-1}|_{X_{\mu(t)}},
{\Psi}_{\beta}\otimes(\pi^{*}ds)^{-1}|_{X_{\mu(t)}}
\right)_{L^{2}}
\\
&=
|t|^{-2\delta_{W}^{q}}
\det\left(
{\Theta}_{\alpha}\otimes(f^{*}dt)^{-1}|_{Y_{t}},
{\Theta}_{\beta}\otimes(f^{*}dt)^{-1}|_{Y_{t}}
\right)_{L^{2}}
\geq
C\,|t|^{-2\delta_{W}^{q}}.
\end{aligned}
\end{equation}
Since $\delta_{W}^{q}\geq0$, there exists $\epsilon>0$ such that for all $s\in S^{o}$,
\begin{equation}\label{eqn:(7.13)}
\det(g_{\alpha\bar{\beta}}(s))
=
\det\left(
{\Psi}_{\alpha}\otimes(\pi^{*}ds)^{-1}|_{X_{s}},{\Psi}_{\beta}\otimes(\pi^{*}ds)^{-1}|_{X_{s}}
\right)_{L^{2}}
\geq
\epsilon>0.
\end{equation}
Since $g_{\alpha\bar{\beta}}(s)\in{\mathcal B}(S)$, we get
$\det(g_{\alpha\bar{\beta}}(0))\geq\epsilon>0$ by \eqref{eqn:(7.13)}.
\end{pf}

\begin{corollary}\label{Corollary7.3}
There exist $c\in{\bf R}$, $r\in{\bf Q}_{>0}$ and $\nu\in{\bf Z}_{\geq0}$ such that as $s\to0$
$$
\log\|\sigma_{W}^{q}(s)\|_{L^{2}}^{2}=c+O\left(|s|^{r}(\log|s|)^{\nu}\right).
$$
\end{corollary}

\begin{pf}
Since 
${\rm Hom}_{G}(W,R^{q}\pi_{*}\omega_{X/S}(\xi))\otimes W\subset R^{q}\pi_{*}\omega_{X/S}(\xi)$
is a holomorphic subbundle, the result follows from Theorem~\ref{Theorem7.2}. 
\end{pf}

{\bf Proof of Theorem~\ref{Theorem1.4} (2) }
By the definition of equivariant Quillen metrics,
\begin{equation}\label{eqn:(7.15)}
\begin{aligned}
\,&
\log\tau_{G}(X_{s},\omega_{X_{s}}(\xi_{s}))(g)
\\
&=
\log\|\sigma(s)\|_{\lambda_{G}(\omega_{X/S}(\xi)),Q}^{2}(g)
-
\sum_{q\geq0,\,W\in\widehat{G}}(-1)^{q}\frac{\chi_{W}(g)}{\dim W}
\log\|\sigma^{q}_{W}(s)\|_{L^{2}}^{2}.
\end{aligned}
\end{equation}
The result follows from \eqref{eqn:(7.15)}, Theorem~\ref{Theorem4.3} and 
Corollary~\ref{Corollary7.3}.
\qed

\begin{corollary}\label{Corollary7.4}
If $X_{0}$ is reduced, normal and has only canonical singularities and 
if $(\xi,h_{\xi})$ is semi-negative in the dual Nakano sense on $X$, 
then there exist $c\in{\bf C}$, $r\in{\bf Q}_{>0}$, $l\in{\bf Z}_{\geq0}$ such that as $s\to0$
$$
\log\tau_{G}(X_{s},\xi_{s})(g)
=
(-1)^{n+1}{\frak a}_{g}(X_{0},\xi)\log|s|^{2}+c+O\left(|s|^{r}(\log|s|)^{l}\right).
$$
\end{corollary}

\begin{pf}
Since $(\xi^{\lor},h_{\xi^{\lor}})$ is Nakano semi-positive on $X$, the result follows from 
Theorem~\ref{Theorem1.4} (2) and \eqref{eqn:(6.20)}.
\end{pf}

\subsection
{The case of log-canonical singularities}
\label{sect:7.4}
\par
We assume that 
{\em $X_{0}$ is reduced and the pair $({\mathcal X},X_{0})$ has only log-canonical singularities}. 
As before, we set $l_{q}=h^{q}(X_{s},\Omega_{X_{s}}^{n}(\xi_{s}))$.

\begin{lemma}\label{Lemma7.5}
As $s\to0$
$$
\log\det\left(
{\Psi}_{\alpha}\otimes(\pi^{*}ds)^{-1}|_{X_{s}},{\Psi}_{\beta}\otimes(\pi^{*}ds)^{-1}|_{X_{s}}
\right)_{L^{2}}
=
O\left(\log(-\log|s|)\right).
$$
In particular, as $s\to0$,
$$
\log\|\sigma_{W}^{q}(t)\|_{L^{2}}^{2}=O\left(\log(-\log|s|)\right).
$$
\end{lemma}

\begin{pf}
Recall that $g_{\alpha\bar{\beta}}(s)$ was defined by \eqref{eqn:(7.9)}.
By \eqref{eqn:(7.11)} and the smoothness of the top form
$h_{\xi}(\psi_{\alpha}\wedge\overline{\psi}_{\beta})\wedge\kappa_{\mathcal X}^{q}$,
we get the following estimate by Lemma~\ref{Lemma7.1} (2)
\begin{equation}\label{eqn:(7.16)}
|g_{\alpha\bar{\beta}}(s)|\leq C\,(-\log|s|^{2})^{n}.
\end{equation}
The result follows from \eqref{eqn:(7.13)} and \eqref{eqn:(7.16)}.
\end{pf}

{\bf Proof of Theorem~\ref{Theorem1.4} (1) }
By Proposition~\ref{Theorem6.8}, Lemmas~\ref{Lemma6.9} and \ref{Lemma7.5}, we get 
${\Bbb L}_{g}\left(\mu^{*}R\pi_{*}\omega_{X/S}(\xi)/Rf_{*}\omega_{Y/T}(F^{*}\xi)\right)=0$.
The result follows from Theorem~\ref{Theorem1.2}.
\qed

\section
{Examples and questions}
\label{sect:8}
\par

\begin{example}\label{Example8.1}
Recall that  $\beta_{g}$ was defined in \eqref{eqn:(6.19)}.
We give examples of one-parameter families with $\alpha_{g}(X_{0},\omega_{X/S})\not=\beta_{g}$. 
Assume $G=\{1\}$ and write $\alpha$, $\beta$ for 
$\alpha_{\{1\}}(X_{0},\omega_{X/S})$, $\beta_{\{1\}}$, respectively.
We set $C:={\bf P}^{1}$. The inhomogeneous coordinate of $C$ is denoted by $s=s_{1}/s_{0}$.
For $d\in{\bf Z}_{>0}$, set
\begin{equation}\label{eqn:(8.1)}
{\mathcal X}':=\{([x],s)\in{\bf P}^{n+1}\times C;\, x_{0}^{d}+\cdots+x_{n}^{d}-s\,x_{n+1}^{d}=0\}
\end{equation}
and let $\mu\colon{\mathcal X}\to{\mathcal X}'$ be a projective resolution of the singularities 
of ${\mathcal X}'$. We set $\pi:={\rm pr}_{2}|_{{\mathcal X}'}\circ\mu$.
Then the family $\pi\colon{\mathcal X}\to C$ is smooth over $C\setminus\{0,\infty\}$.
Let $({\bf C}^{n+1},(z_{0},\cdots,z_{n}))$ be the inhomogeneous coordinates 
of ${\bf P}^{n+1}$ defined by $z_{i}=x_{i}/x_{n+1}$. 
Since ${\mathcal X}\cap({\bf C}^{n+1}\times{\bf C})$ is the graph of
the holomorphic function $z_{0}^{d}+\cdots+z_{n}^{d}$
from ${\bf C}^{n+1}$ to ${\bf C}$, we may regard ${\bf C}^{n+1}$
as a coordinate neighborhood of ${\mathcal X}$ such that $\pi(z)=z_{0}^{d}+\cdots+z_{n}^{d}$.
Then $X_{0}$ has a unique singular point at ${\frak p}:=((0:\cdots:0:1),0)$ and 
${\mathcal O}_{X_{0},{\frak p}}\cong{\bf C}\{z_{0},\ldots,z_{n}\}/(z_{0}^{d}+\cdots+z_{n}^{d})$.
We set $S:=\{s\in C;\,|s|<1\}$. Then $\Delta\cap S=\{0\}$. We compute $\alpha$ and $\beta$.
\par
By \eqref{eqn:(8.1)}, we regard $X_{s}\subset{\bf P}^{n+1}$ for $s\not=\infty$.
Since the ideal sheaf ${\mathcal I}_{X_{s}}$ of $X_{s}$ in ${\bf P}^{n+1}$ 
is isomorphic to ${\mathcal O}_{{\bf P}^{n+1}}(-d)$ and hence
$$
\omega_{X_{s}}
\cong
\omega_{{\bf P}^{n+1}}\otimes N_{X_{s}/{\bf P}^{n+1}}|_{X_{s}}
\cong
{\mathcal O}_{{\bf P}^{n+1}}(d-n-2)|_{X_{s}}
$$
by the adjunction formula, we get for all $s\in S$ and $q\geq0$
\begin{equation}\label{eqn:(8.2)}
H^{q}(X_{s},\omega_{X_{s}})
\cong
\begin{cases}
\begin{array}{ll}
H^{0}({\bf P}^{n+1},{\mathcal O}_{{\bf P}^{n+1}}(d-n-2))
&(q=0),
\\
0&(0<q<n),
\\
H^{n+1}({\bf P}^{n+1},\Omega_{{\bf P}^{n+1}})\cong{\bf C}
&
(q=n).
\end{array}
\end{cases}
\end{equation}
By \eqref{eqn:(8.2)}, the $L^{2}$-metric on $R^{q}\pi_{*}\omega_{{\mathcal X}/C}|_{S}$
is smooth on $S$ for $q>0$, i.e.,
\begin{equation}\label{eqn:(8.3)}
\log\|\sigma^{q}(s)\|_{L^{2}}=O(1)
\qquad
(q>0).
\end{equation}
By the relative Serre duality, 
$\pi_{*}\omega_{{\mathcal X}/C}\cong({\rm pr}_{2})_{*}\omega_{({\bf P}^{n+1}\times C)/C}({\mathcal X})$,
where the isomorphism is induced by the long exact sequence of 
direct images associated to the following short exact sequence of sheaves on ${\bf P}^{n+1}\times C$
$$
0\longrightarrow
\omega_{({\bf P}^{n+1}\times C)/C}
\longrightarrow
\omega_{({\bf P}^{n+1}\times C)/C}({\mathcal X})
\longrightarrow
\omega_{{\mathcal X}/C}
\longrightarrow 0
$$
and by the vanishing $R^{q}\pi_{*}\omega_{({\bf P}^{n+1}\times C)/C}=0$ for $q<n+1$.
\par 
Set $z^{e}:=z_{0}^{e_{0}}\cdots z_{n}^{e_{n}}$ and $|e|=\sum e_{i}$. 
A basis of 
$H^{0}({\bf P}^{n+1},\Omega_{{\bf P}^{n+1}}^{n+1}(X_{s}))
\cong H^{0}({\bf P}^{n+1},{\mathcal O}_{{\bf P}^{n+1}}(d-n-2))$
is given by
$\{z^{e}dz_{0}\wedge\cdots\wedge dz_{n}/(z_{0}^{d}+\cdots+z_{n}^{d}-s)\}_{|e|\leq d-(n+2)}$
and the corresponding basis of $\pi_{*}\omega_{{\mathcal X}/C}|_{S}$ 
as a free ${\mathcal O}_{S}$-module is given by
$$
\{\omega^{e}(s):=\frac{z^{e}}{z_{0}^{d-1}}dz_{1}\wedge\cdots\wedge dz_{n}|_{X_{s}}\}_{|e|\leq d-(n+2)}.
$$ 
We put 
$\sigma^{0}(s):=\bigwedge_{|e|\leq d-(n+2)}\omega^{e}(s)\in\det H^{0}(X_{s},\Omega_{X_{s}}^{n})$.
Then $\sigma^{0}$ gives a basis of the rank one ${\mathcal O}_{S}$-module $\det\pi_{*}\omega_{X/S}$. 
We define
$$
\psi_{s}\colon
{\bf P}^{n+1}\ni(z_{0}:\cdots:z_{n+1})
\to(s^{\frac{1}{d}}z_{0}:\cdots:s^{\frac{1}{d}}z_{n}:z_{n+1})
\in{\bf P}^{n+1}.
$$
Since $\psi_{s}\colon X_{1}\to X_{s}$ is an isomorphism with
$\psi_{s}^{*}\omega^{e}(s)=s^{\frac{n+|e|+1-d}{d}}\omega^{e}(1)$, 
we get
\begin{equation}\label{eqn:(8.4)}
\psi^{*}\sigma^{0}(s)
=
s^{\sum_{|e|\leq d-(n+2)}(n+|e|+1-d)/d}\sigma^{0}(1)
=
s^{-\binom{d}{n+2}/d}\sigma^{0}(1).
\end{equation}
Since 
$$
\begin{aligned}
\|\sigma^{0}(s)\|_{L^{2}(X_{s})}^{2}
&=
\det(
\int_{X_{s}}i^{n^{2}}\omega^{e}(s)\wedge\overline{\omega^{e'}(s)}
)_{|e|,|e'|\leq d-(n+2)}
\\
&=
\det(
\int_{X_{1}}i^{n^{2}}\psi_{s}^{*}\omega^{e}(s)\wedge
\overline{\psi_{s}^{*}\omega^{e'}(s)}
)_{|e|,|e'|\leq d-(n+2)}
=
\|\psi_{s}^{*}\sigma^{0}(s)\|_{L^{2}(X_{1})}^{2},
\end{aligned}
$$
we get by \eqref{eqn:(8.4)}
\begin{equation}\label{eqn:(8.5)}
\log\|\sigma^{0}(s)\|_{L^{2}}^{2}
=
-\frac{1}{d}\binom{d}{n+2}\log|s|^{2}
+
\log\|\sigma^{0}(1)\|_{L^{2}}^{2}.
\end{equation}
Substituting \eqref{eqn:(8.3)} and \eqref{eqn:(8.5)} into the first equality of \eqref{eqn:(6.17)}, we get
\begin{equation}\label{eqn:(8.6)}
\beta=\alpha+\frac{1}{d}\binom{d}{n+2}.
\end{equation}
By \cite[Th.\,8.1]{Yoshikawa07} and \eqref{eqn:(4.15)}, we get
\begin{equation}\label{eqn:(8.7)}
\alpha
=
(-1)^{n+1}\frac{(-1)^{n}}{(n+2)!}\mu(X_{0},{\frak p})
=
-\frac{(d-1)^{n+1}}{(n+2)!},
\end{equation}
where $\mu(X_{0},{\frak p})$ is the Milnor number of $(X_{0},{\frak p})$.
By \eqref{eqn:(8.6)}, \eqref{eqn:(8.7)}, we get
\begin{equation}\label{eqn:(8.8)}
\beta
=
\frac{1}{d}\binom{d}{n+2}-\frac{(d-1)^{n+1}}{(n+2)!}.
\end{equation}
In particular, if $d\geq n+2$, we get $\alpha\not=\beta$ in this example.
\par
Let us compute $\beta$ by using Lemma~\ref{Lemma6.9}. 
If $d\leq n+1$, then $\delta^{0}=0$ by \eqref{eqn:(8.2)}. Assume $d\geq n+2$.
Let $r\colon\widetilde{X}\to X$ be the blowing-up at $\frak o$ and set 
$\widetilde{\pi}:=r\circ\pi\colon\widetilde{X}\to S$. Then
$\widetilde{\pi}^{-1}(0)=\widetilde{X}_{0}+d\,E$, where $\widetilde{X}_{0}$ is the proper transform of $X_{0}$ 
and $E=r^{-1}({\frak o})$ is the exceptional divisor intersecting $\widetilde{X}_{0}$ transversally.
Define $\mu\colon T\to S$ by $\mu(t)=t^{d}$. Then the singular locus of $\widetilde{X}\times_{S}T$ 
is locally isomorphic to the product of the isolated two-dimensional singularity 
$\{(x,y,s)\in{\bf C}^{3};\,xy^{d}=s^{d}\}$ and ${\bf C}^{n-1}$. 
Let $p\colon Y\to \widetilde{X}\times_{S}T$ be the minimal resolution and set $f:={\rm pr}_{2}\circ p$. 
Then $f\colon(Y,Y_{0})\to(T,0)$ is a semistable reduction of $\pi\colon(X,X_{0})\to(S,0)$. 
We set $F:=r\circ{\rm pr}_{1}\circ p$.
Since $\|\omega^{e}(s)\|_{L^{2}}^{2}=|s|^{2(n+|e|+1-d)}\|\omega^{e}(1)\|_{L^{2}}^{2}$ 
by the relations $\psi_{s}^{*}\omega^{e}(s)=s^{\frac{n+|e|+1-d}{d}}\omega^{e}(1)$ and $s=t^{d}$, 
$t^{k}\mu^{*}\omega^{e}$ is an element of $\Gamma(T,f_{*}\omega_{Y/T})$
if and only if $k\geq d-(n+1)-|e|$ by Proposition~\ref{Proposition6.6}. Hence
\begin{equation}\label{eqn:(8.9)}
\mu^{*}(\pi_{*}\omega_{X/S})/f_{*}\omega_{Y/T}
\cong
\bigoplus_{|e|\leq d-(n+2)}
{\mathcal O}_{T}/{\frak m}_{0}^{d-(n+1)-|e|}.
\end{equation}
By \eqref{eqn:(8.9)},  $\delta^{0}=\binom{d}{n+2}$. 
Since $\delta^{q}=0$ for $q>0$ by \eqref{eqn:(8.2)}, we get \eqref{eqn:(8.6)} by Lemma~\ref{Lemma6.9}.
\par
Let $p\colon\widetilde{X}_{0}\to X_{0}$ be the blowing-up at ${\frak p}$ and set $E:=p^{-1}({\frak p})$. 
By \cite[(3.8.1)]{Kollar97}, $K_{\widetilde{X}_{0}}=p^{*}K_{X_{0}}+(n-d)E$. 
Since the isolated singularity $(X_{0},{\frak p})$ is log-canonical if and only if $d\leq n+1$,
Eq.\eqref{eqn:(8.6)} is compatible with Theorem~\ref{Theorem1.4}.
\end{example}

\begin{question}\label{Question8.2}
If all the direct images $R^{q}\pi_{*}\omega_{X/S}(\xi)$, $q\geq0$ are locally free on $S$, 
does Theorem~\ref{Theorem1.2} remain valid without assuming the Nakano semi-positivity of $\xi$?
In general, what can one say about the singularity of $\log\tau_{G}(X_{s},\omega_{X_{s}}(\xi_{s}))(g)$
without the assumption of the Nakano semi-positivity of $\xi$? 
Does the logarithmic divergence still hold? If it is the case, is the coefficient of $\log|s|^{2}$ expressed by
the topological term $\alpha_{g}(X_{0},\omega_{X/S}(\xi))$ and the $g$-action on the torsion modules
$\ker(\varphi)$ and ${\rm coker}(\varphi)$, where
$\varphi\colon R^{q}f_{*}\omega_{Y/T}(F^{*}\xi)\to\mu^{*}R^{q}\pi_{*}\omega_{X/S}(\xi)$ 
is the natural map?
\end{question}

\begin{question}\label{Question8.3}
If $\mu^{*}R^{q}\pi_{*}\omega_{X/S}(\xi)=R^{q}f_{*}\omega_{Y/T}(F^{*}\xi)$ for all Nakano semi-positive
vector bundles on $X$ and for all $q\geq0$, then is the pair $({\mathcal X},X_{0})$ log-canonical?
\end{question}

\begin{question}\label{Question8.4} 
Is the Lefschetz trace
${\Bbb L}_{g}(\mu^{*}R\pi_{*}\omega_{X/S}(\xi)/Rf_{*}\omega_{Y/T}(F^{*}\xi))$
expressed by some other geometric data like the discrepancies of $({\mathcal X},X_{0})$, 
the monodromy of the family $\pi\colon X\to S$ etc.?
\end{question}

\begin{question}\label{Question8.5}
Are the coefficients $\nu_{g},c_{g}$ in Theorem~\ref{Theorem1.2} expressible in terms of
some local or global geometric data associated to the family $\pi\colon(X,X_{0})\to(S,0)$?
\end{question}

\section
{Appendix}
\label{sect:9}
\par
In this section, we prove some technical results. We keep the notation in Sect.\,\ref{sect:6}.

\subsection
{A sequence of K\"ahler forms on $Z$ approximating $\kappa_{Z}$}
\label{subsect:9.2}
\par
The construction of a sequence of K\"ahler forms $\{\kappa_{Z,k}\}$ in the proof of 
Proposition~\ref{Proposition:Takegoshi} is as follows.
(Although the construction of such a sequence of K\"ahler forms is standard, we give it here 
for the completeness reason and its use in Sect.\,\ref{subsect:9.3}.)
Since $X\times_{T}S$ is a hypersurface of a complex manifold $X\times S$
and since $\rho\colon Z\to X\times_{S}T$ is a composite of blowing-ups with non-singular centers,
there is a sequence by e.g. \cite[Sect.\,1.2]{AKMW02}, \cite[Th.\,13.2]{BierstoneMilman97}
$$
\begin{matrix}
W &\to &W_{m-1} &\to\cdots\to &W_{1} &\to &X\times S
\\
\cup &\, &\cup &\, &\cup &\, &\cup
\\
Z &\to &Z_{m-1} &\to\cdots\to &Z_{1} &\to &X\times_{T}S
\end{matrix}
$$
where all $W_{i}$ are smooth projective algebraic manifolds, each morphism
$\pi_{i}\colon W_{i+1}\to W_{i}$ is given by the blowing-up with non-singular center $C_{i}\subset W_{i}$,
and $Z_{i+1}$ is the proper transform of $Z_{i}$ under the blowing-up $\pi_{i}\colon W_{i+1}\to W_{i}$.
We set $\rho_{i}:=\pi_{i}|_{Z_{i+1}}\colon Z_{i+1}\to Z_{i}$. Then 
$\rho=\rho_{0}\circ\cdots\circ\rho_{m-1}$.
Let $E_{i}$ be the exceptional divisor of $\pi_{i}\colon W_{i+1}\to W_{i}$. Then 
$E_{i}\cong{\bf P}(N_{C_{i}/W_{i}})$, where $N_{C_{i}/W_{i}}$ is the normal bundle of $C_{i}$ in $N_{i}$.
Let ${\mathcal O}_{W_{i}}(E_{i})$ be the line bundle on $W_{i}$ defined by the effective divisor $E_{i}$
and let $\sigma_{i}$ be the canonical section of ${\mathcal O}_{W_{i}}(E_{i})$ such that
$E_{i}={\rm div}(\sigma_{i})$.
\par
Since ${\rm pr}_{2}\circ\pi_{0}\circ\cdots\circ\pi_{i}(E_{i})=\{0\}$, there is a small open neighborhood
${\mathcal V}_{i}^{(k)}\subset W_{i+1}$ of $E_{i}$ such that
${\rm pr}_{2}\circ\pi_{0}\circ\cdots\circ\pi_{i}({\mathcal V}_{i}^{(k)})\subset T(\frac{1}{k})$.
Since ${\mathcal O}_{W_{i}}(E_{i})|_{E_{i}}\cong{\mathcal O}_{{\bf P}(N_{C_{i}/W_{i}})}(-1)$,
there is a Hermitian metric $h_{{\mathcal O}_{W_{i}}(E_{i})}$ on ${\mathcal O}_{W_{i}}(E_{i})$
such that $h_{{\mathcal O}_{W_{i}}(E_{i})}(\sigma_{i},\sigma_{i})=1$ on $W_{i}\setminus{\mathcal V}_{i}$
and such that 
$c_{1}({\mathcal O}_{W_{i}}(-E_{i}),h_{{\mathcal O}_{W_{i}}(E_{i})}^{-1})|_{E_{i}}=
dd^{c}\log h_{{\mathcal O}_{W_{i}}(E_{i})}(\sigma_{i},\sigma_{i})|_{E_{i}}$ 
is a positive $(1,1)$-form on the relative tangent bundle $T{\bf P}(N_{C_{i}/W_{i}})/C_{i}$.
For a K\"ahler form $\omega_{W_{i}}$ on $W_{i}$, there exists $A_{i}>0$ such that 
$A_{i}\pi_{i}^{*}\omega_{W_{i}}+c_{1}({\mathcal O}_{W_{i}}(-E_{i}),h_{{\mathcal O}_{W_{i}}(E_{i})}^{-1})$
is a positive $(1,1)$-form on $TW_{i}|_{E_{i}}$ and hence on a neighborhood of $E_{i}$.
Choosing $A_{i}$ large enough, we may assume that  
$A_{i}\pi_{i}^{*}\omega_{W_{i}}+c_{1}({\mathcal O}_{W_{i}}(-E_{i}),h_{{\mathcal O}_{W_{i}}(E_{i})}^{-1})$
is a positive $(1,1)$-form on ${\mathcal V}_{i}^{(k)}$. Since
$A_{i}\pi_{i}^{*}\omega_{W_{i}}+c_{1}({\mathcal O}_{W_{i}}(-E_{i}),h_{{\mathcal O}_{W_{i}}(E_{i})}^{-1})
=A_{i}\pi_{i}^{*}\omega_{W_{i}}$ on $W_{i}\setminus{\mathcal V}_{i}^{(k)}$, 
$\widetilde{\omega}_{W_{i}}^{(k)}:=
\pi_{i}^{*}\omega_{W_{i}}+\frac{1}{A_{i}}c_{1}({\mathcal O}_{W_{i}}(-E_{i}),h_{{\mathcal O}_{W_{i}}(E_{i})}^{-1})$
is a K\"ahler form on $W_{i}$ such that $\widetilde{\omega}_{W_{i}}^{(k)}=\pi_{i}^{*}\omega_{W_{i}}$
on $W_{i}\setminus{\mathcal V}_{i}^{(k)}$ and 
$[\widetilde{\omega}_{W_{i}}^{(k)}]=
\pi_{i}^{*}[\omega_{W_{i}}]+\frac{1}{A_{i}}c_{1}({\mathcal O}_{W_{i}}(-E_{i}))$.
(See e.g. \cite[Prop.\,12.4]{Demailly09} for more details.)
\par
We set $\omega_{W_{0}}:=\kappa_{X}+\kappa_{S}$. For $k>1$, we get a sequence of K\"ahler forms
$\{\omega_{W_{i}}^{(k)}\}$ by the procedure as above. We define $\kappa_{Z,k}:=\omega_{W}^{(k)}|_{Z}$.
Write $\rho^{-1}(\varSigma)=\bigcup_{\lambda\in\Lambda}D_{\lambda}$.
By the construction of $\omega_{W}^{(k)}$,  there is a Hermitian 
metric $h_{{\mathcal O}_{Z}(D_{\lambda}),k}$ on ${\mathcal O}_{Z}(D_{\lambda})$ 
and a real number $a_{\lambda}^{(k)}\geq0$ for every $\lambda\in\Lambda$ such that
\begin{equation}\label{eqn:(9.6)}
\kappa_{Z,k}
=
\kappa_{Z}
-
\sum_{\lambda\in\Lambda}
a_{\lambda}^{(k)}\,c_{1}({\mathcal O}_{Z}(D_{\lambda}),h_{{\mathcal O}_{Z}(D_{\lambda}),k})
\end{equation}
is a K\"ahler form on $Z$ with $\kappa_{Z,k}=\kappa_{Z}$ on
$Z\setminus\varpi^{-1}(T(\frac{1}{k}))$. Here $a_{\lambda}^{(k)}=0$ if $D_{\lambda}$
is not an exceptional divisor of $\rho\colon Z\to X\times_{T}S$ by construction.
By \eqref{eqn:(9.6)}, we get
\begin{equation}\label{eqn:(9.7)}
[\kappa_{Z,k}]|_{Z\setminus Z_{0}}=[\kappa_{Z}]|_{Z\setminus Z_{0}},
\end{equation}
since ${\mathcal O}_{Z}(D_{\lambda})|_{Z\setminus Z_{0}}\cong{\mathcal O}_{Z}|_{Z\setminus Z_{0}}$ 
for all $\lambda\in\Lambda$.
By \eqref{eqn:(9.2)}, \eqref{eqn:(9.7)}, we again get \eqref{eqn:(9.4)}.

\subsection
{Takegoshi's theorem with respect to the degenerate K\"ahler form}
\label{subsect:9.3}
\par
Since $\rho\colon Z\to X\times_{T}S$ is a resolution obtained as a composite of blowing-ups with 
non-singular centers, there is an extension of Takegoshi's theorem (Theorem~\ref{Theorem5.1}) for $Z$
by Mourougane-Takayama \cite[Prop.\,4.4]{MourouganeTakayama09} with respect to the degenerate 
K\"ahler form $\kappa_{Z}=\rho^{*}(\kappa_{X}+\kappa_{S})$. However, it is not immediate from 
their proof if the operator $L^{q}$ in \cite[Prop.\,4.4]{MourouganeTakayama09} is given by 
the multiplication by $\kappa_{Z}^{q}$. We make this point clear by proving the following:

\begin{proposition}\label{Proposition9.1}
For $u\in H^{q}(Z,\Omega_{Z}^{n+1}(\rho^{*}\xi))$, there exists
$\sigma\in H^{0}(Z,\Omega_{Z}^{n+1-q}(\rho^{*}\xi))$ such that
$u=[\sigma\wedge\kappa_{Z}^{q}]$ and $(\varpi^{*}dt)\wedge\sigma=0$.
\end{proposition}

\begin{pf}
We write $Z_{0}=\bigcup_{\lambda\in\Lambda}D_{\lambda}$.
By the construction in Sect.~\ref{subsect:9.2}, there is a K\"ahler form $\omega_{Z}$ on $Z$ of the form
$\omega_{Z}=\kappa_{Z}-\sum_{\lambda\in\Lambda}
a_{\lambda}\,c_{1}({\mathcal O}_{Z}(D_{\lambda}),h_{{\mathcal O}_{Z}(D_{\lambda})})$,
where $a_{\lambda}>0$ is a constant. By Theorem~\ref{Theorem5.1}, there exists 
$\sigma\in H^{0}(Z,\Omega_{Z}^{n+1-q}(\rho^{*}\xi))$ such that
$u=[\sigma\wedge\omega_{Z}^{q}]$ and $(\varpi^{*}ds)\wedge\sigma=0$.
Since $\sigma\wedge(\omega_{Z}^{q}-\kappa_{Z}^{q})=
\sigma\wedge(\omega_{Z}-\kappa_{Z})\wedge\sum_{i+j=q-1}\omega_{Z}^{i}\kappa_{Z}^{j}$,
it suffices to prove the equality of cohomology classes on $Z$
$$
[\sigma\wedge(\omega_{Z}-\kappa_{Z})]=0.
$$
Let $s_{\lambda}$ be the canonical section of ${\mathcal O}_{Z}(D_{\lambda})$ such that
${\rm div}(s_{\lambda})=D_{\lambda}$. Since
$$
\sigma\wedge(\omega_{Z}-\kappa_{Z})
=
\bar{\partial}\{
\sigma\wedge\sum_{\lambda\in\Lambda}
a_{\lambda}\,\partial\log h_{{\mathcal O}_{Z}(D_{\lambda})}(s_{\lambda},s_{\lambda})/2\pi i
\},
$$
it suffices to prove that
$\sigma\wedge\partial\log h_{{\mathcal O}_{Z}(D_{\lambda})}(s_{\lambda},s_{\lambda})$
is a $C^{\infty}$ differential form on $Z$.
\par
Set ${\frak D}:=\bigcup_{\lambda\not=\lambda'}D_{\lambda}\cap D_{\lambda'}$ and
$D_{\lambda}^{o}:=D_{\lambda}\setminus{\frak D}$.
Let $p\in D_{\lambda}^{o}$. There is a system of coordinates $({\mathcal V},(z_{0},\ldots,z_{n}))$ of $Z$ 
centered at $p$ such that $\varpi(z)=\epsilon(z)\,z_{0}^{e_{\lambda}}$ for some $e_{\lambda}\in{\bf Z}_{>0}$ 
and a nowhere vanishing $\epsilon(z)\in{\mathcal O}({\mathcal V})$. 
Since $\varpi(z)=\epsilon(z)\,z_{0}^{e_{\lambda}}$ and hence 
$\varpi^{*}ds=e_{\lambda}\epsilon(z)z_{0}^{e_{\lambda}-1}dz_{0}+z_{0}^{e_{\lambda}}d\epsilon$,
the condition $(\varpi^{*}ds)\wedge\sigma=0$ implies that 
\begin{equation}\label{eqn:(9.8)}
\frac{dz_{0}}{z_{0}}\wedge\sigma
=
-\frac{1}{e_{\lambda}\epsilon(z)}d\epsilon(z)\wedge\sigma
\in
\Omega_{\mathcal V}^{n+2-q}(\rho^{*}\xi).
\end{equation}
\par
Let $p\in{\frak D}$. 
There is a system of coordinates $({\mathcal U},(z_{0},\ldots,z_{n}))$ of $Z$ centered at $p$
such that $\varpi(z)=\epsilon(z)\,z_{0}^{e_{0}}\cdots z_{n}^{e_{n}}$ for some $e_{i}\in{\bf Z}_{>0}$ and
a nowhere vanishing $\epsilon(z)\in{\mathcal O}({\mathcal U})$. 
If $e_{i}>0$, then we deduce from \eqref{eqn:(9.8)} that
$(\frac{dz_{i}}{z_{i}}\wedge\sigma)|_{{\mathcal U}\setminus{\frak D}}\in 
\Omega_{{\mathcal U}\setminus{\frak D}}^{n+2-q}(\rho^{*}\xi)$.
Since ${\frak D}$ has codimension $2$ in ${\mathcal U}$, we get by the Hartogs extension theorem
\begin{equation}\label{eqn:(9.10)}
\frac{dz_{i}}{z_{i}}\wedge\sigma
\in 
\Omega_{\mathcal U}^{n+2-q}(\rho^{*}\xi).
\end{equation}
\par
Let $p\in Z_{0}$. There is a system of coordinates $({\mathcal U},(z_{0},\ldots,z_{n}))$ as above
such that $\varpi|_{\mathcal U}(z)=\epsilon(z)\,z_{0}^{e_{0}}\cdots z_{n}^{e_{n}}$. 
There exist $\phi\in C^{\infty}({\mathcal U})$ and $\nu_{i}\in{\bf R}$ such that
\begin{equation}\label{eqn:(9.11)}
\log h_{{\mathcal O}_{Z}(D_{\lambda})}(s_{\lambda},s_{\lambda})|_{\mathcal U}
=
\sum_{i=0}^{n}\nu_{i}\log|z_{i}|^{2}+\phi.
\end{equation}
Here $e_{i}=0$ implies $\nu_{i}=0$. By \eqref{eqn:(9.10)} and \eqref{eqn:(9.11)}, we get
\begin{equation}\label{eqn:(9.12)}
\sigma\wedge\partial\log h_{{\mathcal O}_{Z}(D_{\lambda})}(\sigma_{\lambda},\sigma_{\lambda})|_{U}
=
\sum_{i=0}^{n}\nu_{i}\,\sigma\wedge\frac{dz_{i}}{z_{i}}+\sigma\wedge\partial\phi
\in
A^{n+2-q,0}_{\mathcal U}(\rho^{*}\xi).
\end{equation}
Since $p\in Z_{0}$ is an arbitrary point, we get 
$\sigma\wedge\partial\log h_{{\mathcal O}_{Z}(D_{\lambda})}(\sigma_{\lambda},\sigma_{\lambda})
\in A^{n+2-q,0}_{Z}(\rho^{*}\xi)$.
This completes the proof.
\end{pf}

As a consequence of Proposition~\ref{Proposition9.2}, we get Takegoshi's theorem
(Theorem~\ref{Theorem5.1}) for $(Y,F^{*}\xi)$ with respect to the degenerate K\"ahler form 
$\kappa_{Y}=r^{*}(\kappa_{X}+\kappa_{T})$.

\begin{proposition}\label{Proposition9.2}
For $v\in H^{q}(Y,\Omega_{Y}^{n+1}(F^{*}\xi))$, there exists 
$\tau\in H^{0}(Y,\Omega_{Y}^{n+1-q}(F^{*}\xi))$ such that
$v=[\tau\wedge\kappa_{Y}^{q}]$ and $(f^{*}dt)\wedge\tau=0$.
\end{proposition}

\begin{pf}
Set $u:=q^{*}v\in H^{q}(Z,\Omega_{Z}^{n+1}(\rho^{*}\xi))$. By Proposition~\ref{Proposition9.1},
there exists $\sigma\in H^{0}(Z,\Omega_{Z}^{n+1-q}(\rho^{*}\xi))$ such that
$u=[\sigma\wedge\kappa_{Z}^{q}]$ and $(\varpi^{*}ds)\wedge\sigma=0$.
By the proof of Proposition~\ref{Proposition:Takegoshi}, 
there exists $\tau\in H^{0}(Y,\Omega_{Y}^{n+1-q}(F^{*}\xi))$ 
such that $q^{*}\tau=\sigma$. Since $q^{*}\kappa_{Y}=\kappa_{Z}$, we get
$q^{*}(v-[\tau\wedge\kappa_{Z}^{q}])=0$. Since the map of cohomologies
$q^{*}\colon H^{q}(Y,\Omega_{Y}^{n+1}(F^{*}\xi))\to H^{q}(Z,\Omega_{Z}^{n+1}(\rho^{*}\xi))$
is an isomorphism by \cite[Th.\,6.9 (i)]{Takegoshi95}, we get $v=[\tau\wedge\kappa_{Y}^{q}]$.
We get the equality $(\varpi^{*}dt)\wedge\sigma=0$ as before in the proof of 
Proposition~\ref{Proposition:Takegoshi}.
\end{pf}


\end{document}